\documentclass[12pt]{article} 
\usepackage{amssymb,latexsym}
\usepackage{mathtools}
\usepackage{graphicx,color,xcolor}

\usepackage[bookmarks]{hyperref}

\font\lagom=cmr10 at 10pt

\newcommand\ekv[2]{\begeq\label{#1}#2\endeq}
\newcommand\eekv[3]{\begin{eqnarray}\label{#1}#2 \\ #3
\nonumber\end{eqnarray}}

\newcommand\eeekv[4]{\begin{eqnarray}\label{#1}#2 \\ #3
\nonumber\\#4\nonumber\end{eqnarray}}

\newcommand\eeeekv[5]{\begin{eqnarray}\label{#1}#2 \\ #3
\nonumber\\#4\nonumber\\#5\nonumber\end{eqnarray}}

\newcommand\3{\vert \hskip -1pt\vert\hskip -1pt\vert }

\newcommand\aby{arbitrary} 
 
\newcommand\an{analytic}
\newcommand\asy{asymptotic} 
\newcommand\bdd{bounded} 
\newcommand\bdy{boundary}
 
\newcommand\coef{coefficient}

\newcommand\dop{differential operator}
\newcommand\ef{eigenfunction} 
\newcommand\ev{eigenvalue} 
\newcommand\e{equation}
\newcommand\fu{function} 
\newcommand\fy{family} 
\newcommand\F{Fourier} 
 
\newcommand\hol{holomorphic}
\newcommand\hs{Hilbert-Schmidt}

\newcommand\indep{independent}
\newcommand\inv{^{-1}} 
\newcommand\lhs{left hand side} 
\newcommand\mfld{manifold}
 
\newcommand\neigh{neighborhood}
\newcommand\no[1]{(\ref{#1})}

\renewcommand\3{\vert\hskip -1pt\vert\hskip -1pt\vert } 
\newcommand\on{orthonormal}
\newcommand\op{operator} 
\newcommand\og{orthogonal}
\newcommand\proba{probability}
\newcommand\pb{problem}

\newcommand\pert{perturbation}

\newcommand\pop{pseudodifferential operator}

\newcommand\rhs{right hand side}
\newcommand\rv{random variable} 
\newcommand\sa{selfadjoint}

\newcommand\sm{\setminus }

\newcommand\sufly{sufficiently}
\newcommand\tf{transformation}

\renewcommand\tf{transform}

\newcommand\trans{^t\hskip -2pt}
\newcommand\tr{{\rm tr\,}}

\newcommand\uf{uniform}
\newcommand\ufly{uniformly}

\newcommand\wrt{with respect to}

\renewcommand\Re{{\mathrm {Re\,}}}
\renewcommand\Im{{\mathrm {Im\,}}}

\newcommand\ran{\operatorname{ran}}

\newcommand{\al}{\sqrt{\alpha}}

\renewcommand\Im{\operatorname{Im}}  

\newcommand\wrtext[1]{\relax\ifmmode{\leavevmode\hbox{#1}}\else{#1}\fi}

\newcommand\begeq{\begin{equation}}

\def\endeq{\end{equation}}

\newcommand{\eps}{\epsilon}
\renewcommand{\phi}{\varphi}

\renewcommand\part[1]{\frac{\partial}{\partial #1}}

\renewcommand{\exp}{\mbox{\rm exp\,}}

\newtheorem{dref}{Definition}[section]
\newtheorem{lemma}[dref]{Lemma}
\newtheorem{theo}[dref]{Theorem}
\newtheorem{prop}[dref]{Proposition}
\newtheorem{remark}[dref]{Remark}
\newtheorem{ex}[dref]{Example}
\newtheorem{coro}[dref]{Corollary}

\newenvironment{proof}{\par\noindent{{\bf Proof:}}}{\hfill$\Box$\medskip}

\title{Eigenvalue asymptotics for randomly 
perturbed non-\sa{} \op{}s} 
\author{Mildred Hager\\ \lagom Dept of
Mathematics\\ \lagom  University of California\\
\lagom Berkeley, CA 94720\\ \lagom mhager@math.berkeley.edu
\and  Johannes
Sj{\"o}strand\\ \lagom CMLS\\ \lagom Ecole Polytechnique\\ \lagom FR  91120
Palaiseau c\'edex\\ \lagom johannes@{}math.polytechnique.fr \\ \lagom and 
UMR7640--CNRS}  
\date{}

\begin{document}
\maketitle

\begin{abstract} We consider quite general $h$-\pop{}s on ${\bf R}^n$ with small 
random perturbations and show that in the limit $h\to 0$ the \ev{}s are
distributed according to a Weyl law with a probabality that tends to 1. 
The first author has previously obtained a similar result in dimension 1. 
Our class of perturbations is different.\medskip
\par 
\centerline{\bf R\'esum\'e}\medskip
\par Nous consid\'erons des op\'erateurs $h$-pseudodiff\'erentiels assez 
g\'en\'eraux et nous montrons que dans la limite $h\to 0$, les valeurs 
propres se distribuent selon une loi de Weyl, avec une probabilit\'e qui 
tend vers 1. Le premier auteur a d\'ej\`a obtenu un r\'esultat semblable 
en dimension 1. Notre classe de perturbations est diff\'erente.
\end{abstract}

\tableofcontents
\section{Introduction}\label{int}
\setcounter{equation}{0}

\par This work can be viewed as a continuation of \cite{Ha2}, where one 
of us studied random perturbations of non-\sa{} $h$-\pop{}s on ${\bf R}$ 
and showed that Weyl \asy{}s holds with a \proba{} that is very close to 
1. In the present work we consider the multidimensional case and weaken 
some of the assumptions in \cite{Ha2} (like independence of the 
differentials  and analyticity of the symbol). Our random perturbations are 
slighly different however, in \cite{Ha2} they are given by a random 
potential while here they are rather given by a random integral \op{}. 

\par Before continuing the general discussion, we fix the framework more in 
detail.  We will work in the semi-classical limit on ${\bf R}^n$.  Write $\rho
=(x,\xi )$ and let $m\ge 1$ be an order \fu{} on the phase space ${\bf
R}^{2n}_{x,\xi }$:
\eekv{int.1}
{\exists C_0\ge 1,\ N_0>0\hbox{ such that }m(\rho )\le C_0\langle \rho 
-\mu \rangle ^{N_0}m(\mu ),}
{\forall \rho ,\mu \in{\bf R}^{2n},\ \langle \rho -\mu \rangle 
=\sqrt{1+\vert \rho -\mu \vert ^2}.}
The corresponding symbol space (cf \cite{DiSj}) is then 
\ekv{int.2}
{
S({\bf R}^{2n},m)=\{ a\in C^\infty ({\bf R}^{2n});\, \vert \partial _\rho 
^\alpha a(\rho )\vert \le C_\alpha m(\rho ),\, \rho \in {\bf R}^{2n},\, 
\alpha \in{\bf N}^{2n}\}.
}
Let 
\ekv{int.3}
{
P(\rho ;h)\sim p(\rho )+hp_1(\rho )+...\hbox{ in } S({\bf R}^{2n},m).
}
Assume $\exists$ $z_0\in{\bf C},\, C_0>0$ such that 
\ekv{int.4}
{
\vert p(\rho )-z_0\vert \ge m(\rho )/C_0,\ \rho \in{\bf R}^{2n}.
}
Let $\Sigma $ denote the closure of $p({\bf R}^{2n})$ so that $\Sigma 
=p({\bf R}^{2n})\cup \Sigma _\infty $, where $\Sigma _\infty \subset {\bf 
C}$ is the set of accumulation points of $p$ in the limit $(x,\xi 
)=\infty $. 

\par For $h>0$ small enough, we also let $P$ denote the $h$-Weyl 
quantization,
$$
Pu(x)=P^w(x,hD_x;h)u(x)={1\over (2\pi h)^n}\iint e^{{i\over h}(x-y)\cdot 
' }P({x+y\over 2},\eta ;h)u(y)dyd\eta .
$$

\par Let $\Omega \Subset  {\bf C}\sm \Sigma _\infty $ be open simply
connected
and not entirely contained in $\Sigma $. Then, as we shall see, 
\smallskip
\par\noindent $1^o$ $\sigma (P)\cap \Omega $ is discrete for $h>0$ small enough,
\smallskip
\par\noindent $2^o$ $\forall\,\epsilon >0$, $\exists\, h(\epsilon )>0$, such that
$$
\sigma (P)\cap \Omega \subset \Sigma +D(0,\epsilon ),\ 0<h\le h(\epsilon ).
$$
Here $D(0,\epsilon )$
denotes the open disc in ${\bf C}$ with center $0$ and radius $\epsilon 
>0$ and we equip the \op{} $P$ with the domain $H(m):= 
(P-z_0)^{-1}(L^2({\bf R}^n))$, where the \op{} to the right is the 
pseudodifferential inverse of $P-z_0$ (see \cite{DiSj} and \cite{Ha2}).
\par If $P$ is \sa{} (so that $p$ is real-valued) we have Weyl \asy{}s:

\par For every interval $I\subset \Omega $ with ${\rm vol}_{{\bf 
R}^{2n}}(p\inv (\partial I))=0$, the number $N(P,I)$ of \ev{}s of $P$ in 
$I$ satisfies
\ekv{int.5}
{ N(P,I)={1\over (2\pi h)^n}({\rm vol\,}(p\inv (I))+o(1)),\ h\to 0.  } This
result has been proved with increasing generality and precision by
J.~Chazarain, B.~Helffer--D.~Robert, and V.~Ivrii.  (We here
follow the presentation of
\cite{DiSj} where references to original works can be found. The corresponding
developement for \sa{} partial differential \op{}s in the high energy limit
has a long and rich history starting with the work of H.~Weyl \cite{We}.) 
A very simple and explicit example is given by the harmonic oscillator
$P={1\over 2}((hD)^2+x^2):L^2({\bf R})\to L^2({\bf R})$, $P(x,\xi )=p(x,\xi
)={1\over 2}(x^2+\xi ^2)$.  In this case the \ev{}s are given by $\lambda
_k(h)=(k+{1\over 2})h$, $k=0,1,2,...$

\par In the non-\sa{} case, Weyl \asy{}s does not always hold. If $P$ is a 
\dop{} with \an{} \coef{}s on the real line, then often the spectrum is 
determined by action integrals over complex cycles, having nothing to do 
with volumes of subsets of real phase space. A simple example of this is 
given by the non-\sa{} harmonic oscillator,
\ekv{int.6}
{
P={1\over 2}((hD)^2+ix^2):L^2({\bf R})\to L^2({\bf R}),
} 
whose spectrum is equal to $\{ e^{i\pi /4}(k+{1\over 2})h;\, k\in{\bf 
N}\}$; This is easy to see by the method of complex scaling, or by 
applying the general multidimensional result of \cite{Sj}. In this case, we 
have $\Sigma _\infty =\emptyset $, and $\Sigma $ is the closed 1st 
quadrant. Clearly the number of \ev{}s in an open set 
$\Gamma \Subset  {\bf C}$ intersecting the 1st quadrant, whose closure 
avoids the ray given $\arg z={\pi \over 4}$ is equal to zero while the 
corresponding 
Weyl \coef{} ${\rm vol}(p^{-1}(\Gamma ))$ is not. (Further results about 
the non-\sa{} harmonic oscillator have been obtained by E.B. Davies and L. 
Boulton, see \cite{Da} and further references given there). 

\par However, in this case and for quite a general class of $h$-\pop{}s in 
one dimension, it was shown by one of us in \cite{Ha2} that if we replace 
the \op{} $P$ by $P+\delta Q_\omega $, where $0<\delta \ll 1$ 
varies in a 
suitable parameter range and $Q_\omega $ is a random potential of a 
suitable type then we do 
have Weyl \asy{} in the interior of $\Sigma $ with a \proba{} that is 
close to 1. The book \cite{EmTr} of M. Embree and L.N. Trefethen as well as 
the paper \cite{TrCh} by L.N. Trefethen and S.J. Chapman contain (in our 
opinion) numerical examples where one can see the onset  of Weyl-\asy{}s 
after adding small random \pert{}s. 

\par In this work we establish similar results in arbitrary dimension that we 
shall now describe. Let $0<\widetilde{m},\widehat{m}\le 1$ be 
square integrable order \fu{}s on ${\bf R}^{2n}$
such that $\widetilde{m}$ or $\widehat{m}$ is integrable, and let 
$\widetilde{S}\in S(\widetilde{m})$, $\widehat{S}\in S(\widehat{m})$ be 
elliptic symbols. We use the same symbols to denote the $h$-Weyl 
quantizations. The \op{}s $\widetilde{S}$, $\widehat{S}$ are then \hs{} 
with 
$$
\Vert \widetilde{S}\Vert _{{\rm HS}}, \Vert \widehat{S}\Vert _{{\rm HS}}\backsim 
h^{-{n\over 2}},
$$
where $\backsim$ indicates same order of magnitude.
Let $\widetilde{e}_1,\widetilde{e}_2,...$, and 
$\widehat{e}_1,\widehat{e}_2,...$ be \on{} bases for $L^2({\bf R}^n)$. Our 
random perturbation will be 
\ekv{int.7}
{ Q_\omega =\widehat{S}\circ \sum_{j,k}\alpha _{j,k}(\omega
)\widehat{e}_j\widetilde{e}_k^*\circ \widetilde{S}, } where $\alpha _{j,k}$
are \indep{} complex ${\cal N}(0,1)$ \rv{}s, and 
$\widehat{e_j}\widetilde{e}_k^*u=(u\vert \widetilde{e}_k)\widehat{e}_j$, 
$u\in L^2$.  In the appendix Section 
\ref{ap} we show that up to a change of the set of \indep{} ${\cal
N}(0,1)$-laws, the representation \no{gp.1} is \indep{} of the choice of 
bases $\widehat{e}_j$ and $\widetilde{e}_j$.

 Let 
\ekv{int.8}
{
M=C_1h^{-n},
}
for some $C_1\gg 1$. Then, as we shall see in Section \ref{gp}, we have 
the following estimate on the \proba{} that $Q$ be large in Hilbert-Schmidt norm:
\ekv{int.9}
{
P(\Vert Q\Vert _{{\rm HS}}^2\ge M^2)\le C\exp (-h^{-2n}/C),
}
for some new constant $C>0$. In the following discussion we may restrict 
the attention to the case when $\Vert Q_\omega \Vert _{{\rm HS}}\le M$. We 
wish to study the \ev{} distribution of $P+\delta Q_\omega $ for $\delta $ 
in a suitable range.

\par Let $\Gamma \Subset \Omega  $ be open with $C^2$ \bdy{} and assume that for every $z\in\partial 
\Gamma $:
\eekv{int.10}
{&&\Sigma _z:=p\inv (z)\hbox{ is a smooth sub-\mfld{} of }T^*{\bf R}^n\hbox{ 
on}}
{&&\hbox{which }dp,d\overline{p}\hbox{ are linearly \indep{} 
at every point.}} The following result will be established in Section 
\ref{sa}.

\begin{theo}\label{int1}
Let $\Gamma \Subset \Omega $ be open with $C^2$ \bdy{} and make the
assumption \no{int.10}.  For $0<h\ll 1$, let $\delta >0$ be a small parameter such that  
$$
0<\delta\ll h^{3n+1/2}.
$$
For some small parameter $0<\eps\ll 1$ assume $h\ln \delta \inv \ll \epsilon \ll 1$ (or
equivalently $\delta \ge e^{-\epsilon /(Dh)}$ for some $D\gg 1$,
implying also that $\epsilon \ge \widetilde{C}h\ln h\inv$ for some
$\widetilde{C}>0$).  Then there is a constant $C>0$ (that is independent of $h$, $\delta$ and $\eps$) such that the number
$N(P_\delta ,\Gamma )$ of \ev{}s of $P_\delta $ in $\Gamma $ satisfies
\ekv{int.11}
{
\vert N(P_\delta ,\Gamma )-{1\over (2\pi h)^n}{\rm vol\,}(p^{-1}(\Gamma 
))\vert \le C{\sqrt{\epsilon }\over h^n}
}
with \proba{} 
$$\ge 1 -{C\over \sqrt{\epsilon }}e^{-{\epsilon/2 \over (2\pi 
h)^n}}.
$$
\end{theo}

This is a restatement of Theorem \ref{sa1}. In Theorem \ref{sa2} we give a 
similar statement about the  simultaneous Weyl \asy{}s for all $\Gamma $s 
in a \fy{} of sets that satisfy the assumptions of the above 
theorem \ufly{}. The lower bound on the \proba{} becomes slightly worse 
but is still very close to 1 for suitable values of $\epsilon $.

\par The condition \no{int.10} says that $\partial \Gamma $ does not 
intersect the set of critical values of $p$ and this is clearly not a 
serious restriction when $\overline{\Gamma }$ is contained in the interior 
of $\Sigma $. However, we also would like to study the \ev{} distribution 
near the \bdy{} of $\Sigma $, and we then need a weaker assumption. 

\par Let $\Gamma \Subset  \Omega $ be open with 
$C^\infty $ \bdy{}. For $z$ in a \neigh{} of $\partial \Gamma $ and 
$0<s,t\ll 1$, we put
\ekv{int.12}
{
V_z(t)={\rm Vol\,}\{\rho \in{\bf R}^{2n};\, \vert p(\rho )-z\vert ^2\le 
t\}.
}
Our weak assumption, replacing \no{int.10} is
\ekv{int.13}
{
\exists \kappa\in ]0,1], \hbox{ such that }V_z(t)={\cal O}(t^{\kappa}),\hbox{ \ufly{} for }z\in{\rm neigh\,}(\partial \Gamma ),\ 0\le t\ll 1.
}
Here we have written in an informal way ``${\rm neigh\,}(\partial \Gamma)$'' for 
some neighbourhood of $\partial \Gamma$. 
Notice that \no{int.10} implies \no{int.13} with $\kappa=1$.

\par Generically, we will have $\{ p,\{ p,\overline{p}\}\}\ne 0$ when 
$p(\rho )\in\partial \Sigma $ and if we assume that
\ekv{int.14}
{\hbox{at every point of }p\inv (0),\hbox{ we have }
\{ p,\overline{p}\}\ne 0 \hbox{ or } \{ p,\{ p,\overline{p}\}\} \ne 0,} 
then as shown in Example \ref{gc0}, we have \no{int.13} with $\delta 
_0=3/4$. (When $p$ is analytic, we believe that \no{int.13} will always 
hold with some $\kappa>0$ but have not consulted with experts in 
analytic geometry.) Under this more general assumption, we have  
\begin{theo}\label{int2} Assume \no{int.13} and let $\delta $ satisfy 
$$ 0<\delta \ll h^{3n+1/2}.  $$ Let
$N(P+\delta Q_\omega ,\Gamma )$ be the number of \ev{}s of $P+\delta
Q_\omega $ in $\Gamma $.  Then for every fixed $K>0$ and for $0<r\ll 1$:
\eekv{int.15}
{
&&\vert N(P+\delta Q_\omega ,\Gamma )-{1\over (2\pi h)^n}\iint_{p\inv 
(\Gamma )}dxd\xi \vert \le} {&&{C\over h^n}\Big({\epsilon \over r}+C_K\big(r^K+\ln 
({1\over r})\iint_{p\inv (\partial \Gamma +D(0,r)}  dxd\xi \big) \Big),\ 0<r\ll 1,
}
with \proba{} 
\ekv{int.16}
{
\ge 1-{C\over r}e^{-{\epsilon\over 2}(2\pi 
h)^{-n}}} provided that 
\ekv{int.17}{
h^{\kappa}\ln{1\over \delta }
\ll \epsilon \ll 1,
}
or equivalently, 
$$
e^{-{\epsilon \over Ch^{\kappa}}}\le \delta ,\ C\gg 1,\ \epsilon \ll 1,
$$
implying that $\epsilon \ge \widetilde{C}h^{\kappa}\ln{1\over h}$, for 
some $\widetilde{C}>0$.

\end{theo}

\par This is a restatement of Theorem \ref{gc2} and as explained after that
theorem, when
$\kappa>1/2$, the integral in the \rhs{} of \no{int.15} is ${\cal 
O}(r^{2\kappa-1})$ and it follows that we have Weyl \asy{}s with \proba{} 
close to 1, if we let $r$ be a suitable power of $\epsilon $. To have the same conclusion when $\kappa\le 1/2$ we can 
assume that the integral is ${\cal O}(r^{\alpha _0})$ for some $\alpha _0>0$.

\par Again we have a similar theorem about the simultaneous \asy{}s for 
$N(P+\delta Q_\omega ,\Gamma )$ when $\Gamma $ varies in a bounded family 
of domains satisfying all the assumptions \ufly{}. See Theorem \ref{gc3}.

\par The proofs follow the same general strategy as those of \cite{Ha2} 
with some essential differences:
\smallskip
\par\noindent We do not use any non-vanishing assumption on the Poisson 
bracket ${1\over i}\{ p,\overline{p}\}$. 
Instead we work systematically with the \op{}s $P^*P$ and $PP^*$ and their 
\ef{}s in order to set up a Grushin-\pb{}. 
\smallskip
\par\noindent As in \cite{Ha2} we reduce ourselves to the study of a 
random \hol{} \fu{}, but in the present work this \fu{} appears as the 
determinant of the full \op{} (essentially) and we need to make some 
estimates for determinants of random matrices, and especially to prove 
that such a determinant is not too small with a \proba{} close to 1. Those 
estimates were \sufly{} elementary to be carried out by hand, but we think 
that future generalizations and improvements will require a careful study 
of the existing results on random determinants and possibly the 
derivation of new results in that direction. See the book \cite{Gi} of 
V.I.~Girko.

\medskip
\par \noindent
{\it Acknowledgements.} We are grateful to F. Klopp for helping us to find
some references.  The first author was supported by a postdoctorial
fellowship from Ecole Polytechnique.  She also thanks Y.~Colin de 
Verdi\`ere for an interesting discussion around random functions.  The second
author is grateful to the Japan Society for the promotion of Sciences and
to the Dept of Mathematics of Tokyo Unversity for offering excellent
working conditions during the month of July, 2005.  He also thanks
E.~Servat for a very interesting remark. We also thank the referee for 
many detailed remarks that have helped to improve the presentation.

\section{Determinants and Grushin \pb{}s}\label{dg}
\setcounter{equation}{0}

Here we mainly follow \cite {SjZw} and give a more explicit formulation of 
one of the results there. Let ${\cal H}_1$, ${\cal H}_2$, ${\cal G}_1$, ${\cal 
G}_2$ be complex Hilbert spaces and let $A_{j,k}:\,{\cal H}_k\to {\cal G}_j$ 
be \bdd{} \op{}s depending in a $C^1$ fashion on the real parameter $t\in 
]a,b[$. We also assume that $\dot{A}_{j,k}$ are of trace class and 
continuously dependent of $t$ in the space of such \op{}s. Here ``over-dot'' 
means derivative with respect to $t$.

\begin{prop}\label{dg1} (\cite{SjZw}) Assume in addition that 
$A=(A_{j,k}):{\cal H}_1\times {\cal H}_2\to {\cal G}_1\times {\cal G}_2$ 
has a \bdd{} inverse $B=(B_{j,k})$, and that $A_{2,2}$ and $B_{1,1}$
are invertible. (The invertibility of one of $A_{2,2}$, $B_{1,1}$ implies 
that of the other.)
Then 
\ekv{dg.1}
{{\rm tr\,}\dot{A}B={\rm tr\,}\dot{A}_{2,2}A_{2,2}^{-1}-{\rm 
tr\,}B_{1,1}^{-1}\dot{B}_{1,1}.}
\end{prop}
\begin{proof}
We expand
$$
\dot{B}_{j,k}=-\sum_\nu \sum_\mu  B_{j,\nu } \dot{A}_{\nu ,\mu }B_{\mu ,k},
$$
that are of the trace class too.
In particular,
$$
-\dot{B}_{1,1}=B_{1,1}\dot{A}_{1,1}B_{1,1}+B_{1,1}\dot{A}_{1,2}B_{2,1}+ 
B_{1,2}\dot{A}_{2,1}B_{1,1}+B_{1,2}\dot{A}_{2,2}B_{2,1}.
$$
Rewrite the \rhs{} of \no{dg.1}:
\begin{eqnarray*}
&&{\rm tr\,}\dot{A}_{2,2}A_{2,2}^{-1}-{\rm tr\,}B_{1,1}^{-1}\dot{B}_{1,1}\\
&=&\tr \dot{A}_{2,2}A_{2,2}^{-1}+\tr \dot{A}_{1,1}B_{1,1}+\tr 
\dot{A}_{1,2}B_{2,1} +\tr B_{1,1}^{-1}B_{1,2} \dot{A}_{2,1}B_{1,1}+
\tr B_{1,1}^{-1}B_{1,2}\dot {A}_{2,2}B_{2,1}\\
&=&\tr \dot{A}_{2,2}A_{2,2}^{-1}+\tr \dot{A}_{1,1}B_{1,1}+\tr 
\dot{A}_{1,2}B_{2,1} +\tr \dot{A}_{2,1}B_{1,2} +
\tr B_{1,1}^{-1}B_{1,2}\dot{A}_{2,2}B_{2,1}\\
&=&\tr \dot{A}_{2,2}(A_{2,2}^{-1}+B_{2,1}B_{1,1}^{-1}B_{1,2})+\tr 
\dot{A}_{1,1}B_{1,1}+\tr \dot{A}_{1,2}B_{2,1}+\tr \dot{A}_{2,1}B_{1,2}\\
&=&\tr \dot{A}B.
\end{eqnarray*}
Here we used the cyclicity of the trace and for the last equality the fact 
that 
\ekv{dg.2}
{
B_{2,2}=A_{2,2}^{-1}+B_{2,1}B_{1,1}^{-1}B_{1,2}.
}
To check \no{dg.2} we proceed by equivalences:
\begin{align*}
\hbox{\no{dg.2}}
&\iff 
A_{2,2}^{-1}=B_{2,2}-B_{2,1}B_{1,1}^{-1}B_{1,2}\\
&\iff
1 = A_{2,2}B_{2,2}-A_{2,2}B_{2,1}B_{1,1}^{-1}B_{1,2}\\
&\iff 
1=(1-A_{2,1}B_{1,2})+A_{2,1}B_{1,1}B_{1,1}^{-1}B_{1,2}\\
&\iff 
1= 1.
\end{align*}
Here and in the following, we often denote the identity operator by $1$ when the 
meaning is clear from the context.
\end{proof}

\par Consider the case ${\cal H}_1={\cal G}_1={\cal H}$, ${\cal 
H}_2={\cal G}_2={\bf C}^N$,
$$
A={\cal P}=\begin{pmatrix}P &R_-\cr R_+ &0\end{pmatrix}.
$$
Assume also that $P$, ${\cal P}$ are invertible. 
(In the proposition we can permute the 
indices $1$ and $2$ and think of $P$ as $A_{2,2}$.) We look for 
$$
\widetilde{{\cal P}}=\begin{pmatrix}P &\widetilde{R}_-\cr \widetilde{R}_+ 
&\widetilde{R}_{+-}\end{pmatrix}:{\cal H}\times {\bf C}^N\to {\cal H}\times {\bf C}^N,
$$
that is also invertible, i.e.\ we should be able to solve uniquely
\ekv{dg.3}
{
\begin{cases}Pu+\widetilde{R}_-\widetilde{u}_-=v,\cr 
\widetilde{R}_+u+\widetilde{R}_{+-}\widetilde{u}_-=\widetilde{v}_+.
\end{cases}
}
Let 
$$
\begin{pmatrix}E &E_+\cr E_- &E_{-+}\end{pmatrix}=\begin{pmatrix}P &R_-\cr R_+ &0\end{pmatrix}^{-1}.
$$
Rewrite the first \e{} in \no{dg.3} as 
$Pu=v-\widetilde{R}_-\widetilde{u}_-$. The general solution to that \e{} is 
$$
u=E(v-\widetilde{R}_-\widetilde{u}_-)+E_+v_+,
$$
where $v_+$, $\widetilde{u}_-$ should be determined so that 
\ekv{dg.4}
{
0=E_-(v-\widetilde{R}_-\widetilde{u}_-)+E_{-+}v_+.
}

\par The second \e{} in \no{dg.3} becomes
\ekv{dg.5}
{
\widetilde{v}_+=\widetilde{R}_+E(v-\widetilde{R}_-\widetilde{u}_-)+
\widetilde{R}_+E_+v_++\widetilde{R}_{+-}\widetilde{u}_-.
}
Hence we get the following system that is equivalent to \no{dg.3}:
\ekv{dg.6}
{
\begin{cases}
E_{-+}v_+-E_-\widetilde{R}_-\widetilde{u}_-=-E_-v,\cr 
\widetilde{R}_+E_+v_++(\widetilde{R}_{+-}-\widetilde{R}_+E\widetilde{R}_-)\widetilde{u}_-
=\widetilde{v}_+-\widetilde{R}_+Ev,
\end{cases}
}
so \no{dg.3} is well-posed iff
\ekv{dg.7}
{
\begin{pmatrix}E_{-+}&-E_-\widetilde{R}_-\cr \widetilde{R}_+E_+ 
&\widetilde{R}_{+-}-\widetilde{R}_+E\widetilde{R}_-\end{pmatrix}:{\bf C}^{2N}\to{\bf 
C}^{2N} 
}
is invertible.

\par Choose $\widetilde{R}_+=tR_+$, $\widetilde{R}_-=sR_-$, 
$\widetilde{R}_{+-}=r{\rm id}_{{\bf C}^N}$ with $s,t,r\in{\bf C}$, and use 
that $R_+E_+=1$, $E_-R_-=1$, $R_+E=0$, to see that the matrix \no{dg.7} is 
equal to
\ekv{dg.8}
{
\begin{pmatrix}E_{-+} &-s\cr t & r\end{pmatrix}.
}
This matrix is invertible precisely when $(s,t,r)$ belongs to the set
\ekv{dg.9}
{
\{ (s,t,0);\, st\ne 0\}\cup \{ (s,t,r);\, r\ne 0,\, -{st\over r}\not\in 
\sigma (E_{-+})\} .
}

\par Since $P$ is invertible, we know that $0\not\in \sigma (E_{-+})$. We 
can therefore find a $C^1$-curve
$$
[0,1]\ni \tau \mapsto (s(\tau ),t(\tau ),r(\tau ))\in\,\hbox{the set 
\no{dg.9}},
$$
with 
$$(s(0),t(0),r(0))=(1,1,0),\quad (s(1),t(1),r(1))=(0,0,1).$$
This means that we have a $C^1$ deformation
$$
{\cal P}(\tau )=\begin{pmatrix}P &s(\tau )R_-\cr t(\tau )R_+ &r(\tau )1\end{pmatrix}:\ {\cal 
H}\times {\bf C}^N\to {\cal H}\times {\bf C}^N
$$
 of bijective \op{}s with 
$${\cal P}(0)={\cal P},\ {\cal P}(1)=\begin{pmatrix}P &0\cr 0 &1\end{pmatrix}.$$  
Applying \no{dg.1} with the indices ``1'' and ``2'' permuted gives
$$
\tr \dot{{\cal P}}{\cal P}^{-1}=\tr \dot{P}P^{-1}-\tr 
E_{-+}^{-1}(\tau )\dot{E}_{-+}(\tau )=-\tr 
E_{-+}^{-1}(\tau )\dot{E}_{-+}(\tau ),
$$
where now ``over-dot'' means derivative with respect to $\tau$.
If we integrate from $\tau =0$ to $\tau =1$, we get with a suitable 
choice of branches for $\ln$:
$$
\ln \det P-\ln \det {\cal P}=\ln \det E_{-+}(0).
$$
For this relation to make sense we also assume that 
\ekv{dg.9.5}
{
P-1\hbox{ is of trace class}.
}

\par Then for the original \op{} and its inverse we have 
\ekv{dg.10}
{
\ln \det P=\ln \det {\cal P}+\ln \det E_{-+},
}
or equivalently,
\ekv{dg.11}
{
\det P=\det {\cal P}\det E_{-+}.
}
\section{General frame-work and reduction to trace class \op{}s}\label{gf}
\setcounter{equation}{0}

Let $m\ge 1$ be an order \fu{} on ${\bf R}^{2n}$ in the sense that there 
exist constants $C_0\ge 1$, $N_0>0$, such that
$$
m(\rho )\le C_0\langle \rho -\mu \rangle ^{N_0}m(\mu ),\ \forall \rho ,\mu 
\in{\bf R}^{2n},
$$
where we write $\langle \rho \rangle =\sqrt{1+\vert \rho \vert ^2}$. We 
consider a symbol
$$
P(\rho ;h)\sim p(\rho )+hp_1(\rho )+...\hbox{ in }S({\bf R}^{2n},m),
$$
where 
$$
S({\bf R}^{2n},m)=\{ a\in C^\infty ({\bf R}^{2n});\, \vert \partial _{x,\xi 
}^\alpha a(x,\xi )\vert \le C_\alpha m(x,\xi ),\, \forall (x,\xi )\in{\bf 
R}^{2n},\, \alpha \in{\bf N}^{2n}\} .
$$
Put 
$$
\Sigma =\overline{p({\bf R}^{2n})},\ \Sigma _\infty= \{ z\in{\bf 
C};\,\exists \rho _j\in{\bf R}^{2n},\ j=1,2,3,...,\, \rho _j\to \infty ,\, 
p(\rho _j)\to z,\, j\to \infty \}.
$$
Assume $\exists z_0\in{\bf C}\setminus \Sigma ,\, C_0>0$, such that
\ekv{gf.1}
{
\vert p(\rho )-z_0\vert \ge m(\rho )/C_0,\ \forall \rho \in{\bf R}^{2n}.
}
Then as pointed out in \cite{Ha2}, for every $z\in {\bf C}\setminus \Sigma 
$, there exists $C>0$ such that 
\ekv{gf.2}
{
\vert p(\rho )-z\vert \ge m(\rho )/C,\ \forall \rho \in{\bf R}^{2n},
}
and for every $z\in {\bf C}\setminus \Sigma _\infty 
$, there exists $C>0$ such that 
\ekv{gf.3}
{
\vert p(\rho )-z\vert \ge m(\rho )/C,\ \forall \rho \in{\bf R}^{2n}\hbox{ 
with }\vert \rho \vert \ge C.
}

\par Let $\Omega \Subset  {\bf C}\setminus \Sigma _\infty $ be open 
simply connected containing at least one point $z_0\in{\bf C}\setminus 
\Sigma $.

\begin{lemma}\label{gf1}
For every compact set $K\subset \Omega $, there exists a smooth map 
$\kappa :\Omega \setminus\{ z_0\}\to \Omega \setminus\{ z_0\}$ with 
$\kappa (z)=z$ for all $z$ in a \neigh{} of $\partial \Omega $, such that 
$\kappa (\Sigma \cap \Omega )\cap K=\emptyset$.
\end{lemma}
\begin{proof}
$\Omega $ is diffeomorphic to the open unit disc $D(0,1)$ in such a way 
that $z_0$ corresponds to $0$. Now consider $\widetilde{\kappa 
}:D(0,1)\setminus\{ 0\}\to D(0,1)\setminus\{ 0\}$ defined by 
$\widetilde{\kappa }(z)=f(\vert z\vert )z/\vert z \vert $, where $f$ is 
a smooth \fu{} on $]0,1]$ with $1-\epsilon \le f(r)\le 1$, such that 
$f(r)=r$ for $1-r\le \epsilon /2$. Choosing $\epsilon >0$ small enough and 
conjugating with the diffeomorphism above we 
get the desired map $\kappa $.
\end{proof}

\par Let $\widetilde{\Omega }\Subset  \Omega $ be open. 
Take $\kappa $ as in the lemma with $K$ containing the closure of 
$\widetilde{\Omega }$. Extend $\kappa $ to be the identity in ${\bf
C}\setminus\Omega $ and put $\widetilde{p}=\kappa \circ p$. Then  
$\widetilde{p}(\rho )-z$ is elliptic in the sense of \no{gf.2}, \ufly{} 
for $z\in \widetilde{\Omega }$ and 
\ekv{gf.4}
{
\widetilde{p}-p\in C_0^\infty ({\bf R}^{2n}).
}
Put
$$
\widetilde{P}=\widetilde{p}+hp_1+h^2p_2+ ...\ .
$$

\par Now pass to \op{}s and denote by the same letters symbols and their 
$h$-Weyl quantizations. We shall consider $P$ as a closed \op{}:
$L^2\to L^2$ with domain $H(m):=(P-z_0)\inv L^2$ (see \cite{Ha2}).
From the discussion above, we get
\begin{itemize}
\item For every compact set $K\subset {\bf C}\setminus \Sigma $, we have 
$\sigma (P)\cap K=\emptyset$, when $h>0$ is small enough.
\item $\sigma (P)\cap \widetilde{\Omega }$ is discrete when $h>0$ is small 
enough.
\item $\sigma (\widetilde{P})\cap \widetilde{\Omega }=\emptyset$ when $h>0$  
is small enough.
\end{itemize}

\par In view of the last property and \no{gf.4}, we also have (for 
$h>0$ small enough),
\begin{prop}\label{gf2}
For $z\in \widetilde{\Omega }$, we have that 
$$
P(z):=(\widetilde{P}-z)^{-1}(P-z)=1+K(z),
$$
where $K(z)$ is a trace class \op{}. Moreover, 
$$
z\in \sigma (P)\Leftrightarrow 0\in \sigma (P(z)).
$$
\end{prop}

\par Notice that $K(z)=(\widetilde{P}-z)^{-1}(P-\widetilde{P})$ is the 
quantization of a symbol belonging to the intersection of $S(\widetilde{m})$ 
for all order \fu{}s $\widetilde{m}$. 

\section{Some functional calculus}\label{fc}
\setcounter{equation}{0}

Let $P=1+K$, $K\in {\rm Op}_h(S(m))$, where $m\in C^\infty ({\bf 
R}^{2n};]0,\infty [)$ is an integrable order \fu{}. We also assume that 
$K=k_0+hk_1+...$ in $S(m)$ on the symbol level. We shall review some 
functional calculus for $Q=P^*P$ and more generally for a \sa{} \op{} 
$Q\ge 0$ with $Q\sim q+hq_1+...$ on the symbol level, with $Q-1\in S(m)$, 
$q\ge 0$.

\par Let $\psi \in C_0^\infty ({\bf R})$. For $h\ll \alpha \ll 1$ we shall 
study the properties of $\psi (\alpha \inv Q)$.

\par To this end we shall consider $\alpha \inv Q$ as a symbol with 
$h/\alpha $ as a new semiclassical parameter, after a suitable dilation in 
phase space. More precisely, we make the change of variables 
$$
x=\alpha ^{1\over 2}\widetilde{x},\ D_x=\alpha ^{-{1\over 2}}D_{\widetilde{x}}
$$
and write 
\ekv{fc.25}
{{1\over \alpha }Q(x,hD_x;h)={1\over \alpha }Q(\alpha ^{1\over 
2}(\widetilde{x},{h\over \alpha }D_{\widetilde{x}});h),}
with symbol $\alpha \inv Q(\alpha ^{1/2}(\widetilde{x},\widetilde{\xi
});h)$ for the $h/\alpha $-quantization.
The lower order terms are ${\cal O}(h/\alpha )$ \ufly{} with all their 
derivatives, so we shall just look at the leading symbol
\ekv{fc.26}
{
{q(\alpha ^{1\over 2}(x,\xi ))\over \alpha },
}
where we dropped the tildes on the new variables. The (new) associated order 
\fu{} will be
\ekv{fc.27}
{
m(x,\xi ):=1+{q(\alpha ^{1\over 2}(x,\xi ))\over \alpha }\ge 1.
}
We have 
$$
\nabla m={(\nabla q)(\alpha ^{1\over 2}(x,\xi ))\over \alpha ^{1\over 
2}}\le C{q^{1\over 2}(\alpha ^{1\over 2}(x,\xi ))\over \alpha ^{1\over 
2}}\le Cm(x,\xi )^{1\over 2},
$$
$$
\nabla ^2m={\cal O}(1),
$$
so by Taylor's formula,
$$m(\rho )
=m(\mu )+{\cal O}(1)m(\mu )^{1\over 2}\vert \rho -\mu \vert +{\cal 
O}(1)\vert \rho -\mu \vert ^2, 
$$
and since $m(\mu )\ge 1$:
\ekv{fc.28}
{
m(\rho )\le C\langle \rho -\mu \rangle ^2 m(\mu ).
}
Hence $m$
is an order \fu{}, \ufly{} \wrt{} $\alpha $.

\par Similarly, we have the improved symbol estimates
\ekv{fc.29}
{
\nabla {q(\alpha ^{1\over 2}\rho )\over \alpha }={\cal O}(1)m^{1\over 2},
}
\ekv{fc.30}
{
\nabla ^2{q(\alpha ^{1\over 2}\rho )\over \alpha }={\cal O}(1),
}
\ekv{fc.31}
{
\nabla ^k {q(\alpha ^{1\over 2}\rho )\over \alpha }={\cal O}(\alpha 
^{{k\over 2}-1}),\ k\ge 2.
}
In particular, we have the standard symbol estimates 
\ekv{fc.32}
{
\nabla ^k {q(\alpha ^{1\over 2}\rho )\over \alpha }={\cal O}(1)m(\rho ).
}

It is therefore clear that we can apply the functional calculus in the 
version of \cite{HeSj} (see also \cite{DiSj}), to see that if $\psi \in 
C_0^\infty ({\bf R})$, and if we interpret $Q/\alpha $
as the \rhs{} of \no{fc.25}, then 
\ekv{fc.33}
{
\psi (\alpha \inv Q)= {\rm Op}_{{h\over \alpha 
},\widetilde{x}}(\widetilde{f}),
}
where 
\ekv{fc.34}
{
\widetilde{f}=\sum_{0}^\infty  ({h\over \alpha })^\nu f_\nu 
(\widetilde{x},\widetilde{\xi }),\hbox{ in }S(m\inv ),
}
with $f_0(\widetilde{x},\widetilde{\xi })=\psi (\alpha \inv q(\alpha 
^{1/2}(\widetilde{x},\widetilde{\xi })))$,
\ekv{fc.35}
{
f_\nu =\sum_{j\le j(\nu )}a_{j,\nu }(\widetilde{x},\widetilde{\xi 
},\alpha )\psi ^{(j)}(\alpha \inv q(\alpha 
^{1/2}(\widetilde{x},\widetilde{\xi }))),\ a_{j,\nu }\in S(1).
}
\begin{prop}\label{fc5}
Let $\widetilde{m}=\widetilde{m}_\alpha (\widetilde{x},\widetilde{\xi })$ 
be an order \fu{}, \ufly{} \wrt{} $\alpha $, such that 
$\widetilde{m}(\widetilde{x},\widetilde{\xi })=1$ for $\alpha \inv 
q(\alpha ^{1/2}(\widetilde{x},\widetilde{\xi }))\le \sup {\rm supp\,}\psi 
+1/C$, for some $C>0$ that is \indep{} of $\alpha $. Then \no{fc.34} holds 
in $S(\widetilde{m})$, for $h$ and $h/\alpha $ \sufly{} small.
\end{prop}
\begin{proof}
Write $q_\alpha =q(\alpha ^{1/2}(\widetilde{x},\widetilde{\xi 
}))/\alpha $, $Q_\alpha =\alpha \inv Q(\alpha ^{1/2}(\widetilde{x},{h\over 
\alpha }D_{\widetilde{x}});h)$ and drop the tildes. Let 
$\widehat{q}_\alpha \in S(m)$ be such that
$\sup{\rm 
supp\,}\psi +1/(5C)\le \widehat{q}_\alpha $, and be equal to 
$q_\alpha $ when $q_\alpha \ge \sup{\rm supp\,}\psi +2/(5C)$.
Let $\chi 
_\alpha \in S(1)$ be equal to 1 when $q_\alpha \le \sup{\rm supp\,}\psi 
+3/(5C)$ and equal to 0
 when $q_\alpha \ge \sup{\rm supp\,}\psi 
+4/(5C)$. We use the same symbols $q_\alpha $, $\widehat{q}_\alpha $, 
$\chi _\alpha $ to denote the $h/\alpha $
quantizations.

\par Let $\widetilde{\psi }$ be an almost \hol{} extension of $\psi $ and 
recall the Cauchy-Green-Riemann-Stokes formula, in the \op{} sense 
(\cite{HeSj}, \cite{DiSj}, \cite{Dy}):
$$
\psi (q_\alpha )={1\over \pi }\int (z-q_\alpha )\inv {\partial \widetilde{\psi 
}(z)\over \partial \overline{z}} 
L(dz),
$$
where $L(dz)$ denotes Lebesgue-measure.
For $z$ in a \neigh{} of ${\rm supp\,}\widetilde{\psi }$, we write 
$$
(z-q_\alpha )\inv =(z-q_\alpha )\inv \circ \chi _\alpha 
+(z-\widehat{q}_\alpha )\inv \circ (1-\chi _\alpha )-(z-q_\alpha )\inv 
(\widehat{q}_\alpha -q_\alpha )(z-\widehat{q}_\alpha )\inv (1-\chi 
_\alpha ).
$$
Then, since the middle term is \hol{} near the support of 
$\widetilde{\psi }$,
$$
\psi (q_\alpha )=\psi (q_\alpha )\circ \chi _\alpha -{1\over \pi }\int 
(z-q_\alpha )\inv (\widehat{q}_\alpha -q_\alpha )(z-\widehat{q}_\alpha 
)\inv (1-\chi _\alpha ) {\partial \widetilde{\psi }\over \partial 
\overline{z}}L(dz).
$$
Here the symbol of $\psi (q_\alpha )\circ \chi _\alpha $ has the \asy{} 
expansion \no{fc.34} in $S(\widetilde{m})$, thanks to the extra factor 
$\chi _\alpha $ and with the same terms given by \no{fc.35}. For $z\in 
{\rm neigh\,}({\rm supp\,}\widetilde{\psi })$, $(z-\widehat{q}_\alpha 
)\inv \in{\rm Op\,}(1/m)$ depends \hol{}ally on $z$ and thanks to the 
factor $\widehat{q}_\alpha -q_\alpha $, whose support on the symbol level 
is separated from that of $1-\chi _\alpha $ by some fixed positive 
distance, we know that $(\widehat{q}_\alpha -q_\alpha )(z-\widehat{q}_\alpha 
)\inv (1-\chi _\alpha )\in (h/\alpha )^N{\rm Op\,}(S(\widetilde{m}))$ for 
any $N\ge 0$ and any $\widetilde{m}$ as in the proposition. Combining this 
with the estimates for the symbol of $(z-q_\alpha )\inv$ from the
Beals lemma as in \cite{HeSj} (see also 
\cite{DiSj}, Proposition 8.6), we get the proposition. 
\end{proof}

\par We next apply the functional result to the study of certain 
determinants. Let $\chi \in C_0^\infty ([0,+\infty [;[0,+\infty [)$, $\chi 
(0)>0$ and extend $\chi $ to $C_0^\infty ({\bf R};{\bf C})$ in such a 
way that $\chi (t)>0$ near 0 and $t+\chi (t)\ne 0$, $\forall t\in{\bf R}$. 
We want to study $\ln\det (Q+\alpha \chi (\alpha \inv Q))$, when $h\ll 
\alpha \ll 1$. Let us first recall from \cite{MeSj} that if 
$\widetilde{Q}={\rm Op}_h(\widetilde{q})$ with $\widetilde{q}\in S(1)$, 
$\widetilde{q}>0$, $\widetilde{q}-1\in S(\widetilde{m})$, where 
$\widetilde{m}$ is an integrable order \fu{}, then 
\ekv{fc.36}
{
\ln\det \widetilde{Q}={1\over (2\pi h)^n}(\iint\ln \widetilde{q}(x,\xi 
)dxd\xi +{\cal O}(h)).
}
In fact, let $\widetilde{Q}_t=(1-t)1+t\widetilde{Q}$, so that 
$\widetilde{Q}_0=1$, $\widetilde{Q}_1=\widetilde{Q}$. Then by standard 
elliptic calculus, with $\widetilde{q}_t=(1-t)+t\widetilde{q}$, we have 
\begin{eqnarray*}
{d\over dt}\ln\det \widetilde{Q}_t&=&\tr \widetilde{Q}_t\inv {d\over 
dt}\widetilde{Q}_t\\
&=&{1\over (2\pi h)^n}(\iint \widetilde{q}_t\inv {d\over 
dt}\widetilde{q}_t\, dxd\xi 
+{\cal O}(h))\\
&=&{1\over (2\pi h)^n}({d\over dt}\iint \ln \widetilde{q}_t(x,\xi ) dxd\xi 
+{\cal O}(h)),
\end{eqnarray*}
and integrating from $t=0$ to $t=1$, we get \no{fc.36}.

\par For $\alpha =\alpha _1>0$ fixed with $\alpha _1\ll 1$, this applies 
to $Q+\alpha _1\chi (\alpha _1\inv Q)$ and we get 
\ekv{fc.37}
{
\ln\det (Q+\alpha _1\chi (\alpha _1\inv Q))={1\over (2\pi h)^{n}}(\iint 
\ln (q+\alpha _1\chi (\alpha _1\inv q))dxd\xi +{\cal O}(h)).}

\par We have for $t>0$ and $E\ge 0$:
$$
{d\over dt}\ln (E+t\chi ({E\over t}))={1\over t}\psi ({E\over t}),
$$
with 
$$
\psi (E)={\chi (E)-E\chi '(E)\over E+\chi (E)}.
$$ 
Now for $h\ll\alpha \le t\le \alpha _1$, we get from Proposition \ref{fc5} by dilatation:
\begin{eqnarray}\label{fc.37.5}
&&\hskip -2truecm {d\over dt}\ln \det (Q+t\chi (t\inv Q))=\tr t\inv \psi (t\inv Q)\\
&=&\Big({t\over 2\pi h}\Big)^n\Big(\iint {1\over t}\psi ({q(t^{1\over 
2}(\widetilde{x},\widetilde{\xi }))\over
t})d\widetilde{x}d\widetilde{\xi }
+{\cal O}({h\over t}){1\over t}\iint 
\widehat{\chi }({q(t^{1\over 
2}(\widetilde{x},\widetilde{\xi }))\over t})d\widetilde{x}d\widetilde{\xi }
\nonumber\\
&&
+{\cal O}(1){1\over t}({h\over t})^\infty \iint (1+{\rm 
dist\,}((\widetilde{x},\widetilde{\xi });{\rm supp\,}\widehat{\chi }
({q(t^{1\over 
2}(\cdot ))\over t})))^{-N}d\widetilde{x}d\widetilde{\xi }\Big),
\nonumber
\end{eqnarray}
where $0\le \widehat{\chi }\in C_0^\infty ({\bf R})$ is equal to one on 
some interval containing $[0,\sup {\rm supp\,}\psi ] $ and the last term is coming from the 
``remainder'' in the asymptotic development (\ref{fc.34}). 

\par We are interested in the integral of this quantity from $t=\alpha $ 
to $t=\alpha _1$. Let us first treat the leading term
\begin{eqnarray*}
({t\over 2\pi h})^n\iint {1\over t}\psi ({q(t^{1\over 
2}(\widetilde{x},\widetilde{\xi }))\over 
t})d\widetilde{x}d\widetilde{\xi }
&=&{1\over (2\pi h)^n}\iint {1\over t}\psi ({q(x,\xi )\over t})dxd\xi \\
&=&{1\over (2\pi h)^n}\iint {d\over dt}\ln (q+t\chi ({q\over t}))dxd\xi ,
\end{eqnarray*}
and integrating this from $t=\alpha $ to $t=\alpha _1$, we get 
\ekv{fc.38}
{
\Big[ 
{1\over (2\pi h)^n}\iint \ln (q+t\chi ({q\over t}))dxd\xi 
\Big]_{t=\alpha }^{\alpha _1}.
}

\par The second term in \no{fc.37.5} is
\begin{eqnarray*}
{\cal O}(1)({t\over h})^n{h\over t^2}\iint \widehat{\chi }\left({1\over 
t}q\big(t^{1\over 2}(\widetilde{x},\widetilde{\xi 
})\big) \right)d\widetilde{x}d\widetilde{\xi }&=& {\cal O}(1){1\over h^n}{h\over 
t^2}\iint \widehat{\chi }\big({1\over t}q(x,\xi )\big)dxd\xi \\
&\le& {\cal O}(1)h^{-n}{h\over t^2}\iint_{q(x,\xi )\le C^2t}dxd\xi .
\end{eqnarray*}
Integrating this from $t=\alpha $ to $t=\alpha _1$, we get
\begin{eqnarray}\label{fc.39}
&\displaystyle {\cal O}(1)h^{1-n}\iint \int_{\max (\alpha ,q/C)\le t\le \alpha _1}{1\over 
t^2}dtdxd\xi &\\
= &\displaystyle {\cal O}(1)h^{1-n}\iint_{q(x,\xi )\le C\alpha _1}\left( {1\over \max \big(\alpha 
,q(x,\xi )/C\big)}-{1\over \alpha _1}\right)dxd\xi& \nonumber \\
\le &\displaystyle {\cal O}(1)h^{-n} \iint_{q(x,\xi )\le C\alpha _1}{h\over \alpha +q(x,\xi 
)}dxd\xi&  .\nonumber
\end{eqnarray}

\par When estimating the third term in \no{fc.37.5} we consider separately 
the regions $\vert (x,\xi )\vert \le C$ and $\vert (x,\xi )\vert > C$ for 
some large $C\ge 1$. Consider first the region $\vert (x,\xi )\vert \le C$. 
Put 
$$
d_t(\widetilde{x},\widetilde{\xi })={\rm 
dist\,}((\widetilde{x},\widetilde{\xi }),\{ (y,\eta );\, {1\over 
t}q(t^{1\over 2}(y,\eta ))\le \widehat{C}\}),\ \widehat{C}=\sup {\rm 
supp\,}\widehat{\chi }.
$$
For $(y,\eta )$ with $q(t^{1\over 2}(y,\eta ))/t\le \widehat{C}$, we 
have $\nabla (q(t^{1\over 2}(y,\eta ))/t)={\cal O}(1)$ and 
$\nabla ^2 (q(t^{1\over 2}(y,\eta ))/t)={\cal O}(1)$, so by Taylor 
expanding at $(y,\eta )$ we get 
\begin{eqnarray*}
{1\over t}q(t^{1\over 2}(\widetilde{x},\widetilde{\xi }))\le {\cal 
O}(1)(1+d_t(\widetilde{x},\widetilde{\xi 
})+d_t(\widetilde{x},\widetilde{\xi })^2)\le {\cal 
O}(1)(1+d_t(\widetilde{x},\widetilde{\xi }))^2.
\end{eqnarray*}
The contribution to the third term in \no{fc.37.5} from $t^{1/2}\vert 
(\widetilde{x},\widetilde{\xi })\vert \le C$ is therefore
\begin{eqnarray}\label{fc.40}
&&{\cal O}_N(1)({h\over t})^\infty {1\over t}\iint_{\vert 
(\widetilde{x},\widetilde{\xi })\vert \le Ct^{-1/2}} (1+{1\over 
t}q(t^{1\over 2}(\widetilde{x},\widetilde{\xi 
})))^{-N}d\widetilde{x}d\widetilde{\xi }\\
&&=
{\cal O}_{M,N}(1){1\over h^n}({h\over t})^M {1\over t}\iint_{\vert (x,\xi 
)\vert \le C}(1+{1\over t}q(x,\xi ))^{-N} dxd\xi .
\nonumber
\end{eqnarray}
We integrate this from $t=\alpha $ to $\alpha _1$, so we want to estimate
$$
h^M\int_\alpha ^{\alpha _1}{1\over t^{M+1}(1+{1\over t}q(x,\xi ))^N}dt.
$$
If $q\le \alpha $, we get 
$$
{\cal O}(1)h^M\int_\alpha ^{\alpha _1}{1\over t^{M+1}}dt={\cal O}(1)(({h\over 
\alpha })^M).
$$
If $\alpha <q\le \alpha _1$, we get 
\begin{eqnarray*}
h^M\int_\alpha ^q {1\over t^{M+1}}(1+{q\over t})^{-N}dt +h^M\int_q^{\alpha 
_1} {1\over t^{M+1}}(1+{q\over t})^{-N}dt
\backsim h^M\int_\alpha ^q {t^N\over t^{M+1}q^N}dt+({h\over q})^M,
\end{eqnarray*}
with the symbol $\backsim$ indicating same order of magnitude. 
Choose 
$N=M+1$, to get 
$$
\backsim ({h\over q})^M.
$$
For $\alpha _1\le q$, we get 
$${\cal O}(1)h^M\int_\alpha ^{\alpha _1}{t^N\over 
t^{M+1}q^N}dt={\cal O}(1)h^M$$
with the same choice of $N$.
Thus the expression \no{fc.40} is 
$$
{{\cal O}(1)\over h^n}\iint_{\vert (x,\xi )\vert \le C}\Big( {h\over \alpha 
+q(x,\xi )}\Big) ^M dxd\xi ,\ \forall M\ge 0.
$$

We next look at the contribution to the last term in \no{fc.37.5} from the 
region $t^{1\over 2}\vert (\widetilde{x},\widetilde{\xi })\vert >C$, which 
is 
\begin{eqnarray*}
{\cal O}(1)\Big({h\over t}\Big)^\infty {1\over t}\iint_{t^{1\over 2}\vert 
(\widetilde{x},\widetilde{\xi })\vert \ge C}(1+\vert 
(\widetilde{x},\widetilde{\xi })\vert
)^{-N}d\widetilde{x}d\widetilde{\xi }={\cal O}(1)({h\over 
t})^Mt^{\widetilde{N}},\ \forall M,\widetilde{N},\\
={\cal O}(1)h^M,\ \forall M.
\end{eqnarray*}

\par Summing up, we have proved
\begin{prop}\label{fc6}
Let $\chi \in C_0^\infty ({\bf R})$ with $\chi (0)>0$ and $\chi (t)\ge 0$ 
for $t\ge 0$. Then for $h\ll \alpha \ll 1$:
\begin{eqnarray}\label{fc.41}
&&\ln \det (Q+\alpha \chi (\alpha \inv Q))=
\\&&{1\over (2\pi h)^n}\Big(\iint \ln 
(q+\alpha \chi ({q\over \alpha }))dxd\xi +{\cal O}(1)\iint_{\vert (x,\xi 
)\vert \le C}{h\over \alpha +q(x,\xi )}dxd\xi \Big) +{\cal O}(h^\infty 
).\nonumber
\end{eqnarray}
\end{prop}

\par Most of the proof was based on \no{fc.37.5}, of which the second part is 
valid for any $\psi \in C_0^\infty ({\bf R})$. The estimates leading to 
the preceding proposition, also give
\begin{prop}\label{fc7}
Let $\chi \in C_0^\infty ({\bf R})$, and choose $\widehat{\chi }\in 
C_0^\infty ({\bf R};[0,1])$ equal to $1$ on an interval containing 0 
and ${\rm supp\,}\chi $. Then for $0<h\ll \alpha \ll 1$, we have
\begin{eqnarray}\label{fc.42}
\hskip -1truecm \tr \chi (\alpha \inv Q)&=&{1\over (2\pi h)^n}\Big(\iint \chi ({q(x,\xi )\over 
\alpha })dxd\xi +{\cal O}({h\over \alpha })\iint \widehat{\chi }({q(x,\xi 
)\over \alpha })dxd\xi\\&&+{\cal O}_{N,M}(1)({h\over \alpha })^M\iint_{\vert 
(x,\xi )\vert \le C}(1+{q\over \alpha })^{-N}dxd\xi +{\cal O}(h^\infty 
)\Big).\nonumber 
\end{eqnarray}
\end{prop}

\par Put 
\ekv{fc.43}
{
V(t)=\iint_{q(x,\xi )\le t}dxd\xi ,\ 0\le t\le {1\over 2},
}
so that $0\le V(t)$ is an increasing \fu{}. By introducing an assumption 
on $V(t)$, we shall make \no{fc.42} more explicit and replace \no{fc.41} by 
a more explicit estimate.

\par  We assume that there exists $\kappa\in ]0,1]$ such that
\ekv{fc.44}
{
V(t)={\cal O}(1)t^{\kappa},\ 0\le t\le {1\over 2}.
}
 
\par The first integral in \no{fc.41} can be written
\begin{eqnarray*}
\iint \ln (q+\alpha \chi ({q\over \alpha }))dxd\xi &=& \iint \ln
(q) dxd\xi +{\cal O}(1)\iint_{q\le C\alpha }\ln {1\over q}dxd\xi \\
&=& \iint \ln (q)dxd\xi +{\cal O}(1)\int_0^{C\alpha}\ln {1\over
q}dV(q),
\end{eqnarray*}
so \no{fc.41} gives
\eekv{fc.45}
{&& \ln \det (Q+\alpha \chi ({1\over \alpha }Q))}
{&=& {1\over (2\pi h)^n}(\iint \ln (q)dxd\xi +{\cal
O}(1)\int_0^{C\alpha }\ln {1\over q}dV(q)+{\cal O}(1)\int_0^{1\over
2}{h\over \alpha +q}dV(q)+{\cal O}(h^\infty )). }

\par Similarly \no{fc.42} can be written
\eekv{fc.46}
{{\rm tr\,}\chi (\alpha \inv Q)&=&{1\over (2\pi h)^n}\Big( \int
\chi ({q\over \alpha })dV(q)+{\cal O}({h\over \alpha
})\int\widehat{\chi }({q\over \alpha })dV(q)} 
{&&+{\cal O}_{N,M}(1)({h\over \alpha })^M\int_0^{\alpha
_1}(1+{q\over \alpha })^{-N}dV(q)+{\cal O}(h)\Big) .}
In particular, the number $N(\alpha )$ of \ev{}s of $Q$ in
$[0,\alpha ]$ satisfies
\ekv{fc.47}
{
N(\alpha )\le {\cal O}(1)\big( h^{-n}\int_0^{\alpha _1}(1+{q\over
\alpha })^{-N}dV(q) +h^{1-n}).
}
\begin{prop}\label{fc8}
Under the assumption \no{fc.44}, we have
\ekv{fc.48}
{
\int_0^\alpha  \ln (q)dV(q)={\cal O}(\alpha ^{\kappa}\ln \alpha
), }
\ekv{fc.49}
{
\int_0^{\alpha _1}{h\over \alpha +q}dV(q)=
\begin{cases}{\cal O}(\alpha ^{\kappa} 
{h\over \alpha }),\hbox{ for }\kappa<1,\cr \cr {\cal O}(h\ln{1\over 
\alpha }),\hbox{ when }\kappa=1,\end{cases} }
\\
\ekv{fc.50}
{
N(\alpha )={\cal O}(\alpha ^{\kappa}h^{-n}+h^{1-n}).
}
\end{prop}
\begin{proof}
This follows by straight forward calculations, starting with an
integration by parts.

$$
\int_0^\alpha  \ln (q)dV(q)=[\ln (q)V(q)]_0^\alpha  -\int_0^\alpha 
{1\over q} V(q)dq={\cal O}(\alpha ^{\kappa}\ln (\alpha )),
$$
\begin{eqnarray*}
\int_0^{\alpha _1}{h\over \alpha +q}dV(q)&=&[{h\over \alpha
+q}V(q)]_0^{\alpha _1}+\int_0^{\alpha _1}{h\over (\alpha
+q)^2}V(q)dq\\ \medskip
&=&{\cal O}(h)+{\cal O}(1)\int_0^{\alpha _1}{hq^{\kappa}\over
(\alpha +q)^2}dq\\ \medskip
&=&{\cal O}(h)+{\cal O}(1)h\alpha ^{\kappa-1}\int_0^{\alpha
_1/\alpha }{\widetilde{q}^{\kappa}\over
(1+\widetilde{q})^2}d\widetilde{q}\\ \medskip
&=&\begin{cases}{\cal O}(h\alpha ^{\kappa-1}),\ 0<\kappa<1,\cr  \cr {\cal 
O}(h\ln{1\over \alpha }),\ \kappa=1.\end{cases}
\end{eqnarray*}
To get \no{fc.50}, we use \no{fc.47} and the estimate
\begin{eqnarray*}
\int_0^{\alpha _1}(1+{q\over \alpha })^{-N}dV(q)&=& [(1+{q\over
\alpha })^{-N}V(q)]_0^{\alpha _1}+N\int_0^{\alpha _1}(1+{q\over
\alpha })^{-(N+1)}V(q){dq\over \alpha }\\
&=&{\cal O}(1)\alpha ^N+{\cal O}(1)\int_0^\infty 
(1+\widetilde{q})^{-(N+1)}\widetilde{q}^{\kappa }d\widetilde{q}\alpha ^{\kappa}\\
&=&{\cal O}(1)\alpha ^{\kappa}.
\end{eqnarray*}
\end{proof}

\par Since we will always assume that $h\ll \alpha \ll 1$, and since
$\kappa\le 1$, we can simplify \no{fc.50} to
\ekv{fc.53}
{N(\alpha )={\cal O}(1)\alpha ^{\kappa}h^{-n}.}

\begin{prop}\label{fc9}
Assume \no{fc.44}. Under the assumptions of Proposition \ref{fc6}, we have 
\ekv{fc.51}
{
\ln \det (Q+\alpha \chi (\alpha \inv Q))={1\over (2\pi h)^n}(\iint
\ln (q) dxd\xi +{\cal O}(1) \alpha ^{\kappa}\ln \alpha ). }

\par Under the assumptions of Proposition \ref{fc7}, we have
\ekv{fc.52}
{
{\rm tr\,}\chi (\alpha \inv Q)={1\over (2\pi h)^n}(\iint \chi
({q(x,\xi )\over \alpha })dxd\xi +{\cal O}(1)\alpha ^{\kappa 
}{h\over \alpha }). }

\par For $0<h\ll \alpha \ll 1$, the number $N(\alpha )$ of \ev{}s
of $Q$ in $[0,\alpha ]$ satisfies \no{fc.50}.
\end{prop}

\par Notice that when $Q\ge 0$:
$$
\ln \det Q\le \ln \det (Q+\alpha \chi (\alpha \inv Q)),
$$
so \no{fc.51} with $\alpha =Ch$, $C\gg 1$, gives an upper bound which is 
more precise than the one in \cite{MeSj}.
\section{Grushin \pb{} for the unperturbed \op{}}\label{gu}
\setcounter{equation}{0}

\par In Section \ref{gf} we introduced the \op{}
$$
P(z)=1+K(z),\quad K(z)\in {\rm Op}_h(S(m))
$$
where $m$ is an integrable order \fu{}, so that $K(z)$ is a trace class 
\op{}. Here $P(z)$ depends \hol{}ally on $z\in \widetilde{\Omega 
}\Subset  {\bf C}$, where $\widetilde{\Omega }$ is open. 
Also recall that $P(z)=(\widetilde{P}-z)\inv (P-z)$. We are 
interested in the spectrum of small random perturbations of $P$; $P_\delta 
=P+\delta Q$. Correspondingly, we get $P_\delta (z)=(\widetilde{P}-z)\inv 
(P_\delta -z)=1+K_\delta (z)$, and the main work in later sections will be to 
study $\vert \det (1+K_\delta (z))\vert $. The upper bounds will be fairly 
simple to get, and the delicate point will be to get lower bounds. As a 
preparation for this more delicate step, we here study a Grushin \pb{} for 
the unperturbed \op{} $P(z)$. In this section $z\in\widetilde{\Omega }$ 
will be fixed and we simply write $P$ instead $P(z)$.

\par Let $e_1,e_2,...$ be an orthonormal (ON) basis of eigenvectors of $P^*P$ 
and let $0\le \lambda _1\le \lambda _2\le ...$ be the corresponding \ev{}s. 
(Strictly speaking, if we want the \ev{}s to form an increasing sequence, 
the set of indices $j$ should be of the form $J=J_0\cup J_1\cup J_2$, with 
\begin{itemize}
\item $J_0={\bf N}$ or a finite set, $0\le \lambda _j<1$, for $j\in J_0$, 
\item $J_1={\bf N}$ or a finite set, $\lambda _j=1$, $j\in J_1$
\item $J_2=-{\bf N}$ or a finite set, $\lambda _j>1$, $j\in J_2$.\end{itemize}
We will only be concerned with finitely many indices from $J_0$.)

\par Since $P=P(z)$ is a Fredholm operator of index zero by Proposition \ref{gf2}, we know that $PP^*$ and $P^*P$ have the same 
number $N_0$ of \ev{}s equal to $0$.  Let $f_1,...,f_{N_0}$ be an ON basis of
$\ker(PP^*)$.  For $j>N_0$, we have $\lambda _j>0$ and $Pe_j$ is an
eigenvector of $PP^*$ with \ev{} $\lambda _j$: $$ PP^*Pe_j=\lambda _jPe_j. 
$$ 
Using standard notations for norms and scalar products, we
put $f_j=\Vert Pe_j\Vert \inv Pe_j$.  Then $\{ f_j\}_{j\in J}$ is an ON
system of eigenvectors of $PP^*$, with $PP^*f_j=\lambda _jf_j$.  Let $f\in
L^2({\bf R}^n)$ with $(f\vert f_j)=0$ for all $j\in J$.  Then $(P^*f\vert
e_j)=(f\vert Pe_j)$.  If $j\le N_0$, we get $(P^*f\vert e_j)=(f\vert 0)=0$,
and if $j\ge N_0+1$, we get $(P^*f\vert e_j)=\Vert Pe_j\Vert (f\vert
f_j)=0$. Hence $P^*f=0$, so $f\in \ker(PP^*)$ and hence $f=0$ is zero 
since $\ker(PP^*)$ is the span of $f_1,...,f_{N_0}$. We conclude 
that $\{ f_j\}_{j\in J}$ is an ON \it basis \rm of eigenvectors of $PP^*$. 
By construction, we have $Pe_j=w_{j}f_j$ with $0\le w_{j} =\Vert Pe_j\Vert $. Then 
$$
w_{j}^2=(P^*Pe_j\vert e_j)=\lambda _j,
$$
so $w_{j}=\sqrt{\lambda _j}$ and it follows that,
\ekv{gu.1}
{Pe_j=\sqrt{\lambda _j}f_j,}
\ekv{gu.2}{P^*f_j=\sqrt{\lambda _j}e_j}
for all $j\in J$. 

\par Let $0<\alpha \ll 1$ and let $N=N(\alpha )$ be given by 
\ekv{gu.3}
{
\lambda _j\le \alpha \Leftrightarrow j\le N(\alpha ).
}
Define 
$$
R_+:L^2({\bf R}^n)\to {\bf C}^N,\quad R_-:{\bf C}^N\to L^2({\bf R}^n)
$$
by 
$$
R_+u=\sqrt{\alpha }((u\vert e_j))_{j=1}^N,\quad R_-u_-=\sqrt{\alpha 
}\sum_1^N u_-(j)f_j,
$$
and put 

\ekv{gu.4}
{ {\cal P}=\begin{pmatrix}P &R_-\cr R_+ &0\end{pmatrix}:L^2({\bf R}^n)\times {\bf C}^N\to
L^2({\bf R}^n)\times {\bf C}^N .}
If $u=\sum_{j\in J}u_je_j$, $u_-=(u_-(j))_{j=1}^N,$
we get
$$
{\cal P}
\begin{pmatrix}u \cr u_-\end{pmatrix}=
\begin{pmatrix}\sum_J\sqrt{\lambda 
_j}u_jf_j+\sum_1^N \sqrt{\alpha }u_-(j)f_j \cr (\sqrt{\alpha }u_j)_{j=1}^N\end{pmatrix},
$$
and we conclude that 
\eeekv{gu.5}
{\vert \det {\cal P}\vert &=&(\prod_1^N \vert \det 
\begin{pmatrix}
\sqrt{\lambda 
_j}&\sqrt{\alpha } \cr \sqrt{\alpha } &0\end{pmatrix}\vert )(\prod_{N<j\in 
J}\sqrt{\lambda _j})}
{&=&\alpha ^N\prod_{N<j\in J}\sqrt{\lambda _j}}{&=&\alpha ^{N\over 2}\prod_J \max 
(\sqrt{\alpha },\sqrt{\lambda _j}).}
Notice that
\ekv{gu.6}
{
\vert \det P\vert =\prod_J \sqrt{\lambda _j}.
}

\par Let $\delta _j(k)=\delta _{j,k},\ 1\le j,k\le N$. Then ${\cal P}$ 
maps ${\bf C}e_j\times {\bf C}\delta _j$ to ${\bf C}f_j\times {\bf C}\delta 
_j$ and has the corresponding matrix 
$$
\begin{pmatrix}\sqrt{\lambda _j}&\sqrt{\alpha }\cr \sqrt{\alpha } &0\end{pmatrix}.$$ 
The 
inverse is given by 
$$
\begin{pmatrix} 0 &{1\over \sqrt{\alpha }}\cr {1\over \sqrt{\alpha }} 
&-{\sqrt{\lambda _j}\over \alpha } \end{pmatrix},
$$
so if $v=\sum_J v_j f_j$, $v_+=\sum_1^N v_+(j)\delta _j$, then (writing as before ${\cal E}$ for ${\cal P}^{-1}$)
\ekv{gu.7}
{
{\cal E}\begin{pmatrix}v\cr v_+\end{pmatrix}=
\begin{pmatrix}\sum_{N+1}^\infty  {1\over 
\sqrt{\lambda _j}}v_je_j+{1\over \sqrt{\alpha }}\sum_1^N v_+(j)e_j\cr 
\sum_1^N{1\over \sqrt{\alpha }}v_j\delta _j-\sum_1^N {\sqrt{\lambda 
_j}\over \alpha }v_+(j)\delta _j\end{pmatrix}
 =\begin{pmatrix}
 E &E_+\cr E_- &E_{-+}\end{pmatrix}
 \begin{pmatrix}v\cr v_+\end{pmatrix},}
where
\eeeekv{gu.8}
{&&E_+v_+={1\over \sqrt{\alpha }}\sum_1^N v_+(j)e_j,}
{&&E_-v={1\over \sqrt{\alpha }}\sum_1^N v_j\delta _j,}
{&&E_{-+}=-{1\over \alpha }{\rm diag\,}(\sqrt{\lambda _j}),}
{&&\Vert E\Vert ,\Vert E_+\Vert ,\Vert E_-\Vert ,\Vert E_{-+}\Vert 
\le {1\over 
\sqrt{\alpha }}.}
From \no{gu.5}--\no{gu.8}, we see that 
\ekv{gu.9}
{
\vert \det P\vert=\vert \det{\cal P}\vert \vert \det{E_{-+}}\vert , 
}
as we already know from \no{dg.11}. 

\par We next study $\vert \det {\cal P}\vert $, when $h\ll \alpha \ll 1$. 
The formula \no{gu.5} can be written
\ekv{gu.10}
{
\vert \det{\cal P}\vert ^2=\alpha ^N\det 1_\alpha (P^*P),
}
where $1_\alpha (t)=\max (\alpha ,t)$. Let $\chi \in C_0^\infty 
([0,2[;[0,1])$ be equal to 1 on $[0,1]$. Then for $t\ge 0$,
\ekv{gu.11}
{
t+{\alpha \over 4}\chi ({4t\over \alpha })\le 1_\alpha (t)\le t+\alpha 
\chi ({t\over \alpha }).
}

\par In the following, we assume that $Q=P^*P$ satisfies the assumptions 
of Section \ref{fc}, including \no{fc.44}, and choose $h\ll 
\alpha \ll 1$. Then we know that 
\ekv{gu.12}
{N(\alpha )={\cal O}(\alpha ^{\kappa}h^{-n}),}
and Proposition \ref{fc9} in combination with \no{gu.10}--\no{gu.12} show 
that
\ekv{gu.13}
{
\ln \vert \det {\cal P}\vert ^2={1\over (2\pi h)^n}(\iint \ln (q)dxd\xi 
+{\cal O}(1)\alpha ^{\kappa}\ln \alpha ).
}
As noticed after Proposition \ref{fc9}, we also have the upper bound
\ekv{gu.14}
{
\ln\det P^*P\le {1\over (2\pi h)^n}(\iint \ln(q)dxd\xi +{\cal 
O}(1)\alpha ^{\kappa}\ln {1\over\alpha}).
}

\section{The Hilbert-Schmidt norm of 
a Gaussian random matrix.
}\label{hs}
\setcounter{equation}{0}

Let $\alpha (\omega )$ be a complex Gaussian random variable with density
\ekv{hs.1}
{
{1\over \pi \sigma ^2}e^{-\vert \alpha \vert ^2/\sigma ^2}L(d\alpha ),\ 
L(d\alpha )=d\Re \alpha\, d\Im \alpha ,  
}
that is a ${\cal N}(0,\sigma ^2 )$-law with $\sigma ^2 $ denoting the variance. 
The distribution of $\vert \alpha 
(\omega )\vert ^2$ is 
\ekv{hs.2}
{\mu d\alpha ={1\over s}e^{-r/s}H(r)dr,}
where $s=\sigma ^2$ and $H(r)$ denotes the standard Heaviside \fu{}. 
Notice that 
$$
\Vert \vert \alpha \vert ^2\Vert _{L^1}=\langle \vert \alpha \vert 
^2\rangle =\sigma ^2.
$$
\par  Let 
$\alpha _j(\omega )$, $j=1,2,...$ be independent random variables as above 
with variance $\sigma _j^2$ and assume for simplicity that $\sigma _1\ge 
\sigma _j$ for all $j$. We also assume that 
\ekv{hs.3}
{
\sum_1^\infty  \sigma _j^2<\infty ,
}
implying the a.s. convergence of $\sum_1^\infty \vert \alpha _j(\omega 
)\vert ^2$.

\par We want to estimate the \proba{} that $\sum \vert \alpha _j(\omega 
)\vert ^2\ge a$. The \proba{} distribution of $\sum_1^\infty  \vert \alpha 
_j(\omega )\vert ^2$ is equal to $(\mu _1*\mu _2*....)dx$, where $\mu _j$
is given in \no{hs.2} with $s=s_j=\sigma _j^2$, so that 
\ekv{hs.4}{\sum_1^\infty  s_j<\infty .}

\par The \F{} \tf{} of $\mu _j$ is given by 
\ekv{hs.5}
{\widehat{\mu }_j(\rho )={1\over 1+is_j\rho },}
which has a simple pole at $\rho =i/s_j$. The probability that we are 
after, is 
\ekv{hs.6}
{
\int_a^\infty (\mu _1*\mu _2*...)dr={1\over 2\pi }\int \prod_1^\infty  
(\widehat{\mu }_j(\rho ))\overline{\widehat{1_{[a,\infty [}}}(\rho )d\rho ,
}
by Parseval's identity. Here 
\ekv{hs.7}
{
\widehat{1_{[a,\infty [}}(\rho )={1\over i(\rho -i0)}e^{-ia\rho },
}
so the \proba{} \no{hs.5} becomes
\ekv{hs.8}
{
{i\over 2\pi }\int_{-\infty }^\infty (\prod_1^\infty {1\over 
1+is_j\rho }){1\over \rho +i0}e^{ia\rho }d\rho .
} 
The assumption \no{hs.4} implies that the infinite product converges away 
from the poles $i/s_j$. For 
$\rho $ in a half plane $\Im \rho \le b<{1\over 2s_1}$, we have 
$$
\vert {1\over 1+is_1\rho }\vert \le {1\over ((1-bs_1)^2+s_1^2(\Re \rho 
)^2)^{1\over 2}},
$$
$$
\vert \prod _2^\infty {1\over 1+is_j\rho }\vert \le \prod _2^\infty {1\over 
1-bs_j}\le \exp (C_0\sum_2^\infty bs_j),
$$
where $C_0$
is a universal constant appearing in the estimate, 
$$
{1\over 1-t}\le e^{C_0t}, \ 0\le t\le {1\over 2}.
$$

\par Shifting the contour in \no{hs.8} from ${\bf R}$ to ${\bf R}+ib$ and 
choosing $b=1/(2s_1)$, we can estimate the \proba{} \no{hs.6} from above 
by 
\ekv{hs.9}
{
C(s_1)\exp [{C_0\over 2s_1}\sum_1^\infty  s_j -{1\over 2s_1}a],
}
where $C_0>0$ is the universal constant introduced above and $C(s_1)$ can be chosen 
\ufly{} \bdd{} on any compact subset of $]0,+\infty ]$.

\begin{remark}\label{hs1}\rm
W. Bordeaux Montrieux has used a more elementary argument in the case of
real matrices, by means of the Markov-Chebyschev inequality.  With $\mathsf P$
denoting the \proba{} and $\langle \rangle $ the expectation values, it
gives for every $a>0$:
\ekv{hs.10}
{ a\mathsf P(\sum_1^\infty \vert \alpha _j(\omega )\vert ^2\ge a)\le \langle
\sum_1^\infty \vert \alpha _j(\omega )\vert ^2\rangle =
\sum_1^\infty \langle \vert \alpha _j(\omega )\vert ^2\rangle  
=\sum_1^\infty \sigma _j^2.} 
We will prefer \no{hs.9} however, since it gives an exponential 
decay \wrt{} $a$. 
\end{remark}
\begin{remark}\label{hs2}\rm
If $Q=(\alpha _{j,k}(\omega ))_{j,k\in{\bf N}}$ is a random matrix where 
$\alpha _{j,k}(\omega )$
are \indep{} ${\cal N}(0,\sigma _{j,k}^2)$ laws, and 
\ekv{hs.11}
{
\sum_{j,k}\sigma _{j,k}^2 <\infty, 
}
then \no{hs.9} gives an estimate on the \proba{} that the Hilbert-Schmidt 
norm is  $\ge a^{1/2}$:
\ekv{hs.12}
{
\mathsf P(\Vert (\alpha _{j,k}(\omega ))\Vert _{{\rm HS}}^2\ge a)\le C(s_1)\exp [{C_0\over 
2s_1}\sum_{j,k\in{\bf N}^2}\sigma _{j,k}^2-{1\over 2s_1}a]
}
where $C_0$, $C(s_1)$
is the same constants as in \no{hs.9} and $s_1=\max \sigma _{j,k}^2$.  
\end{remark}

\section{Estimates on determinants of Gaussian random matrices}\label{er}
\setcounter{equation}{0}

Consider first a random vector
\ekv{er.1}
{\trans{u(\omega )}=(\alpha _1(\omega ),...,\alpha _N(\omega ))\in{\bf C}^N,
}
where $\alpha _1,...,\alpha _N$ are \indep{} complex Gaussian random variables 
with a ${\cal N}(0,1)$ law and $\omega $ is the random parameter living 
in a probability space with 
probability $\mathsf P$. The law of $\alpha _j$, i.e. the direct image of $\mathsf P$ under $\alpha _j$, 
is given by 
\ekv{er.2}
{
 (\alpha _j)_*(\mathsf P)={1\over \pi }e^{-\vert z\vert ^2}L(dz)=:f(z)L(dz)
}
and $L(dz)=L_{\bf C}(dz)$ is the Lebesgue measure on 
${\bf C}$. 

\par The distribution of $u$ is 
\ekv{er.4}
{
u_*(\mathsf P)={1\over \pi ^N}e^{-\vert u\vert ^2}L_{{\bf C}^N}(du).
}
If $U:{\bf C}^N\to {\bf C}^N$ is unitary, then $Uu$ has the same 
distribution as $u$. 

\par We next compute the distribution of $\vert u(\omega )\vert ^2$. The 
distribution of $\vert \alpha _j(\omega )\vert ^2$ is $\mu (r)dr$, where 
$$
\mu (r)=-H(r){d\over dr}e^{-r}=e^{-r}H(r),
$$
where $H(r)=1_{[0,\infty [}(r)$. We have $\widehat{\mu }(\rho )={1\over 
1+i\rho }$.

\par We have $\vert u(\omega )\vert ^2=\sum_1^N \vert \alpha _j(\omega 
)\vert ^2$ and since $\vert \alpha _j(\omega )\vert ^2$ are \indep{} and 
identically distributed, the distribution of $\vert u(\omega )\vert ^2$ is 
$\mu *...*\mu\, dr=\mu ^{*N}dr$, where $*$ indicates convolution. For 
$r>0$, we get by straight forward calculation the $\chi_{2N}^2$ distribution (for the variable $2r$)
\ekv{er.4.5}
{
\mu ^{*N}dr={r^{N-1}e^{-r}\over (N-1)!}H(r) dr.
}
Recall here that 
$$
\int_0^\infty r^{N-1}e^{-r}dr=\Gamma (N)=(N-1)!,
$$
so $\mu ^{*N}$ is indeed normalized. 

\par The expectation value of each $\vert \alpha _j(\omega )\vert ^2$ is $1$ so:
\ekv{er.5}
{\langle \vert u(\omega )\vert ^2\rangle =N.}

\par We next estimate the probability that $\vert u(\omega )\vert ^2$ is 
very large in a fashion that is slightly different from that of Section 
\ref{hs}. It will be convenient to pass to the variable $\ln (\vert 
u(\omega )\vert ^2)$, which has the distribution obtained from \no{er.4.5} 
by replacing $r$ by $t=\ln r$, so that $r=e^t$, $dr/r=dt$. Thus $\ln 
(\vert u(\omega )\vert ^2) $ has the distribution 
\ekv{er.6}
{
{r^Ne^{-r}\over (N-1)!}H(r){dr\over r}={e^{Nt-e^t}\over (N-1)!}dt=:\nu 
_N(t)dt.
}

%

\par Now consider a random matrix
\ekv{er.8}
{
(u_1...u_N)
}
where $u_k(\omega )$ are random vectors in ${\bf C}^N$ (here viewed as 
column vectors) of the form 
$$
\trans{u_k(\omega )}=(\alpha _{1,k}(\omega ),...,\alpha _{N,k}(\omega )),
$$
and all the $\alpha _{j,k}$ are \indep{} with the same law \no{er.2}. 

\par
Then
\ekv{er.10}
{
\det (u_1\,u_2...u_N)=\det (u_1\,\widetilde{u}_2...\widetilde{u}_N),
}
where $\widetilde{u}_j$ are obtained in the following way (assuming the $u_j$ 
to be linearly \indep{}, as they are almost surely): $\widetilde{u}_2$ is 
the \og{} projection of $u_2$ in the \og{} complement $(u_1)^{\perp}$, 
$\widetilde{u}_3$ is the \og{} projection of $u_3$ in 
$(u_1,u_2)^\perp=(u_1,\widetilde{u_2})^\perp$, etc.

\par If $u_1$ is fixed, then $\widetilde{u}_2$ can be viewed as a random 
vector in ${\bf C}^{N-1}$ of the type \no{er.1}, \no{er.2}, and with 
$u_1,u_2$ fixed, we can view $\widetilde{u}_3$ as a random vector of the same 
type in ${\bf C}^{N-2}$ etc. On the other hand
\ekv{er.9'}
{
\vert \det (u_1\, u_2...u_N)\vert ^2=\vert u_1\vert ^2\vert 
\widetilde{u}_2\vert ^2\cdot ..\cdot \vert \widetilde{u}_N\vert ^2.
} 
The squared lengths $\vert u_1\vert^2 , \vert 
\widetilde{u}_2\vert^2 ,...,\vert \widetilde{u}_N\vert^2 $ are \indep{} random 
variables with distributions $\mu ^{*N}dr, \mu ^{*(N-1)}dr,..., \mu dr$. 
This reduction plays an important role in \cite{Gi}.
The following lemma will not be used directly.
\begin{lemma}\label{er2}
Let $\alpha ,\beta >0$ be \indep{} random variables with distributions 
$\mu _\alpha (r){dr\over r}$, $\mu _\beta (r) {dr\over r}$. Then the 
product $\alpha \beta $ has the distribution $\mu _{\alpha \beta 
}{dr\over r}$, with 
\ekv{er.10'}
{
\mu _{\alpha \beta }=\mu _\alpha \sharp \mu _\beta :={\cal M}\inv (({\cal 
M}\mu _\alpha ) ({\cal M}\mu _\beta )).
}
Here 
$$
{\cal M}\mu (\tau )=\int r^{-i\tau }\mu (r){dr\over r}
$$
is the Mellin \tf{} of $\mu $.
\end{lemma}
\begin{proof} Recall that the Mellin \tf{} of $\mu (r)$ is the \F{} \tf{} 
of $\mu (e^t)$; $r=e^t$, $r\inv dr=dt$. The distribution of $\ln \alpha $ 
is related to that of $\alpha $ by the same change of variables $\mu 
_\alpha (r){dr\over r}\to \mu _\alpha (e^t)dt=\nu _\alpha (t)dt$. Since 
multiplication on the \F{} \tf{} side corresponds to convolution, 
\no{er.10'} 
is equivalent to the fact that the distribution of the sum of two \indep{} 
random variables is equal to the convolution of the distributions of the 
two variables.\end{proof}

\par The proof also shows that the multiplicative convolution in the lemma 
is given by
\ekv{er.11}
{
\mu _\alpha \sharp \mu _\beta (r)=\int_0^\infty  \mu _\alpha ({r\over 
\rho })\mu _\beta (\rho ){d\rho \over \rho }.}

\par As already mentioned we shall not use the lemma directly but rather its 
proof by taking logarithms and use that the distribution 
of the random variable $\ln \vert \det (u_1\, u_2...u_N)\vert ^2$ is equal 
to 
\ekv{er.12}
{
(\nu _1*\nu _2*...*\nu _N)dt,
} 
with $\nu _j$ defined in \no{er.6}.

\par We have 
$$
\nu _N(t)\le \widetilde{\nu }_N(t):={1\over (N-1)!}e^{Nt}.
$$
Choose $x(N)\in{\bf R}$ such that 
\ekv{er.13}
{
\int_{-\infty }^{x(N)}\widetilde{\nu }_N(t)dt=1.
}
More explicitely, we have 
\ekv{er.14}
{{1\over N!}e^{Nx(N)}=1,\quad x(N)={1\over N}\ln (N!)={1\over N}\ln 
\Gamma (N+1).}
Using Stirling's formula,
$$
{(N-1)!\over \sqrt{2\pi }}={\Gamma (N)\over \sqrt{2\pi 
}}=e^{-N}N^{N-{1\over 2}}(1+{\cal O}({1\over N})),
$$
 we get
\eeekv{er.15}
{x(N)&=&{1\over N}\big({1\over 2}\ln (2\pi )-(N+1)+(N+{1\over 2})\ln 
(N+1)+{\cal O}({1\over N})\big)}
{&=& {1\over N}\big( (N+{1\over 2})\ln N-N+C_0+{\cal O}({1\over N})\big)
}
{&=& \ln N+{1\over 2N}\ln N-1+{C_0\over N}+{\cal O}({1\over N^2}),
}
where $C_0=(\ln 2\pi )/2>0$.

\par With this choice of $x(N)$, we put 
$$
\rho _N(t)=1_{]-\infty ,x(N)]}(t)\widetilde{\nu }_N(t),
$$
so that $\rho _N(t)dt$ is a probability measure ``obtained from $\nu 
_N(t)dt$, by transfering mass to the left'' in the sense that 
\ekv{er.16}
{
\int f\nu _N dt\le \int f\rho _N dt,
}
whenever $f$ is a \bdd{} decreasing \fu{}. Equivalently,
$$
g*\nu _N\le g*\rho _N,
$$
whenever $g$ is a \bdd{} increasing \fu{}. Now, for such a $g$, both $g*\nu 
_N$ and $g*\rho _N$ are \bdd{} increasing functions, so by induction, we 
get 
$$
g*\nu _1*...*\nu _N\le g*\rho _1*...*\rho _N.
$$
In particular, by taking $g=H$, we get
\ekv{er.17}
{
\int_{-\infty }^t (\nu _1*...*\nu _N )(s)ds\le \int_{-\infty }^t 
(\rho _1*...*\rho _N)(s)ds, \ t\in {\bf R}.
} 

\par We have by \no{er.14}
\eekv{er.18}
{
\widehat{\rho }_N(\tau )&=&\int_{-\infty }^{x(N)}{1\over (N-1)!}e^{t(N-i\tau 
)}dt={1\over (N-1)!(N-i\tau )}e^{Nx(N)-ix(N)\tau}
}
{ &=&{e^{-ix(N)\tau }\over 1-i{\tau \over N}}.}
This \fu{} has a pole at $\tau =-iN$. 

\par Similarly,
\ekv{er.19}
{
\widehat{1_{]-\infty ,a]}}(\tau )={i\over \tau +i0}e^{-ia\tau }.
}
By Parseval's formula, we get 
\begin{eqnarray}\label{er.19bis}
\int_{-\infty }^a \rho _1*..*\rho _Ndt&=&{1\over 2\pi }\int_{-\infty 
}^\infty  {\cal F}(\rho _1*..*\rho _N)(\tau )\overline{{\cal 
F}1_{]-\infty ,a]}}(\tau )dt\\
&=&{1\over 2\pi }\int_{-\infty }^{+\infty }e^{-i\tau 
(\sum_1^Nx(j)-a)}{-i\over \tau -i0}\prod_1^N{1\over (1-{i\tau \over 
j})}d\tau .
\end{eqnarray}
We deform the contour to $\Im \tau =-1/2$ (half-way between ${\bf R}$ and 
the first pole in the lower half-plane).
For $j\ge 2$, we use the estimate
$$
\vert {1\over 1-{i\tau \over j}}\vert \le {1\over 1-{1\over 2j}}=\exp 
({1\over 2j}+{\cal O}({1\over j^2})),
$$
when $\Im \tau =-1/2$. Hence,
$$
\prod_2^N \vert {1\over 1-{i\tau \over j}}\vert\le \exp ({1\over 
2}\sum_2^N({1\over j}+{{\cal O}(1)\over j^2}))\le CN^{1\over 2}.
$$
It follows that for $a\le \sum_1^N x(j):$
\ekv{er.20}
{
\int_{-\infty }^a \rho _1*..*\rho _Ndt\le CN^{1\over 2}\exp (-{1\over 
2}(\sum_1^N x(j)-a)). 
}
In view of \no{er.17}, (\ref{er.19bis}) the \rhs{} is an upper bound for the probability 
that $\ln \vert \det (u_1...u_N)\vert ^2\le a$.

\par From the formula \no{er.15}, we get for some constants $C_1, C_2\in{\bf R}$:
\ekv{er.21}
{
\sum_1^N x(j)\ge C_1+(N+{1\over 2})\ln N-2N+{1\over 4}(\ln N)^2+C_0\ln N 
\ge C_2+(N+{1\over 2})\ln N -2N.
}
Hence, for $a\le C_2+(N+{1\over 2})\ln N-2N$,
\eeekv{er.22}
{&& \hskip -1truecm \mathsf P(\ln \vert \det (u_1...u_N)\vert ^2\le a)}
{&\le& CN^{1\over 2}\exp [-{1\over 2}(C_2+(N+{1\over 2})\ln N-2N-a)]}
{&=&C\exp [-{1\over 2}(C_2+(N-{1\over 2})\ln N-2N-a)].}

\par We shall next extend our bounds on the 
\renewcommand\proba{probability} \proba{} 
for the determinant to be small, to determinants of the form
$$
\det (D+Q)
$$
where $Q=(u_1...u_N)$ is as before, and $D=(d_1...d_N)$ is a fixed 
complex $N\times N$ matrix. As before, we can write 
$$
\vert \det ((d_1+u_1)...(d_N+u_N))\vert ^2=\vert d_1+u_1\vert ^2\vert 
\widetilde{d}_2+\widetilde{u}_2\vert ^2\cdot ..\cdot \vert 
\widetilde{d}_N+\widetilde{u}_N\vert ^2,
$$
where $\widetilde{d}_2=\widetilde{d}_2(u_1)$, 
$\widetilde{u}_2=\widetilde{u}_2(u_1,u_2)$ are the \og{} projections of 
$d_2$, $u_2$ on $(d_1+u_1)^\perp$, 
$\widetilde{d}_3=\widetilde{d}_3(u_1,u_2)$, 
$\widetilde{u}_3=\widetilde{u}_3(u_1,u_2,u_3)$ are the \og{} projections 
of $d_2$, $u_2$ on $(d_1+u_1,d_2+u_2)^{\perp}$ and so on.

\par Let $\nu _d^{(N)}(t)dt$ be the \proba{} distribution of $\ln \vert 
d+u\vert^2 $, when $d\in {\bf C}^N$ is fixed and $u\in{\bf C}^N$ is random 
as in \no{er.1}, \no{er.2}. Notice that $\nu _0^{(N)}(t)=\nu ^{(N)}(t)$ is the 
density we have alreay studied. 
\begin{lemma}\label{er3} For every $a\in{\bf R}$, we have 
$$\int_{-\infty }^a \nu _d^{(N)}(t)dt\le \int_{-\infty }^a \nu 
^{(N)}(t)dt.$$
\end{lemma} 
\begin{proof}
Equivalently, we have to show that $\mathsf P(\vert d+u\vert ^2\le 
\widetilde{a})\le \mathsf P(\vert u\vert ^2\le \widetilde{a})$ for every 
$\widetilde{a}>0$. For this, we may assume that $d=(c,0,...,0)$, $c>0$. We 
then only have to prove that 
$$
\mathsf P(\vert c+\Re u_1\vert ^2\le b^2)\le \mathsf P(\vert \Re u_1\vert ^2\le b^2),\ b>0,
$$
and here we may replace $\mathsf P$ by the corresponding \proba{} density
$$
\mu (t)dt={1\over {\sqrt{\pi }}}e^{-t^2}dt
$$
for $\Re \mu _1$. Thus, we have to show that
\ekv{er.23}
{
{1\over \sqrt{\pi }}\int_{\vert c+t\vert \le b}e^{-t^2}dt\le
{1\over \sqrt{\pi }}\int_{\vert t\vert \le b}e^{-t^2}dt .
}
Fix $b$ and rewrite the \lhs{} as 
$$
I(c)={1\over \sqrt{\pi }}\int_{-b-c}^{b-c}e^{-t^2}dt.
$$
The derivative satisfies (recall that $c>0$)
$$
I'(c)={1\over {\sqrt{\pi }}}(e^{-(b+c)^2}-e^{-(b-c)^2})\le 0.
$$
hence $c\mapsto I(c)$ is decreasing and \no{er.23} follows, since it is 
trivially fulfilled when $c=0$.
\end{proof}

\par Now consider the \proba{} that $\ln \vert \det (D+Q)\vert ^2\le a$. 
If $\chi _a(t)=H(a-t)$, this \proba{} becomes
\begin{eqnarray*}
&& \int ..\int \mathsf P(du_1)...\mathsf P(du_N)\times \\ &&\hskip -7truemm\chi _a(\ln \vert
d_1+u_1\vert ^2+\ln \vert
\widetilde{d}_2(u_1)+\widetilde{u}_2(u_1,u_2)\vert ^2+...+\ln \vert 
\widetilde{d}_N(u_1,..,u_{N-1})+\widetilde{u}_N(u_1,..,u_N)\vert ^2). 
\end{eqnarray*}
Here we first carry out the integration \wrt{} $u_N$, noticing that with 
the other $u_1,..,u_{N-1}$ fixed, we may consider 
$\widetilde{d}_N(u_1,..,u_{N-1})$
as a fixed vector in ${\bf C}\simeq (d_1+u_1,...,d_{N-1}+u_{N-1})^\perp$ 
and $\widetilde{u}_N$ as a random vector in ${\bf C}$. Using also the 
lemma, we get
\begin{eqnarray*}
&&\mathsf P(\ln \vert \det (D+Q)\vert ^2\le a)\\
&=&\int ..\int \nu 
_{\widetilde{d}_N}^{(1)}(t_N)dt_N\mathsf P(du_{N-1})..\mathsf P(du_1)\times \\
&&\chi _a(\ln \vert d_1+u_1\vert ^2+..+\ln \vert 
\widetilde{d}_{N-1}(u_1,..,u_{N-2})+\widetilde{u}_{N-1}(u_1,..,u_{N-1})\vert 
^2+t_N)\\
&\le& \int ..\int \nu^{(1)}(t_N)dt_N\mathsf P(du_{N-1})..\mathsf P(du_1)\times\\
&&\chi _a(\ln \vert d_1+u_1\vert ^2+..+\ln \vert 
\widetilde{d}_{N-1}(u_1,..,u_{N-2})+\widetilde{u}_{N-1}(u_1,..,u_{N-1})\vert 
^2+t_N).
\end{eqnarray*}
We next estimate the $u_{N-1}$- integral in the same way and so on. 
Eventually, we get 
\begin{prop}\label{er4}
Under the assumptions above,
\begin{eqnarray*}
\mathsf P(\ln \vert \det (D+Q)\vert ^2\le a)&\le& \int ..\int \chi 
_a(t_1+...+t_N)\nu ^{(1)}(t_N)\nu ^{(2)}(t_{N-1})..\nu ^{(N)}(t_1)\\
&=&\mathsf P(\ln \vert \det Q\vert ^2\le a).\end{eqnarray*} 
In particular the estimate \no{er.22} 
extends to random perturbations of constant matrices: 
\ekv{er.24}
{\mathsf P(\ln \vert \det (D+Q)\vert ^2\le a)\le 
C\exp [-{1\over 2}(C_2+(N-{1\over 2})\ln N-2N-a)],}
when $a\le C_2 +(N+{1\over 2})\ln N-2N$.
\end{prop}

\section{Grushin \pb{} for the perturbed \op{}}\label{gp}
\setcounter{equation}{0}

\par Let $P$
be as in Section \ref{gf}. Let $0<\widetilde{m},\widehat{m}\le 1$ be 
square integrable order \fu{}s on ${\bf R}^{2n}$
such that $\widetilde{m}$ or $\widehat{m}$ is integrable, and let 
$\widetilde{S}\in S(\widetilde{m})$, $\widehat{S}\in S(\widehat{m})$ be 
elliptic symbols. We use the same symbols to denote the $h$-Weyl 
quantizations. The \op{}s $\widetilde{S}$, $\widehat{S}$ will be \hs{} 
with 
$$
\Vert \widetilde{S}\Vert _{{\rm HS}}, \Vert \widehat{S}\Vert _{{\rm HS}}\backsim 
h^{-{n\over 2}}.
$$
Let $\widetilde{e}_1,\widetilde{e}_2,...$, and 
$\widehat{e}_1,\widehat{e}_2,...$ be \on{} bases for $L^2({\bf R}^n)$. Our 
random perturbation will be 
\ekv{gp.1}
{
Q_\omega =\widehat{S}\circ \sum_{j,k}\alpha _{j,k}(\omega 
)\widehat{e}_j\widetilde{e}_k^*\circ \widetilde{S},
}
where $\alpha _{j,k}$
are \indep{} complex ${\cal N}(0,1)$ \rv{}s. See the appendix, Section \ref{ap} 
for a general discussion.

\par Consider the polar decompositions
\ekv{gp.2}{\widehat{S}=\widehat{D}\widehat{U},\ 
\widetilde{S}=\widetilde{U}\widetilde{D},}
where 
$\widehat{U}$, $\widetilde{U}$ are unitary \pop{}s with symbol in $S(1)$ 
and $\widehat{D}$, $\widetilde{D}$ are positive \sa{} elliptic \pop{}s 
with symbol in $S(\widehat{m})$ and $S(\widetilde{m})$ respectively. After 
replacing $\widehat{e}_j$
by $\widehat{U}\widehat{e}_j$ and $\widetilde{e}_k$ by 
$\widetilde{U}^*\widetilde{e}_k$, we get with the new \on{} bases that 
\ekv{gp.3}
{
Q_\omega =\widehat{D}\circ \sum_{j,k}\alpha _{j,k}(\omega 
)\widehat{e}_j\widetilde{e}_k^*\circ \widetilde{D},
}
Now as we recall in the appendix (Section \ref{ap}), we may replace the 
bases $\widehat{e}_j$ and $\widetilde{e}_j$ by any new \on{} bases we like, 
if we replace the $\alpha _{j,k}(\omega )$ by a new set of \rv{}s (that we 
also denote by $\alpha _{j,k}$) having identical properties. If we 
choose $\widehat{e}_j$ to be an \on{} basis of \ef{}s of $\widehat{D}$ and 
similarly for $\widetilde{e}_j$, then we get 
\ekv{gp'.4}
{ Q_\omega =\sum_{j,k}\widehat{s}_j\widetilde{s}_k\alpha _{j,k}(\omega
)\widehat{e}_j\widetilde{e}_k^*, } where $\widehat{s}_j>0$ and
$\widetilde{s}_j>0$ are the \ev{}s of $\widehat{D}$ and $\widetilde{D}$ 
respectively, i.e. the singular values of $\widehat{S}$ and 
$\widetilde{S}$.

\par We are then precisely in the situation of Section \ref{hs}, noting 
that $\widehat{s}_j\widetilde{s}_k\alpha _{j,k}(\omega )$ are \indep{} 
${\cal N}(0,\widehat{s}_j^2\widetilde{s}_k^2)$-laws, so \no{hs.12} can be applied 
with $\sigma _{j,k}=\widehat{s}_j\widetilde{s}_k$, 
$$
\sum_{j,k}\sigma _{j,k}^2=\Vert \widehat{S}\Vert _{{\rm HS}}^2\Vert 
\widetilde{S}\Vert ^2_{{\rm HS}}\backsim h^{-2n}.
$$
We also know that 
$$
s_1=\max \sigma _{j,k}=\Vert \widehat{S}\Vert \Vert \widetilde{S}\Vert 
\backsim 1.
$$
From \no{hs.12}, we deduce that 
\ekv{gp'.5}
{
\mathsf P(\Vert Q_\omega \Vert _{{\rm HS}}^2\ge a)\le C\exp [Ch^{-2n}-{a\over C}]
}
for some constant $C>0$. Let 
\ekv{gp'.6}
{
M=C_1h^{-n},
}
for some $C_1\gg 1$. Then \no{gp'.5} gives
\ekv{gp'.7}
{
\mathsf P(\Vert Q\Vert _{{\rm HS}}^2\ge M^2)\le C\exp (-h^{-2n}/C),
}
for some new constant $C>0$.

\par We also want to control the trace class norm of 
$Q_\omega $, so we will use the assumption that one of  $\widetilde{m}$ 
and $\widehat{m}$ is integrable. Assume for instance that $\widehat{m}$ is 
integrable. Then $\widehat{m}^{1/2}$ is square integrable and we can 
factorize $\widehat{S}=\widehat{S}_1\widehat{S}_2$, with 
$\widehat{S}_j\in{\rm Op}(\widehat{m}^{1/2})$ being Hilbert-Schmidt operators. Let us write
$$
Q_\omega =\widehat{S}_1\widehat{S}_2\sum_{j,k}\alpha _{j,k}(\omega 
)\widehat{e}_j\widetilde{e}_k^*\widetilde{S}.
$$
Now recall that the composition of two \hs{} \op{}s is of trace class and 
the corresponding trace class norm does not exceed the product of the \hs{} 
norms of the two factors. Knowing that $\Vert \widehat{S}_1\Vert 
_{{\rm HS}}={\cal O}(h^{-n/2})$
and applying \no{gp'.7} to $\widehat{S}_2\sum_{j,k}\alpha _{j,k}(\omega 
)\widehat{e}_j\widetilde{e}_k^*\widetilde{S}$, we get
\ekv{gp'.8}
{
\mathsf P(\Vert Q_\omega \Vert _{\rm tr}\ge M^{3/2})\le C\exp (-h^{-2n}/C).
}  

\par In the following we will restrict the attention to $Q_\omega $'s with 
\ekv{gp.4}
{
\Vert Q_\omega \Vert _{{\rm HS}}\le M,\ \Vert Q_\omega \Vert _{\rm tr}\le M^{3/2},
}
and we have just seen that the \proba{} that this is the case is bounded 
from below by $1-Ce^{-h^{-2n}/C}$.

\par We wish to study the \ev{} distribution of
\ekv{gp.5}
{
P_\delta =P-\delta Q_\omega ,
}
when $\delta >0$ is \sufly{} small. (The minus sign is for notational 
convenience only.)

\par Recall from Section \ref{gf}, that for $z\in\widetilde{\Omega }$,
\ekv{gp.6}
{
P(z)=(\widetilde{P}-z)^{-1}(P-z)
}
is a trace class \pert{} of the identity. We now introduce 
\ekv{gp.7}
{
P_\delta (z)=(\widetilde{P}-z)\inv (P-\delta Q_\omega  -z)=P(z)-\delta 
(\widetilde{P}-z)\inv Q_\omega .
}
The Grushin \pb{} will be used to find lower bounds for $\vert \det P_\delta 
(z)\vert $. First we derive an upper bound: We have with $P_\delta 
(z)=P_\delta $, $P=P(z)$:
\eeekv{gp.8}
{&&\hskip -10truemm P_\delta ^*P_\delta }
{&=&P^*P-\delta (P ^*(\widetilde{P}-z)\inv 
Q_\omega +Q_\omega ^*(\widetilde{P}^*-\overline{z})\inv P -\delta 
Q_\omega ^*(\widetilde{P}^*-\overline{z})\inv (\widetilde{P}-z)\inv 
Q_\omega )}
{&=&P^*P+\delta R,
}
where 
\eekv{gp.9}
{
\Vert R\Vert_{{\rm HS}} &\le& C (\Vert Q_\omega \Vert _{{\rm HS}}+\delta \Vert
Q_\omega \Vert\Vert
Q_\omega \Vert_{{\rm HS}})\le \widetilde{C}M,
}
{
\Vert R\Vert_{{\rm tr}} &\le& C (\Vert Q_\omega \Vert _{{\rm tr}}+\delta \Vert
Q_\omega \Vert\Vert
Q_\omega \Vert_{{\rm tr}})\le \widetilde{C}M^{3/2},
}
provided that $\delta \Vert Q_\omega \Vert\le {\cal O}(1)$, as will
follow from \no{gp.12}.

\par In Section \ref{fc} we studied $P^*P+\alpha \chi (\alpha \inv P^*P)$ 
for $h\ll \alpha \ll 1$. This \op{} is 
$\ge \alpha $ if $1_{[0,1]}\le \chi $, as we may assume. Now assume that 
\ekv{gp.12}
{
\delta M\ll h .
}
Then
$$
P^*P+\alpha \chi (\alpha \inv P^*P)+\delta R\ge {\alpha \over 2},
$$
and 
\begin{eqnarray*}
\ln\det P_\delta ^*P_\delta &\le& \ln \det (P_\delta ^*P_\delta +\alpha 
\chi ({P^*P\over \alpha }))
\\&=& \ln \det (P^*P+\alpha \chi ({P^*P\over \alpha })+\delta R)
\\
&=&\ln \det (P^*P+\alpha \chi ({P^*P\over \alpha }))+\int_0^\delta \tr 
((P^*P+\alpha \chi ({P^*P\over \alpha })+tR)\inv R)dt.
\end{eqnarray*}
The integral is ${\cal O}(1){\delta \over \alpha }\Vert R\Vert _\tr ={\cal 
O}(1)\delta M^{3\over 2}/\alpha $ and combining this with \no{fc.51} 
(assuming now \no{fc.44}),  we get
\ekv{gp.13}
{
\ln \det P_\delta ^*P_\delta  \le {1\over (2\pi h)^n}(\iint \ln \vert 
p\vert ^2 dxd\xi +{\cal O}(1)\alpha ^{\kappa}\ln {1\over \alpha })+{\cal O}(1){\delta M^{3\over 2}\over \alpha }.
}
Here we choose $\alpha =Ch$ , $C\gg 1$ and we can drop the last remainder 
term if we assume that
\ekv{gp.14}
{
\delta M^{3\over 2}\ll h^{1+\kappa-n}\ln {1\over h},\ \delta \ll
h^{1+\kappa+n/2}\ln {1\over h}.
} 
For $n\geq 2$ this follows from \no{gp.12}, but for $n =1$ it might be a stronger assumption depending on the value of $\kappa$. 
Then
\ekv{gp.15}
{
\ln \vert \det P_\delta \vert \le {1\over (2\pi h)^n}(\iint \ln \vert 
p\vert dxd\xi +{\cal O}(1)h^{\kappa}\ln {1\over h}).
}

\par Still with $h\ll \alpha \ll 1$ we define $R_+$, $R_-$, ${\cal 
P}_0={\cal P}$, ${\cal E}_0={\cal E}$ as in Section \ref{gu}. Here $z$ is 
fixed, $P=P(z)$. With $P_\delta =P_\delta (z)$, we put 
\ekv{gp.16}
{
{\cal P}_\delta =\begin{pmatrix}
P_\delta &R_-\cr R_+ &0\end{pmatrix}:\, L^2({\bf 
R}^n)\times {\bf C}^N\to L^2({\bf R}^n)\times {\bf C}^N.
}
Now $\Vert \delta (\widetilde{P}-z)\inv Q_\omega \Vert \le C\delta M$, 
and 
$$
{\delta M\over \sqrt{\alpha }}\ll {h\over \sqrt{\alpha }}\ll 1
$$
under the assumption \no{gp.12}, so ${\cal P}_\delta $ has the inverse
\ekv{gp.17}
{
{\cal E}_\delta ={\cal E}_0\Big( 1-
\begin{pmatrix}\delta (\widetilde{P}-z)\inv 
Q_\omega &0\cr 0 &0\end{pmatrix}{\cal E}_0\Big)\inv
}
of norm $\le {\cal O}(1/\sqrt{\alpha })$. Writing 
\ekv{gp.18}
{
\widetilde{Q}_\omega =(\widetilde{P}-z)\inv Q_\omega ,
}
we have the Neumann series expansion
\eekv{gp.19}
{ {\cal E}_\delta &=&\begin{pmatrix}E^\delta &E_+^\delta \cr E_-^\delta
&E_{-+}^\delta \end{pmatrix} \\ \notag}{&=&
\begin{pmatrix}\sum_{j=0}^\infty  E^0(\delta \widetilde{Q}_\omega E^0)^j 
&\sum_{j=0}^\infty  
(E^0\delta 
\widetilde{Q}_\omega )^jE_+^0 \cr \sum_{j=0}^\infty E_-^0(\delta
\widetilde{Q}_\omega E^0)^j
&E_{-+}^0+\sum_{j=1}^\infty E_-^0(\delta \widetilde{Q}_\omega E^0)^{j-1}\delta 
\widetilde{Q}_\omega E_+^0\end{pmatrix}.
}

\par For $0\le t\le \delta $ we have 
\begin{eqnarray*}
{d\over dt}\ln \det {\cal P}_t &=&-\tr {\cal E}_t 
\begin{pmatrix}\widetilde{Q}_\omega &0\cr 0 &0\end{pmatrix}\\
&=&{\cal O}(1) {1\over \sqrt{\alpha }}\Vert \widetilde{Q}_\omega \Vert _\tr\\
&=&{\cal O}(1){1\over \sqrt{\alpha }}M^{3\over 2},
\end{eqnarray*}
so $$
\ln \det {\cal P}_\delta =\ln \det ({\cal P}) +{\cal O}(1){\delta \over 
\sqrt{\alpha }}M^{3\over 2}.
$$
Applying \no{gu.13}, we get
\ekv{gp.20}
{
\ln \vert \det {\cal P}_\delta \vert ={1\over (2\pi h)^n}(\iint\ln\vert 
p\vert dxd\xi +{\cal O}(1)
\alpha ^{\kappa}\ln \alpha 
)
+{\cal 
O}(1){\delta \over \sqrt{\alpha }}M^{3\over 2}.
}
Again, under the strengthened assumption \no{gp.14}, we get with
$\alpha =Ch$, $C\gg 1$,
\ekv{gp.21}
{
\ln \vert \det {\cal P}_\delta \vert ={1\over (2\pi h)^n}(\iint\ln\vert 
p\vert dxd\xi +{\cal O}(1)
h^{\kappa}\ln {1\over h}
).
}

\par The idea to get a lower bound for $\ln \vert \det P_\delta \vert $ with 
high \proba{} is now to use \no{dg.10}, \no{dg.11} which gives 
\ekv{gp.22}
{
\ln \vert \det P_\delta \vert =\ln \vert \det {\cal P}_\delta \vert +\ln 
\vert \det E_{-+}^\delta \vert ,
}
and to get a lower bound for $\ln \vert \det E_{-+}^\delta \vert $.

\section{Lower bounds on the determinant}\label{lb}
\setcounter{equation}{0}

We keep the assumptions 
formulated in the beginning of Section \ref{gp}, in particular \no{gp.1}. We 
restrict the attention to the case when \no{gp.4} holds with $M$ given by 
\no{gp'.6}, and recall that so is the case with probability 
$\ge 1-Ce^{-h^{-2n}/C}$. The restrictions \no{gp.12}, \no{gp.14} on $\delta $ 
will be further strengthened below.

\par Using a formula of the type \no{gp.22} we shall show that for every 
$z\in \widetilde{\Omega }$, the determinant of $P_\delta (z)$ is very
likely not to be too small.  For that we study the \proba{} distribution of the 
random matrix $E_{-+}^\delta $, and show that we are close enough to 
the Gaussian case to be able to apply the results of Section \ref{er} to the 
determinant. Recall that we work under the assumption \no{gp.4}, which is 
fulfilled with \proba{} $\ge 1-Ce^{-h^{-2n}/C}$. We want to study the map
\begin{eqnarray*}
Q\mapsto E_{-+}^\delta =E_{-+}^0+\sum_1^\infty E_-^0(\delta 
\widetilde{Q}E^0)^{j-1}\delta \widetilde{Q}E_+^0\\
=E_{-+}^0+\delta E_-^0\widetilde{Q}E_+^0+\sum_2^\infty 
\Big( {\delta Ch^{-n}\over 
\sqrt{\alpha }}\Big)^j{1\over \sqrt{\alpha }}R_j,
\end{eqnarray*}
where $\widetilde{Q}=(\widetilde{P}-z)\inv Q$ and $\Vert R_j\Vert _{{\rm HS}}\le 1$. 
Here, we used that $\Vert E_\pm ^0\Vert ,\, \Vert E^0\Vert \le 
1/\sqrt{\alpha }$. We can rewrite this further as 
\begin{eqnarray}\label{lb.1}
E_{-+}^\delta &=&E_{-+}^0+{\delta \over \alpha }(\sqrt{\alpha 
}E_-^0\widetilde{Q}\sqrt{\alpha }E_+^0+\frac{Ch^{-n}{\delta Ch^{-n}
}}{\al}\sum_0^\infty \Big( {\delta Ch^{-n}\over 
\sqrt{\alpha }}\Big)^j R_{j+2})\\
&=:&E_{-+}^0+{\delta \over \alpha }\widehat{Q}. \nonumber
\end{eqnarray} 
\par  We strengthen \no{gp.12}, \no{gp.14} to 
\ekv{lb.5}
{
\frac{\delta M^2}{\al} \ll 1,
}
and recall that by (\ref{gp'.6}), $M=C_{1}h^{-n}$. 

Then
$$
{\delta Ch^{-n}\over \sqrt{\alpha }}\ll {h^n\over C}\ll 1,\
$$ 
and we get
\ekv{lb.6}
{
\widehat{Q}=\sqrt{\alpha }E_-^0\widetilde{Q}\sqrt{\alpha }E_+^0+T,\ \Vert 
T\Vert _{{\rm HS}}\le \frac{C^2h^{-2n}\delta}{\al} \ll 1.
}

\par In view of \no{gp.1} we have 
\ekv{lb.7}
{
\sqrt{\alpha }E_-^0\widetilde{Q}\sqrt{\alpha }E_+^0=\sqrt{\alpha 
}E_-^0(\widetilde{P}-z)\inv 
\widehat{S}\sum_{j,k}\widehat{e}_j\alpha _{j,k}
\widetilde{e}_k^*\widetilde{S}\sqrt{\alpha 
}E_+^0,
}
where we recall from \no{gu.8} that
\ekv{lb.8}
{
\sqrt{\alpha }E_+^0v_+=\sum_1^N v_+(j)e_j,\quad \sqrt{\alpha 
}E_-^0v(j)=(v\vert f_j),\quad 1\le j\le N,
}
where $e_1,..,e_N$
and $f_1,...,f_N$ are \on{} bases for $\ran(1_{[0,\alpha 
]}(P(z)^*P(z)))$ and $\ran(1_{[0,\alpha 
]}(P(z)P(z)^*))$ respectively, writing $\ran(B)$ for the range of $B$.

\par Here, we wish to apply the discussion of Section \ref{ap}. The \op{}s 
$\widetilde{S}\sqrt{\alpha }E_+^0$, $\sqrt{\alpha 
}E_-^0(\widetilde{P}-z)\inv \widehat{S}$ are clearly \hs{} of rank 
$\le $ $N$. Let $\widetilde{t}_j$, $\widehat{t}_j$ denote the 
singular values of these \op{}s so that $\widetilde{t}_j=\widehat{t}_j=0$ 
for $j\ge N+1$.
\begin{lemma}\label{lba}
We have 
\ekv{lb.9}
{
{1\over C}\le \widetilde{t}_j,\, \widehat{t}_j\le C,\ 1\le j\le N,
}
where $C>0$
is \indep{} of $h,\alpha $.
\end{lemma}
\begin{proof}
\no{lb.8} shows that $\Vert \sqrt{\alpha }E_+^0\Vert ,\Vert \sqrt{\alpha 
}E_-^0\Vert \le 1$, and clearly $\Vert 
(\widetilde{P}-z)\inv\widehat{S}\Vert ,\Vert \widetilde{S}\Vert ={\cal 
O}(1)$, so the upper bound in \no{lb.9} is clear.

\par On the other hand, $\sqrt{\alpha }E_+^0v_+$ is confined to a \bdd{} 
region in phase space, and it is easy to show that 
$$
C\Vert \widetilde{S}\sqrt{\alpha }E_+^0v_+\Vert \ge \Vert 
\sqrt{\alpha }E_+^0v_+\Vert =\Vert v_+\Vert ,
$$
which implies that the smallest \ev{} of $((\widetilde{S}\sqrt{\alpha 
}E_+^0)^*(\widetilde{S}\sqrt{\alpha 
}E_+^0))^{1/2}$ is $\ge $ $1/C$. The lower bound on $\widetilde{t}_j$ 
follows. The argument for $\widehat{t}_j$ is essentially the same. 
\end{proof}

\par Let $\widehat{f}_1,...,\widehat{f}_N$ and 
$\widetilde{f}_1,...,\widetilde{f}_N$ be \on{} bases in ${\bf C}^N$
of \ef{}s of $((\sqrt{\alpha }E_-^0(\widetilde{P}-z)\inv 
\widehat{S})(\sqrt{\alpha }E_-^0(\widetilde{P}-z)\inv \widehat{S})^* 
)^{1/2}$ and $((\widetilde{S}\sqrt{\alpha }E_+^0)^*
(\widetilde{S}\sqrt{\alpha }E_+^0))^{1/2}$ respectively, with $\widehat{t}_j$ 
and $\widetilde{t}_k$ as the corresponding \ev{}s. We can then choose the 
\on{} bases $\{ \widehat{e}_j\}$, $\{\widetilde{e}_j\}$ in $L^2$ so that 
\ekv{lb.9.5}
{\widehat{e}_j={1\over \widehat{t}_j}(\sqrt{\alpha 
}E_-^0(\widetilde{P}-z)\inv \widehat{S})^*\widehat{f}_j,\quad 
\widetilde{e}_j={1\over \widetilde{t}_j}(\widetilde{S}\sqrt{\alpha 
}E_+^0)\widetilde{f}_j,
}
for $j=1,2,...,N$. Then from \no{lb.7}, we get
\ekv{lb.10}
{
\sqrt{\alpha }E_-^0\widetilde{Q}\sqrt{\alpha }E_+^0=\sum_{1\le j,k\le 
N}\widehat{t}_j\widetilde{t}_k\alpha _{j,k}\widehat{f}_j\widetilde{f}_k^*.
}

\par Now we will be a a little more specific about the assumption 
\no{gp.4}. We will restrict the attention to the set ${\cal Q}_M$ of 
matrices $(\alpha _{j,k}(\omega ))$ such that 
\ekv{lb.11}
{
\Vert \widehat{S}_2\sum \alpha _{j,k}(\omega 
)\widehat{e}_j\widetilde{e}_k^*\widetilde{S}\Vert _{{\rm HS}}\le M,
}
which implies \no{gp.4} and which is fulfilled with \proba{} $\ge$ 
$1-C\exp (-h^{-2n}/C)$. Here we recall that we assumed $\widehat{m}$ to be 
integrable and wrote $\widehat{S}=\widehat{S}_1\widehat{S}_2$ with 
$\widehat{S}_j\in{\rm Op}(S(\widehat{m}^{1/2}))$. (When $\widetilde{m}$ is 
integrable instead, we make a corresponding factorization of $\widetilde{S}$.) 

\par \no{lb.6} can be reformulated as 
\ekv{lb.12}
{
\widehat{Q}(\alpha )={\rm diag\,}(\widehat{t}_j)\circ \Big( (\alpha 
_{j,k})_{1\le j,k\le N}+\widetilde{T}(\alpha_{\cdot} )\Big) \circ {\rm 
diag\,}(\widetilde{t}_k),
} 
\ekv{lb.13}
{
\Vert \widetilde{T}(\alpha_{\cdot} )\Vert _{{\rm HS}}\le {\cal O}(1)\frac{\delta M^2}{\al},
}
for $(\alpha _{j,k})\in{\cal Q}_M$.

\par Let $\3 (\alpha _{j,k})\3  $ denote the norm in \no{lb.11} and let 
${\cal H}$ be the corresponding Hilbert space of ${\bf N}\times {\bf 
N}$ matrices. We shall view ${\rm HS}({\bf C}^N,{\bf C}^N)=:{\cal H}_N$ as a 
subspace of ${\cal H}$ in the natural way. Note that the two norms are 
\ufly{} equivalent on this subspace.

\par The Cauchy inequality implies (after decreasing $M$ by a constant 
factor) that
the differential of the map $\alpha_{\cdot} \mapsto \widetilde{T}(\alpha_{\cdot} )$ satisfies the 
following estimate on ${\cal Q}_M$:
\ekv{lb.14}
{
\Vert d\widetilde{T}\Vert _{{\cal H}\to{\cal H}_N}={\cal O}(1)\frac{\delta M}{\al}.
} 
On ${\cal H}$, ${\cal H}_N$
we have the basic \proba{} measures,
\ekv{lb.15}
{
\mu _{\cal H}=\prod_{j,k=1}^\infty (e^{-\vert \alpha _{j,k}\vert 
^2}{L(d\alpha _{j,k})\over \pi })\quad \mu _{{\cal H}_N}=
\prod_{j,k=1}^N (e^{-\vert \alpha _{j,k}\vert 
^2}{L(d\alpha _{j,k})\over \pi }).
}
We shall now estimate $\Pi _*(\mu _{\cal H})$ on ${\cal Q}_M$, where 
\ekv{lb.15.5}
{
\Pi ((\alpha_{j,k}) )=(\alpha _{j,k})_{1\le j,k\le 
N}+\widetilde{T}(\alpha_{\cdot} ),
}
and to do so, we identify $\widetilde{T}(\alpha_{\cdot} )$ with its image in 
${\cal H}$ under the natural inclusion ${\cal H}_N\subset {\cal H}$, and 
write 
\ekv{lb.16}
{
\Pi =\Pi _0\circ \kappa ,\quad \kappa (\alpha_{\cdot} )=\alpha_{\cdot} 
+\widetilde{T}(\alpha_{\cdot} ),\quad \Pi _0(\alpha_{\cdot} )=(\alpha _{j,k})_{1\le j,k\le N},
}
for $\alpha_{\cdot} =(\alpha _{j,k})\in{\cal H}$.

\par We first proceed formally, ignoring some technical difficulties due to 
the infinite dimension. We have
\eeeekv{lb.17}
{
\vert \, \Vert \kappa (\alpha_{\cdot} )\Vert _{{\rm HS}}^2-\Vert \alpha_{\cdot} \Vert 
_{{\rm HS}}^2\, \vert &=& \vert \, \Vert \Pi _0\kappa (\alpha_{\cdot} )\Vert _{{\rm HS}}^2-\Vert \Pi 
_0\alpha_{\cdot} \Vert _{{\rm HS}}^2\, \vert 
}
{
&=& \vert \, \Vert \Pi _0\kappa (\alpha_{\cdot} )\Vert _{{\rm HS}}-\Vert \Pi 
_0\alpha_{\cdot} \Vert _{{\rm HS}}\, \vert\times  \vert \, \Vert \Pi _0\kappa (\alpha_{\cdot} )\Vert _{{\rm HS}}+\Vert \Pi 
_0\alpha_{\cdot} \Vert _{{\rm HS}}\, \vert
}
{
&\le& \Vert \widetilde{T}(\alpha_{\cdot} )\Vert_{{\rm HS}} \big( 2\3 \alpha_{\cdot} \3+{\cal O}(\frac{\delta 
M^2}{\al})\big)
 }
{
&\le& {\cal O}(1)\frac{\delta M^3}{\al},
} 
where we strengthen the assumption \no{lb.5} to 
\ekv{lb.20}
{ \frac{{\delta }M^3}{\al}\ll 1,\hbox{ or equivalently }{\delta }\ll h^{3n+1/2}.  }   
As for the Jacobian of $\kappa $, we recall that
if $A:{\cal H}\to {\cal H}$ is linear with $\Vert A\Vert _{\rm tr}\ll 1$
(\ufly{} \wrt{} $M$), then $\det (1+A)=1+{\cal O}(\Vert A\Vert _{\tr})$. 
Also, if $A$ is of rank $\le N^2$, we know that $\Vert A\Vert _\tr \le
N^2\Vert A\Vert $, so in our case we get from \no{lb.14} that $$
\det {\partial \kappa \over \partial x}=1+{\cal 
O}(1)\frac{\delta N^2M}{\al}.
$$ 
Here the remainder term is $\ll 1$ in view of the assumption \no{lb.20} and 
the fact that $N\ll M$. (Recall that $N={\cal O}(\alpha ^{\kappa}h^{-n})$ by \no{fc.50}.) 

\par If $F$ is a locally defined \hol{} map $:{\cal H}\to {\cal H}$, 
then 
$$
L(dF(x))=\vert \det {\partial F\over \partial x}\vert ^2L(dx),
$$
so in our case,
$$
L(d\kappa (x))=(1+{\cal O}(1)\frac{\delta N^2M}{\al})L(dx).
$$ 

\par Combining this with \no{lb.17}, we get 
$$
\kappa _*(\mu _{\cal H})\le (1+{\cal O}(1)\frac{\delta M^3}{\al})\mu _{\cal H}\hbox{ 
on }{\cal Q}_M.
$$ 
Since $(\Pi _0)_*\mu _{\cal H}=\mu _{{\cal H}_N}$, we conclude that 
\ekv{lb.21}
{
\Pi _*(\mu _{\cal H})\le (1+{\cal O}(1)\frac{\delta M^3}{\al})\mu _{{\cal H}_N}\hbox{ 
on }{\cal Q}_M.
} 
At the end of this section we shall complete the proof of \no{lb.21} 
by means of finite dimensional approximations.

\par For $\alpha_{\cdot} =\alpha_{\cdot} (\omega )\in{\cal Q}_M$ we want to estimate the 
\proba{} that $\vert \det E_{-+}^\delta \vert $ is small. According to 
Proposition \ref{er4}, the $\mu _{{\cal H}_N}(d\check{Q})$-measure of the 
set of matrices $\check{Q}$ with 
$$
\vert \det ({\rm diag\,}(\widehat{t}_j)\inv\circ {\alpha \over \delta 
}E_{-+}^0\circ {\rm diag\,}(\widetilde{t}_j)\inv +\check{Q})\vert \le 
e^{a}
$$
is 
$$
\le Ce^{-{1\over 2}(C_2+(N-{1\over 2})\ln N-2N-a)},
$$
if 
\ekv{lb.22}
{ a\le C_2+(N+{1\over 2})\ln N-2N. } 
In view of \no{lb.1}, \no{lb.12},
\no{lb.21} this is also (after a slight increase of $C$) an upper bound for
the \proba{} to have $(\alpha _{j,k})\in{\cal Q}_M$ and $$
\vert \det ({\rm diag\,}(\widehat{t}_j)\inv\circ ({\alpha \over \delta 
}E_{-+}^0+\widehat{Q})\circ {\rm diag\,}(\widetilde{t}_j)\inv )\vert \le 
e^{a},
$$
or equivalently that 
$$
\vert \det E_{-+}^\delta \vert \le e^{{a}}({\delta \over \alpha 
})^N\prod_1^N\widehat{t}_j\prod_1^N\widetilde{t}_j.
$$

\par Write 
\ekv{lb.25}
{
a=C_2+(N-{1\over 2})\ln N-2N-b
}
and restrict the attention to $b\ge 0$. Then 
$$
e^{a}=e^{C_2+(N-{1\over 2})\ln N-2N-b}
$$
and we get
\ekv{lb.26}
{
\mathsf P(\hbox{\no{lb.11} holds and }\vert \det E_{-+}^\delta\vert  \le e^{N\ln{1\over 
\alpha }-N\ln{1\over\delta }+(N-{1\over 2})\ln N-CN+
C_2-b})\le e^{-b}.
}

\par Summing up the discussion so far, we have
\begin{prop}\label{lb1}
Consider the Grushin \pb{} \no{gp.16}. Assume \no{fc.44} 
and choose $\alpha =Ch$, $C\gg 0$. Then there exist positive constants 
$\widetilde{C}_0,\widetilde{C}_1,\widetilde{C}_2,\widetilde{C}$ such that for $b\ge 0$
\eekv{lb.28}
{
&& \mathsf P(\hbox{\no{lb.11} holds and }
\vert \det E_{-+}^\delta \vert 
\ge e^{-\widetilde{C}_0h^{\kappa-n}\ln {1\over h}-\widetilde{C}_1-
\widetilde{C}_2h^{\kappa-n}\ln{1\over \delta }-b})
} 
{
&&\ge 
1-\widetilde{C}e^{-b}-\widetilde{C}e^{-C_0h^{-2n}}.
}
Here $\delta $ is assumed to satisfy \no{lb.20}. 
 \end{prop}

\par In view of \no{gp.21}, \no{gp.22}, we get
\begin{theo}\label{lb2}
We now return to the original $P_\delta (z)$ in \no{gp.7} and we assume \no{fc.44}
\ufly{} for all $z$ in some open set $\widehat{\Omega }\Subset  
\widetilde{\Omega }$. If 
$\delta $ satisfies \no{gp.12},  \no{gp.14} there is a constant $C>0$ such that
\ekv{lb.29}
{
\ln \vert \det P_\delta \vert \le {1\over (2\pi h)^n}(\iint \ln \vert 
p\vert dxd\xi +Ch^{\kappa}\ln{1\over h}),\ \forall 
z\in\widehat{\Omega },
}
with \proba{} $\ge 1-Ce^{-C_0h^{-2n}}$. If $\delta $ satisfies the stronger condition  
\no{lb.20},  
 then there are constants $C, 
 \widetilde{C},C_0>0$ such that for every $z\in\widehat{\Omega }$ and 
 $\epsilon \ge 0$:
\ekv{lb.30}
{
\ln \vert \det P_\delta \vert \ge {1\over (2\pi h)^n}(\iint \ln \vert 
p\vert dxd\xi -Ch^{\kappa}(\ln{1\over h}+\ln{1\over \delta })-\epsilon )
}
with \proba{} $\ge 1-e^{-\epsilon (2\pi 
h)^{-n}}-\widetilde{C}e^{-C_0h^{-2n}}$.\end{theo}

\par Notice that the last term in the lower bound for the \proba{} is much 
smaller than the second term, and can therefore be eliminated. 
\par We end this section by completing the proof of \no{lb.21} by finite 
dimensional approximation. (We suggest the reader to proceed directly to 
Section \ref{sa}.)
\begin{lemma}\label{lb3}
We can choose the \on{} bases $\{\widehat{e}_j\}$, $\{\widetilde{e}_j\}$ 
in $L^2$ so that \no{lb.9.5} is fulfilled for $1\le j\le N$ and such that 
the square of the norm in \no{lb.11} is equivalent to 
\ekv{lb.31}
{
\sum_{j,k}\vert \alpha _{j,k}\vert ^2\widehat{\mu }_2(j)^2\widetilde{\mu 
}(k)^2,
}
where $\widehat{\mu }_2(j)$, $\widetilde{\mu }(k)$ denote the singular 
values of $\widehat{S}_2$ and $\widetilde{S}$ respectively.
\end{lemma}

\par In this lemma we did not try to have any \uf{}ity \wrt{} $h$.

\par Assume the lemma for a while. Then for $\widetilde{N}\ge N+1$, we can 
replace $Q_\omega $ in \no{gp.1} by 
$$
Q_\omega ^{\widetilde{N}}=\widehat{S}\circ \Big( \sum_{1\le j,k\le 
\widetilde{N}}\alpha _{j,k}(\omega )\widehat{e}_j\widetilde{e}_k^* 
+\sum_{j{\rm \, or\,}k\ge \widetilde{N}+1}\beta 
_{j,k}^{\widetilde{N}}(\alpha ^{\widetilde{N}}(\omega 
))\widehat{e}_j\widetilde{e}_k^*\Big)\circ 
\widetilde{S},
$$
which depends on finitely many \rv{}s. Here $\alpha 
^{\widetilde{N}}(\omega )=(\alpha _{j,k}(\omega ))_{1\le j,k\le\widetilde{N}}$ 
and $\beta ^{\widetilde{N}}_{j,k}$ are the linear \fu{}s of $\alpha 
^{\widetilde{N}}$ which minimize 
$$\Vert \widehat{S}_2\circ (\sum_{j,k\le 
\widetilde{N}}\alpha _{j,k}\widehat{e}_j\widetilde{e}_k^*+
\sum_{j\,{\rm or\,}k>\widetilde{N}}\beta _{j,k}\widehat{e}_j\widetilde{e}_k^*)\circ 
\widetilde{S}\Vert _{\rm HS}.$$ 
Here we can use the 
$\widehat{e}_j$, $\widetilde{e}_j$ of Lemma \ref{lb3}. On the set ${\cal 
Q}_M$, we have $\widehat{S}_2\circ Q^{\widetilde{N}}\circ \widetilde{S}\to 
\widehat{S}_2Q\widetilde{S}$ in \hs{} norm, and $\Vert
\widehat{S}_2Q^{\widetilde{N}}\widetilde{S}\Vert
_{\rm HS}\le M$ when $\alpha (\omega )\in{\cal Q}_M$.

\par We get the corresponding matrix $E_{-+}^{\delta 
,\widetilde{N}}=E_{-+}^0+{\delta \over \alpha 
}\widehat{Q}_{\widetilde{N}}$ and $\widehat{Q}_{\widetilde{N}}$ can be 
written as in \no{lb.12} with $\widetilde{T}(\alpha )$ replaced by 
$\widetilde{T}_{\widetilde{N}}(\alpha )$ satisfying \no{lb.13}. Now instead 
of $\mu _{\cal H}$ we have the finite dimensional measure 
$$
\mu _{{\cal H}_{\widetilde{N}}}=\prod_{j,k=1}^{\widetilde{N}}(e^{-\vert 
\alpha _{j,k}\vert ^2}{L(d\alpha _{j,k})\over \pi }),
$$
which we can view as the restriction of $\mu _{\cal H}$ to the tribe 
generated by $\alpha _{j,k}$ with $1\le j,k\le\widetilde{N}$) and we 
define $\Pi ^{\widetilde{N}}$
as in \no{lb.15.5} with $\widetilde{T}$ replaced by 
$\widetilde{T}_{\widetilde{N}}$. The subsequent arguments now become 
rigorous since we are in finite dimension and we get 
$$
\Pi ^{\widetilde{N}}_*(\mu _{{\cal 
H}_{\widetilde{N}}})\le (1+{\cal O}(1)\frac{\delta M^3}{\al})\mu _{{\cal H}_N}\hbox{ 
on }{\cal Q}_M.
$$ 
Since $\Pi ^{\widetilde{N}}\to \Pi $ on ${\cal Q}_M$, we obtain \no{lb.21} 
in the limit.

We next prove Lemma \ref{lb3}. 
\begin{proof} 
We consider first the following simplified problem. Let $m_S\le 1$ be a 
square integrable order \fu{} and let $S\in{\rm Op\,}(S(m_S))$ be elliptic. 
We look for an \on{} basis $e_1,e_2,..$ in $L^2$ such that 
\ekv{lb.32}
{
\Vert \sum u_kSe_k\Vert ^2\backsim \sum \mu _j(S)^2\vert u_k\vert ^2,
}
where $\mu _1(S)\ge \mu _2(S)\ge ...\to 0$ are the singular values of $S$ 
and such that 
$$
e_1,...,e_{N_0}
$$
 is a prescribed \on{} family of \fu{}s in ${\cal S}$.

\par Since $\sum \mu _j^2= {\cal O}(h^{-n})$, we have $N\mu _N^2={\cal 
O}(h^{-n})$, and using also that $\mu _N\le {\cal O}(1)$, we get
$$
\mu _N\le {{\cal O}(1)\over (Nh^n)^{1/2}+1}.
$$
On the other hand there exists a constant $\kappa _0>0$ such that 
$$m_S(\rho )\ge {1\over C_0}\langle \rho \rangle ^{-\kappa _0},$$
and we can use the mini-max principle to compare the \ev{}s of $(S^*S)^{1/2}$ 
with those of $(1+((hD)^2+x^2))^{-\kappa _0/2}$ and deduce that 
$$
\mu _N\ge {1\over {\cal O}(1)}{1\over (1+hN^{1/n})^{\kappa _0/2}}.
$$

\par If 0$<\mu \ll 1$, we have, with $p$ denoting the symbol of 
$(S^*S)^{1/2}$, that
\begin{eqnarray*}
{\rm dist\,}(0,p^{-1}([0,2\mu ]))&\ge &{\rm dist\,}(0,m_S^{-1}([0,{2\mu 
\over C}]))\\
&\ge & {\rm dist\,}(0,\{ \rho ;\, {1\over C_0}\langle \rho \rangle 
^{-\kappa _0}\le {2\mu \over C}\} )\\
&\ge &{1\over C_1}\mu ^{-1/\kappa _0}.
\end{eqnarray*}

\par If $u$ is a corresponding normalized \ef{}, we have $(\mu ^{-1}(S^*S)^{1/2}-1 
)u=0$, and we notice that $\mu ^{-1}(S^*S)^{1/2}\in{\rm Op\,}(S(\mu 
^{-1}m_S))$, where $\mu ^{-1}m_S$ satisfies \ufly{} the axioms of an 
order \fu{}, when $\mu \to 0$. We conclude that 
$$
u={\cal O}(1),\hbox{ in }H(m)
$$ 
\ufly{} \wrt{} $m$ if $m=m_\mu $ belongs to a family of order\fu{}s that 
satisfy \ufly{} the axioms and $m=1$ on $\{ \rho \in T^*{\bf R}^n;\, 
p(\rho )\le 2\mu \}$. From this, we deduce that 
$$
(\phi \vert u)={\cal O}(\mu ^N),\ \forall N,
$$ 
if $\phi \in {\cal S}$ is fixed, and $\mu \to 0$.

\par Let $f_1,f_2,...$ be an \on{} basis of \ef{}s of $(S^*S)^{1/2}$ with 
$\mu _1\ge \mu _2\ge ... $ the corresponding decreasing enumeration of \ev{}s. 
Then $(e_j\vert f_k)={\cal O}(k^{-\infty })$, $1\le j\le N_0$, $k\ge 1$. 
Let $N\gg N_0$. For $j\ge N+1$, put 
\ekv{lb.32a}
{
g_j=f_j-\Pi _{E_{N_0}}f_j=f_j+r_j,\ E_{N_0}\ni r_j={\cal O}(j^{-\infty }).
}
Here, we let $E_{N_0}=(e_1,...,e_{N_0})$
be the span of $e_1,...,e_{N_0}$ and $\Pi _{E_{N_0}}:L^2\to E_{N_0}$ be the 
corresponding \og{} projection. Then $g_j\in E_{N_0}^\perp$ and for 
$j,k>N$:
\ekv{lb.32b}
{
(g_j\vert g_k)=\delta _{j,k}+{\cal O}(j^{-\infty }k^{-\infty }).
}
Here the estimates are \uf{} \wrt{} $N$ and if $N$
is \sufly{} large, we see that $G=((g_j\vert g_k))$ is a positive definite 
matrix of which any real power has elements satisfying \no{lb.32b}. 
Let $(a_{j,k})=G^{-1/2}$, so that $a_{j,k}=\delta _{j,k}+{\cal 
O}(j^{-\infty }k^{-\infty })$, $j,k>N$. Put 
$$
e_j=\sum_{k>N}a_{j,k}g_k,\ j>N.
$$
Then $e_j$, $j>N$
form an \on{} basis in the span $G_N^\perp$ of $g_{N+1},g_{N+2},...$. We 
see that for $j>N$:
\ekv{lb.32c}
{
e_j=f_j+\sum_{k>N}{\cal O}(j^{-\infty }k^{-\infty })f_k+\widetilde{r}_j,\ 
E_{N_0}\ni\widetilde{r}_j={\cal O}(j^{-\infty }). 
}

\par $G_N=(G_N^\perp )^\perp$ is a space of dimension $N$, containing 
$E_{N_0}$. For $1\le j\le N$, we consider
$$
\Pi _{G_N}f_j=f_j-\Pi _{G_N^\perp}f_j,
$$
\begin{eqnarray*}
\Pi _{G_N^\perp}f_j&=&\sum_{N+1}^\infty  (f_j\vert e_k)e_k\\
&=& \sum_{N+1}^\infty  (f_j\vert \widetilde{r}_k)e_k\\
&=&\sum_{N+1}^\infty {\cal O}(k^{-\infty })(f_k+\sum_{\ell =N+1}^\infty 
{\cal O}(k^{-\infty }\ell^{-\infty })f_{\ell}+\widetilde{r}_k )\\
&=& \sum_{N+1}^\infty {\cal O}(k^{-\infty })f_k+{\cal O}(N^{-\infty }),
\end{eqnarray*}
where the last term is in $E_{N_0}$. Thus, we get for $1\le j\le N$:
$$
\Pi _{G_N}f_j=f_j+\sum_{N+1}^\infty  {\cal O}(k^{-\infty 
})f_k+\widehat{r}_j,\ E_{N_0}\ni \widehat{r}_j={\cal O}(N^{-\infty }).
$$
This implies
\ekv{lb.32d}
{
\Pi _{G_N}f_j=f_j+\sum_1^\infty {\cal O}((k+N)^{-\infty })f_k,\ 1\le j\le 
N.
}

\par Now, complete $e_1,...,e_{N_0}$ to an \on{} basis $e_1,...,e_N$ in 
$G_N$. Then $e_1,e_2,.....$ is an \on{} basis in $L^2$. \no{lb.32d} shows 
that $\Pi _{G_N}f_1,...,\Pi _{G_N}f_N$ is very close to being an \on{} 
basis in $G_N$, and we see that 
\ekv{lb.32e}
{
e_j=\sum_{k=1}^N u_{j,k}f_k+\sum_1^\infty  {\cal O}((k+N)^{-\infty })f_k,\ 
1\le j\le N,
}
where $(u_{j,k})_{1\le j,k\le N}$ is a unitary matrix. \no{lb.32c}, 
\no{lb.32e} imply that for all $j\ge 1$:
\ekv{lb.33}
{
e_j=\sum_{k=1}^\infty a_{j,k}f_k+\sum_{k=1}^\infty {\cal O}((j+N)^{-\infty 
}(k+N)^{-\infty })f_k,
}
where $a_{j,k}=u_{j,k}$ for $j,k\le N$ and $a_{j,k}=\delta _{j,k}$ when 
$\max (j,k)>N$.
\par We now fix $N$ \sufly{} large so that the above estimates hold.  Using
\no{lb.33}, we get
$$
e_j=f_j+\sum_k {\cal O}(j^{-\infty }k^{-\infty })f_k,\ 
Se_j=\mu _jf_j+\sum_k {\cal O}(j^{-\infty }k^{-\infty })f_k,
$$
$$
(Se_k\vert Se_j)=\mu _k^2\delta _{j,k}+{\cal O}(k^{-\infty }j^{-\infty }),
$$
\begin{eqnarray*}
\Vert \sum_1^\infty  u_k Se_k\Vert ^2&=&\sum_{k=1}^\infty \mu 
_k^2u_k\overline{u}_k+\sum_{k=1}^\infty \sum_{j=1}^\infty {\cal 
O}(k^{-\infty }j^{-\infty })u_k\overline{u}_j\\
&=& ((\mu ^2+K)u\vert u)_{\ell^2}\\
&=& ((1+\mu \inv K\mu \inv )\mu u\vert \mu u)
\end{eqnarray*}
where $\mu $ denotes the \op{} ${\rm diag\,}(\mu _j)$. Here $\mu \inv K\mu 
\inv$ is compact: $\ell^2\to \ell^2$, so 
$1+\mu \inv K\mu \inv$ is a non-negative \sa{} Fredholm \op{} of index 0. 
If $u(j)={\cal O}(j^{M_0})$ and $(1+\mu \inv K\mu \inv )u=0$, then $u={\cal
O}(j^{-\infty })$.  If $0\ne v\in\ell^2$, $(1+\mu \inv K\mu \inv )v=0$, we
conclude that $$ ((1+\mu \inv K\mu \inv )v\vert v)=((\mu ^2+K)u\vert
u)=\Vert S(\sum_1^\infty u_ke_k)\Vert^2>0 $$ thus $1+\mu \inv K\mu \inv$ is
bijective and we finally conclude that
\no{lb.32} holds.

\par Now we can finish the proof of the lemma. We choose 
$\{\widehat{e}_j\}$, $\{\widetilde{e}_j\}$ in \no{lb.9.5} so that
\no{lb.32} holds with $S=\widehat{S}_2$, $e_j=\widehat{e}_j$ and
$S=\widetilde{S}^*$, $e_j=\widetilde{e}_j$ respectively.  Then the square
of the norm \no{lb.11} is equal to
\ekv{lb.34}
{
\sum_{i,j,k,\ell }(\widehat{S}_2\widehat{e}_i\vert 
\widehat{S}_2\widehat{e}_j)(\widetilde{S}^*\widetilde{e}_k\vert 
\widetilde{S}^*\widetilde{e}_\ell )\alpha _{i,k}\overline{\alpha 
}_{j,\ell}=(\widehat{{\cal S}}\otimes \widetilde{{\cal S}}\alpha \vert 
\alpha )_{\ell^2\otimes \ell^2},
}
where
$$
\widehat{{\cal S}}_{j,i}=(\widehat{S}_2\widehat{e}_i\vert 
\widehat{S}_2\widehat{e}_j),\ \widetilde{{\cal S}}_{\ell 
,k}=(\widetilde{S}^*\widetilde{e_k}\vert 
\widetilde{S}^*\widetilde{e}_\ell ).
$$
From \no{lb.32} we know that 
\begin{eqnarray*}
&&\widehat{{\cal S}}=\widehat{\mu }\widehat{P}\widehat{P}\widehat{\mu },\ 
\widehat{\mu }={\rm diag\,}(\widehat{\mu }_2(j))\\
&&\widetilde{{\cal S}}=\widetilde{\mu }
\widetilde{P}\widetilde{P}\widetilde{\mu },\ 
\widetilde{\mu }={\rm diag\,}(\widetilde{\mu }(j)),
\end{eqnarray*}
where $\widehat{P}$, $\widetilde{P}$ are positive \sa{} \op{}s satisfying
$$
{1\over C}I\le \widehat{P},\,\widetilde{P}\le CI.
$$
Then \no{lb.34} can be written 
$$
\Vert (\widehat{P}\widehat{\mu }\otimes \widetilde{P}\widetilde{\mu 
})\alpha \Vert _{\ell^2\otimes \ell^2}^2,
$$
and the lemma follows.
\end{proof}

\section{Spectral \asy{}s when $dp,\, d\overline{p}$ are \indep{}}\label{sa}
\setcounter{equation}{0}

\par Let $\Gamma \Subset \Omega $ be open with $C^2$
\bdy{} and assume that for every $z\in \partial \Gamma $:
\eekv{sa.1}
{&&\Sigma _z:=p\inv (z)\hbox{ is a smooth sub-\mfld{} of }T^*{\bf R}^n\hbox{ 
on}}
{&&\hbox{which }dp,d\overline{p}\hbox{ are linearly \indep{} 
at every point.}}
This assumption, which is satisfied also in a \neigh{} of $\partial
\Gamma $, 
implies that ${\rm codim\,}(\Sigma _z)=2$. The assumption 
can also be rephrased more briefly by saying that $\partial\Gamma$ 
does not contain 
any critical value of $p:{\mathbb R}^{2n}\to {\mathbb R}^2$. Here $p$ is the leading symbol of the 
original ($z$-\indep{} \op{}.) If $p_z(\rho )=(\widetilde{p}(\rho )-z)\inv 
(p(\rho )-z)$ is the principal symbol of $(\widetilde{P}-z)\inv (P-z)$, we 
introduce 
\ekv{sa.2}
{I(z)=\int_{{\bf R}^2}\ln \vert p_z(\rho )\vert d\rho }
which is the same integral as in \no{lb.29}, \no{lb.30} (where $z$ was fixed). 
It is easy to see that $I(z)$ is a smooth \fu{} on the \neigh{} of
$\partial \Gamma $ where \no{sa.1} holds and as in 
\cite{MeSj} we can compute $\Delta _zI(z)$. Since $z\mapsto p_z(\rho )$ is 
\hol{}, we know that $\Delta _z\ln \vert p_z(\rho )\vert =0$ when 
$p_z(\rho )\ne 0$, ie when $\rho \not\in \Sigma _z$. On the other hand 
$p_z(\rho )=(\widetilde{p}(\rho )-z)\inv (p(\rho )-z)$ where the first 
factor is \hol{} in $z$ and non-vanishing, so
$$
\Delta _z\ln \vert p_z(\rho )\vert =\Delta _z\ln \vert p(\rho )-z\vert 
=2\pi \delta (z-p(\rho )).
$$
If $\phi \in C_0^\infty (\Omega )$, we get
\begin{eqnarray*}
&&\int (\Delta _zI(z))\phi (z)L(dz)=\int\int \Delta _z(\ln \vert p_z(\rho 
)\vert )\phi (z)L(dz)d\rho \\
&&=2\pi \int\int \delta (z-p(\rho ))\phi (z)L(dz)d\rho =2\pi \int \phi 
(p(\rho ))d\rho .
\end{eqnarray*}
Thus we get (as in \cite{Ha1, Ha2} when $n=1$):
\ekv{sa.3}
{
{1\over 2\pi }\Delta (I(z))L(dz)=p_*(d\rho )\hbox{ near }\partial
\Gamma ,
}
where $d\rho $ is the symplectic volume element. Notice that this formula 
is still true without the assumption \no{sa.1} and hence not only in a 
\neigh{} of $\partial \Gamma $, but in $\Omega$; however $I(z)$ is no more smooth in general 
but still well-defined as a distribution. 
This fact will be used in the proof of Theorem \ref{sa1}.

\par In view of \no{sa.1}, we have $V(t)\backsim t$ and \no{fc.44} 
holds \ufly{} with $\kappa=1$, when $z$ varies in 
a \neigh{} of $\partial \Gamma $. Correspondingly 
the conclusions in Theorem \ref{lb2} hold \ufly, when $z$ varies in a
small \neigh{} of $\partial \Gamma $.

\begin{theo}\label{sa1}
Let $\Gamma \Subset  \Omega $ be open with $C^2$ \bdy{} and make the 
assumption \no{sa.1}. Let $\delta >0$ satisfy \no{lb.20} 
 and
assume that $h\ln {1\over \delta }\ll \epsilon \ll 1$ (or equivalently 
$\delta \ge e^{-\epsilon /(Ch)}$, $C\gg 1$, $\epsilon \ll 1$, implying 
also that $\epsilon \ge \widetilde{C}h\ln {1\over h}$ for some 
$\widetilde{C}>0$).
Then 
with $C>0$ large enough, the number $N(P_\delta ,\Gamma )$ of \ev{}s of 
$P_\delta $ in $\Gamma $ satisfies 
\ekv{sa.4}
{
\vert N(P_\delta ,\Gamma )-{1\over (2\pi h)^n}{\rm vol\,}(p^{-1}(\Gamma 
))\vert \le C{\sqrt{\epsilon }\over h^n}
}
with \proba{} 
$$\ge 1 -{C\over \sqrt{\epsilon }}e^{-{\epsilon/2 \over (2\pi 
h)^n}}.
$$
\end{theo}
\begin{proof}
The \ev{}s of $P_\delta $ in $\widetilde{\Omega }$ coincide with the zeros 
of the \hol{} \fu{}
\ekv{sa.5}
{F_\delta (z)=\det P_\delta (z).}
Theorem \ref{lb2} tells us that there exists a \neigh{} $\widehat{\Omega }$ 
of $\partial \Omega $ such that 
\smallskip
\par\noindent (a) With \proba{} $\ge 1-Ce^{-C_0h^{-2n}}$, we have
$$
\ln \vert F_\delta (z)\vert \le {1\over (2\pi h)^n}(I(z)+Ch\ln{1\over h}),\ 
z\in\widehat{\Omega }.
$$
\smallskip
\par\noindent (b) For every $z\in\widehat{\Omega }$ and $\epsilon >0$ we have 
$$
\ln \vert F_\delta (z)\vert \ge {1\over (2\pi h)^n}(I(z)-Ch(\ln{1\over 
h}+\ln{1\over \delta })-\epsilon ),
$$
with \proba{} $\ge 1-e^{-\epsilon (2\pi h)^{-n}}-Ce^{-C_0h^{-2n}}$. Notice 
here that $\ln {1\over \delta }\ge \ln{1\over h}$.\smallskip

\par We can then repeat the arguments of \cite{Ha1, Ha2}. Recall Proposition 
6.1 from \cite{Ha2} proved in a more general form in \cite{Ha1}.
\begin{prop}\label{sa1.5} Let $\widehat{\Omega }$, $\widetilde{\Omega }$ 
be open \neigh{}s of $\partial \Gamma $ and $\overline{\Gamma }$ 
respectively.
Let $\phi \in C^\infty (\widehat{\Omega };{\bf R})$ and let $f$ be a 
\hol{} \fu{} in $\widetilde{\Omega }$ such that 
\ekv{sa.6}
{
\vert f(z;\widetilde{h})\vert \le e^{\phi (z)/\widetilde{h}},\ 
z\in\widehat{\Omega },\, 0<\widetilde{h}\ll 1.
}
Assume that for some $\epsilon >0$, $\epsilon \ll 1$, $\exists 
z_k\in\widehat{\Omega }$, $k\in J$, such that 
$$
\partial \Gamma \subset\bigcup_{k\in J}D(z_k,\sqrt{\epsilon }),\ \# J={\cal 
O}({1\over \sqrt{\epsilon }}),
$$
\ekv{sa.7}
{
\vert f(z_k;\widetilde{h})\vert \ge e^{{1\over \widetilde{h}}(\phi 
(z_k)-\epsilon )},\ k\in J.
}
Then the number of zeros of $f$ in $\Gamma $ satisfies
$$
\# (f\inv (0)\cap \Gamma )={1\over 2\pi \widetilde{h}}\int_\Gamma  \Delta 
\phi L(dz)+{\cal O}({\sqrt{\epsilon }\over \widetilde{h}}),
$$
where we let $\phi $
denote some distribution in ${\cal D}'(\Gamma \cup\widehat{\Omega })$ 
extending the previous function $\phi $.\end{prop}

\par The original statement in \cite{Ha1, Ha2} was with a smooth \fu{} $\phi $
defined in a whole \neigh{} of $\overline{\Gamma }$ satisfying \no{sa.6} 
there, but the proof works without any changes under the weaker 
assumptions above.

\par In view of (a), (b), we can apply the proposition with 
$\widetilde{h}=(2\pi h)^n$ and $\epsilon $ replaced by $2\epsilon $, 
$\phi =I(z)+Ch\ln {1\over 
h}$, $f=F_\delta $. Then \no{sa.6} holds with a \proba{} as in (a), while 
\no{sa.7} holds with a \proba{} 
$$
\ge 1-{C\over \sqrt{\epsilon }}e^{-{\epsilon \over 2}(2\pi h)^{-n}}
-Ce^{-C_0h^{-2n}}.
$$
We can define $\phi $ as a distribution in a full \neigh{} of 
$\overline{\Gamma }$ by \no{sa.2}. Then 
$$
{1\over 2\pi \widetilde{h}}\int_\Gamma \Delta \phi L(dz)={1\over (2\pi 
h)^n}\int_\Gamma {1\over 2\pi }\Delta I(z)L(dz)={1\over (2\pi 
h)^n}\iint_{p\inv (\Gamma )}dxd\xi .
$$
The theorem follows.
\end{proof}

\par We next give a result about the simultaneous Weyl \asy{}s for a 
\fy{} of domains
\begin{theo}\label{sa2}
Let ${\cal G}$ be a \fy{} of domains $\Gamma \Subset  \Omega $ that 
satisfy the assumptions of Theorem \ref{sa1} \ufly{} in the following sense:
Each $\Gamma $ is of the form $g(z)<0$ (with $g=g_\Gamma $) where $g$ 
belongs to a \bdd{} set in $C^2(\overline{\Omega })$ and $g>1/C$ on 
$\partial \Omega $ and $\vert dg\vert >1/C$ on $\partial \Gamma $, where 
$C>0$ is \indep{} of $\Gamma $. We also assume that \no{sa.1} holds for 
all $z\in\partial \Gamma $, $\Gamma \in{\cal G}$, \ufly{} \wrt{} 
$(z,\Gamma )$. 

\par Choose $\delta ,\epsilon $ as in Theorem \ref{sa1}. Then with \proba{}
$$\ge 1 -{{C}\over \epsilon }e^{-{\epsilon /2\over (2\pi
h)^n}},  $$
we have \no{sa.4} with a constant $C$ \indep{} of $\Gamma $.
\end{theo}
\begin{proof}
As in the proof of Theorem \ref{sa1}, we use Proposition \ref{sa2} with an 
appropriate grid of points $z_k$ (see \cite{Ha2} for further details). We 
now need ${\cal O}(1/\epsilon )$ points to achieve that the union of the 
$D(z_k,\sqrt{\epsilon })$ covers the union of all the $\partial \Gamma $, 
rather than ${\cal O}(1/\sqrt{\epsilon })$ points as in the proof of Theorem 
\ref{sa1}. 
\end{proof}

\section{Counting zeros of holomorphic \fu{}s} \label{cz}
\setcounter{equation}{0}
\par Let $\Gamma \Subset {\bf C}$ have smooth \bdy{} $\partial\Gamma $. Assume for
simplicity that $\gamma :=\partial \Gamma $ is connected (or
equivalently that $\Gamma $ is simply connected). This is for
notational convenience only. For $0<r\ll 1$, we put 
\ekv{cz.1}
{\gamma _r=\gamma +D(0,r)=\partial \Gamma +D(0,r).}
Then $\gamma _r$ has smooth \bdy{} and is a thin domain of width
$\approx 2r$. Let $G_r(z,w)$, $P_r(z,w)$ denote the Green and
Poisson kernels of $\gamma _r$, so that the Dirichlet \pb{}
$$
\Delta u=v,\ {u_\vert }_{\partial \gamma _r}=f,\quad u,v\in
C^\infty (\overline{\gamma _r}),\ f\in C^\infty (\partial
\gamma _r),$$ has the unique solution 
$$
u(z)=\int_{\gamma _r}G_r(z,w)v(w)L(dw)+\int_{\partial \gamma
_r}P_r(z,w)f(w)\vert dw\vert. 
$$

\par We recall some properties of the Green kernel: If $\Omega
\Subset  {\bf C}$ has a smooth \bdy{} and $G_\Omega (x,y) $is
the corresponding Green kernel, then 
\ekv{cz.2}{
G_\Omega \le 0
,}
\ekv{cz.3}{
G_\Omega \hbox{ is }C^\infty \hbox{ for }x\ne y,
}
\ekv{cz.4}
{
G_\Omega (x,y)={1+o(1)\over 2\pi }\ln \vert x-y\vert \hbox{ for
}x\approx y,\ x,y\not\in \partial \Omega ,\ x-y\to 0. }
\ekv{cz.5}
{
G_\Omega ({x\over r},{y\over r})=G_{r\Omega }(x,y),\ x,y\in
r\Omega. }
$\Omega ={1\over r}\gamma _r$ is a very long domain of
approximately constant width and \no{cz.4} is valid \ufly{} for
$x,y\in \Omega $, $\vert x-y\vert \le {\cal O}(1)$, ${\rm
dist\,}(x,\partial \Omega ),{\rm dist\,}(y,\partial \Omega )\ge
1/{\cal O}(1)$. Moreover,
\ekv{cz.6}
{
\vert G_{r\inv \gamma _r}(x,y)\vert \le C_0 e^{-\vert x-y\vert
/C_0},\ x,y\in r\inv \gamma _r,\, \vert x-y\vert \ge {1\over {\cal
O}(1)}. }

\par To recover these well-known facts, notice that $r\inv \gamma _r$ is 
given by $-1<\phi (x)<1$, where $\phi (x)$ is the suitably signed distance 
from $r\inv \partial \Gamma $ to $x$, so that $\vert \nabla \phi (x) \vert 
=1$, $\vert \nabla ^2\phi (x)\vert ={\cal O}(r)$. If $u\in H_0^1(r\inv 
\gamma _r)$, we have by integration by parts,
$$
\int_{r\inv\gamma _r}((\nabla \phi )^2+\phi (x)\Delta \phi )\vert u\vert 
^2dx=-2\Re \int_{r\inv \gamma _r}\phi (\nabla \phi \cdot \nabla 
u)\overline{u}dx,
$$
implying
$$
\int (1-{\cal O}(r))\vert u\vert ^2dx\le (2+{\cal O}(r)) \Vert \nabla 
u\Vert \Vert u\Vert ,
$$
$$
\Vert u\Vert \le (2+{\cal O}(r))\Vert \nabla u\Vert ,
$$
$$
-\Delta \ge ({1\over 4}-{\cal O}(r)),
$$
where $\Delta =\Delta _{r\inv \gamma _r}$ is the Dirichlet Laplacian on 
$r\inv \gamma _r$. From this estimate we can develop exponential decay 
estimates for $-\Delta $, since we still have a positive lower bound for 
$\Re (e^\psi (-\Delta )e^{-\psi })=-\Delta -\vert \nabla \psi \vert ^2$, 
if $\vert \nabla \psi (x)\vert^2\le 1/5 $. We drop the ensuing routine arguments.  
\par 
In view of \no{cz.5}, \no{cz.6} we get
\ekv{cz.7}
{\vert G_r(x,y)\vert \le C_0e^{-{1\over C_0r}\vert x-y\vert },\
x,y\in \gamma _r,\ \vert x-y\vert \ge r/{\cal O}(1),}
\eekv{cz.8}
{&&
\vert G_r(x,y)\vert ={1+o(1)\over 2\pi }\ln \vert {x\over r}-{y\over r}\vert ,\hbox{
for }\vert x-y\vert \le r/C_0, }
{&&
{\rm dist\,}(x,\partial \Omega ),{\rm
dist\,} (y,\partial \Omega )\ge r/C_0,\ \vert {x\over r}-{y\over
r}\vert \to 0 .} 

\par Let $\phi $ be a continuous subharmonic \fu{} defined in some \neigh{} of 
$\overline{\gamma _r}$. Let 
\ekv{cz.9}
{
\mu =\mu _\phi =\Delta \phi 
}
be the corresponding locally finite  positive measure.

Let $u$ be a \hol{} \fu{} defined in a \neigh{} of $\overline{\gamma _r}$. 
We assume that 
\ekv{cz.10}
{
h\ln \vert u(z)\vert \le \phi (z),\ z\in\overline{\gamma _r}. 
}

\begin{lemma}\label{cz1}
Let $C_1,C_2>1$ and let $z_0\in\overline{\gamma _{(1-{1\over C_1})r}}$ be 
a point where 
\ekv{cz.11}
{h\ln \vert u(z_0)\vert \ge \phi (z_0)-\epsilon ,\ 0<\epsilon \ll 1.}
Then the number of zeros of $u$ in $D(z_0,C_2r)\cap \gamma _{(1-{1\over 
C_2})r}$ is 
\ekv{cz.12}{\le {C_3\over h}(\epsilon +\int_{\gamma _r}-G_r(z_0,w)\mu (dw)),}
where $C_3$ is independent of $\epsilon ,h$.
\end{lemma}
\begin{proof}
Writing $\phi $ as a uniform limit of an increasing sequence of smooth 
\fu{}s, we may assume that $\phi \in C^\infty $.
Let 
$$
n_u(dz)=\sum 2\pi \delta (z-z_j),
$$
where $z_j$ are the zeros of $u$ counted with their multiplicity. We may 
assume that no $z_j$ are situated on $\partial \gamma _r$. Then, since 
$\Delta \ln \vert u\vert =n_u$,
\eeekv{cz.13}
{
h\ln \vert u(z)\vert &=& \int_{\gamma _r} G_r(z,w)h n_u (dw)+\int_{\partial 
\gamma _r}P_r(z,w)h\ln \vert u(w)\vert \vert dw\vert}   
{&\le& \int_{\gamma _r}G_r(z,w)hn_u(dw)+\int_{\partial \gamma 
_r}P_r(z,w)\phi (w)\vert dw\vert} 
{&=& \int_{\gamma _r}G_r(z,w)hn_u(dw)+\phi (z)-\int_{\gamma 
_r}G_r(z,w)\mu (dw).}
Putting $z=z_0$ in \no{cz.13} and using \no{cz.11}, we get
$$
\int_{\gamma _r}-G_r(z_0,w)hn_u(dw)\le \epsilon +\int_{\gamma 
_r}-G_r(z_0,w)\mu (dw).
$$
Now $$
-G_r(z_0,w)\ge {1 \over 2\pi  C_3},\ C_3>0,
$$
in $D(z_0,C_2r)\cap \gamma _{(1-{1\over C_2})r}$, and we get \no{cz.12}.
\end{proof}

\par Notice that this argument is basically the same as when using Jensen's 
formula to estimate the number of zeros of a \hol{} \fu{} in a disc. We could 
assume the bound $h\ln \vert u(z)\vert \le \phi (z)$, in 
$D(z_0,\widetilde{C}_2r)\cap \gamma _{(1-{1\over 
\widetilde{C}_2})r}=:\widetilde{\Omega }_r$ for some 
$\widetilde{C}_2>C_2$. Then we can replace the bound \no{cz.12} by 
$$
{C_3\over h}(\epsilon +\int_{\widetilde{\Omega }_r}-G_{\widetilde{\Omega 
}_r}(z_0,w)\mu (dw)),
$$
which is sharper, since $-G_{\Omega _1}\le -G_{\Omega _2}$, when $\Omega 
_1\subset \Omega _2$.

\par Now we sharpen the assumption \no{cz.11} and assume 
\ekv{cz.14}
{
h\ln \vert u(z_j)\vert \ge \phi (z_j)-\epsilon ,
}
where $z_1,...,z_N\in \gamma _{(1-{1\over C_1})r}$
are points such that 
\ekv{cz.15}
{
\gamma _r\subset \bigcup_1^N D(z_j,C_1r),\quad N\backsim {1\over r}.
}
We may assume that $z_1,z_2,...,z_N$ are arranged in such a way that 
\ekv{cz.16}
{
\vert z_j-z_k\vert \backsim r{\rm dist\,}(j,k),\ j\ne k,
}
where $j,k$ are viewed as elements of ${\bf Z}/N{\bf Z}$ and we take the 
natural distance on that set. We will also assume for a while that $\phi $ 
is smooth. 

\par According to Lemma \ref{cz1}, we have 
\ekv{cz.17}
{
\# (u\inv (0)\cap (D(z_j,C_1r)\cap \gamma _{(1-{1\over C_1})r}))\le 
{C_3\over h}(\epsilon +\int_{\gamma _r}-G_r(z_j,w)\mu (dw)).
}

\par We introduce 
\ekv{cz.18}
{
\widetilde{r}=(1-{1\over C_1})r
}
and consider the harmonic \fu{}s on $\gamma _{\widetilde{r}}$,
\ekv{cz.19}
{
\Psi (z)=h(\ln \vert u(z)\vert +\int_{\gamma 
_{\widetilde{r}}}-G_{\widetilde{r}}(z,w)n_u(dw)),
}
\ekv{cz.20}
{
\Phi (z)=\phi (z)+\int_{\gamma _{\widetilde{r}}}-G_{\widetilde{r}}(z,w)\mu 
(dw).
}
Then $\Phi (z)\ge \phi (z)$ with equality on $\partial \gamma 
_{\widetilde{r}}$. Similarly, $\Psi (z)\ge h\ln \vert u(z)\vert $ with 
equality on $\partial \gamma _{\widetilde{r}}$.

\par Consider the harmonic \fu{}
\ekv{cz.21}
{
H(z)=\Phi (z)-\Psi (z),\ z\in \gamma _{\widetilde{r}}.
}
Then on $\partial \gamma _{\widetilde{r}}$, we have by \no{cz.10} that
$$
H(z)=\phi (z)-h\ln \vert u(z)\vert \ge 0,
$$
so by the maximum principle,
\ekv{cz.22}
{
H(z)\ge 0,\hbox{ on }\gamma _{\widetilde{r}}.
}
By \no{cz.14}, we have 
\eeekv{cz.23}
{
H(z_j)&=& \Phi (z_j)-\Psi (z_j)
}
{
&=&\phi (z_j)-h\ln \vert u(z_j)\vert +
\int_{\gamma _{\widetilde{r}}}-G_{\widetilde{r}}(z_j,w)\mu 
(dw)-\int_{\gamma _{\widetilde{r}}} -G_{\widetilde{r}}(z_j,w)hn_u(dw)
}
{
&\le & \epsilon +\int_{\gamma _{\widetilde{r}}} -G_{\widetilde{r}}(z_j,w)\mu (dw).
}

\par Harnack's inequality implies that 
\ekv{cz.24}
{
H(z)\le {\cal O}(1)(\epsilon +\int -G_{\widetilde{r}}(z_j,w)\mu (dw))\hbox{ 
on }D(z_j,C_1r)\cap \gamma _{(1-{1\over C_1})\widetilde{r}}.
}

\par Now assume that $u$ extends to a \hol{} \fu{} in a \neigh{} of 
$\Gamma \cup \overline{\gamma _r}$. We then would like to evaluate the 
number of zeros of $u$ in $\Gamma $. Using \no{cz.17}, we first have 
\ekv{cz.25}
{
\# (u^{-1}(0)\cap \gamma _{\widetilde{r}})\le {C\over h}(N\epsilon 
+\sum_{j=1}^N \int_{\gamma _r}(-G_r(z_j,w))\mu (dw)).
}

\par Let $\chi \in C_0^\infty (\Gamma \cup \gamma _{(1-{1\over 
C_1})\widetilde{r}};[0,1])$ be equal to 1 on $\Gamma $. Of course $\chi $ 
will have to depend on $r$ but we may assume that for all $k \in{\bf N}$, and as $r\to 0$,
\ekv{cz.26}
{
\nabla ^k\chi ={\cal O}(r^{-k}).
}
We are interested in 
\ekv{cz.27}
{
\int \chi (z)hn_u(dz)=\int_{\gamma _{\widehat{r}}}h\ln \vert u(z)\vert 
\Delta \chi (z)L(dz),\ \widehat{r}=(1-{1\over C_1})\widetilde{r}.
}
Here we have on $\gamma _{\widetilde{r}}$
\eeeekv{cz.28}
{h\ln \vert u(z)\vert &=&\Psi (z)-\int_{\gamma 
_{\widetilde{r}}}-G_{\widetilde{r}}(z,w)hn_u(dw)}
{&=&\Phi (z)-H(z)-\int_{\gamma 
_{\widetilde{r}}}-G_{\widetilde{r}}(z,w)hn_u(dw)}
{
&=&\phi (z)+\int_{\gamma _{\widetilde{r}}}-G_{\widetilde{r}}(z,w)\mu 
(dw)-H(z)-\int_{\gamma _{\widetilde{r}}}-G_{\widetilde{r}}(z,w)hn_u(dw)
}
{&=& \phi (z)+R(z),}
where the last equality defines $R(z)$.

\par Inserting this in \no{cz.27}, we get 
\ekv{cz.29}
{
\int\chi (z)hn_u(dz)=\int \chi (z)\mu (dz)+\int R(z)\Delta \chi (z)L(dz).
}
(Here we also used some extension of $\phi $ to $\Gamma $ 
with $\mu =\Delta \phi $.) The task is now to estimate $R(z)$ and the 
corresponding integral in \no{cz.29}. Put 
\ekv{cz.30}
{
M_j=\mu (\Omega _j),\ \Omega _j=D(z_j,C_1r)\cap \gamma _r.
}
Using the exponential decay property \no{cz.7} (equally valid for 
$G_{\widetilde{r}}$) we get for $z\in \Omega _j\cap \gamma 
_{\widetilde{r}}$, ${\rm dist\,}(z,\partial (D(z_j,C_1r)\cap \gamma 
_{\widetilde{r}}))\ge r/{\cal O}(1)$:
\ekv{cz.31}
{
\int_{\gamma _{\widetilde{r}}}-G_{\widetilde{r}}(z,w)\mu (dw)\le 
\int_{\Omega _j\cap \gamma _{\widetilde{r}}}-G_{\widetilde{r}}(z,w)\mu 
(dw)+{\cal O}(1)\sum_{k\ne j}M_ke^{-{1\over C_0}\vert j-k\vert }.
}
Similarly from \no{cz.24}, we get 
\ekv{cz.32}
{
H(z)\le {\cal O}(1)(\epsilon +\int_{\Omega _j\cap \gamma 
_{\widetilde{r}}}-G_{\widetilde{r}}(z_j,w)\mu (dw)+\sum_{k\ne j}e^{-{1\over 
C_0}\vert j-k\vert }M_k),
}
for $z\in \Omega _j\cap \gamma _{\widetilde{r}}$. 

\par This gives the following estimate on the contribution from the first 
two terms in $R(z)$ to the last integral in \no{cz.29}:
\eeekv{cz.33}
{&&\int_{\gamma _{\widetilde{r}}}(\int_{\gamma
_{\widetilde{r}}}-G_{\widetilde{r}}(z,w)\mu (dw)-H(z))\Delta \chi (z)L(dz)
} { &=&{\cal O}(1)(N\epsilon +\sum_j (\sup_{z\in \Omega _j\cap \gamma
_{\widehat{r}}}\int_{\Omega _j\cap \gamma _{\widetilde{r}}}
-G_{\widetilde{r}}(z,w)\mu (dw)+\sum_{k\ne j}e^{-{1\over C_0}\vert j-k\vert
}M_k)) } 
{ &=& {\cal O}(1)(N\epsilon +\sum_j \sup_{z\in \Omega _j\cap
\gamma _{\widehat{r}}} \int_{\Omega _j\cap \gamma _{\widetilde{r}}}
-G_{\widetilde{r}}(z,w)\mu (dw)+\mu (\gamma
_r)).  }

\par The contribution from the last term in $R(z)$ (in \no{cz.28}) to the 
last integral in \no{cz.29} is 
\ekv{cz.34}
{
\int_{z\in \gamma _{\widehat{r}}}\int_{w\in \gamma 
_{\widetilde{r}}}G_{\widetilde{r}}(z,w)hn_u(dw)\Delta \chi (z)L(dz).
}
Here 
\begin{eqnarray*}
&& \int_{z\in \gamma _{\widehat{r}}}G_{\widetilde{r}}(z,w)(\Delta \chi )
(z)L(dz)\\
&=& \int_{\widetilde{z}\in \widetilde{r}\inv \gamma _{\widehat{r}}} 
G_{\widetilde{r}}(\widetilde{r}\widetilde{z},\widetilde{r}\widetilde{w})\Delta 
_z\chi (\widetilde{r}\widetilde{z})\widetilde{r}^2L(d\widetilde{z})\\
&=& \int G_{\widetilde{r}\inv \gamma 
_{\widetilde{r}}}(\widetilde{z},\widetilde{w})\Delta _{\widetilde{z}}(\chi 
(\widetilde{r}\widetilde{z}))L(dz)={\cal O}(1),
\end{eqnarray*}
so the expression \no{cz.34} is 
\eeekv{cz.35}
{
&&{\cal O}(h)\# (u\inv (0)\cap \gamma _{\widetilde{r}})
}
{
&=&{\cal O}(1)({\epsilon \over r}+\sum_{j=1}^N \int_{\gamma 
_r}(-G_r(z_j,w))\mu (dw))
}
{&=&
{\cal O}(1)({\epsilon \over r}+\sum_{j=1}^N \int_{\Omega 
_j}-G_r(z_j,w)\mu (dw)+\mu (\gamma _r)).
}
Using all this in \no{cz.29}, we get
\eekv{cz.36}
{&&\hskip -5truemm
\int\chi (z)hn_u(dz)=\int \chi (z)\mu (dz)
}
{
&&+{\cal O}(1)({\epsilon \over r}+\sum_j(\sup_{z\in 
\Omega _j\cap \gamma _{\widehat{r}}}\int_{\Omega _j\cap \gamma
_{\widetilde{r}}} -G_{\widetilde{r}}(z,w)\mu (dw)+
\int_{\overline{\Omega _j}}-G_r(z_j,w)\mu (dw))+\mu (\gamma _r)).
}
We replace the smoothness assumption on $\phi $ by the assumption that 
$\phi $ is continuous near $\Gamma $ and keep \no{cz.14}. Then by 
regularization, we still get \no{cz.36}.

\par In order to simplify this further, we introduce a weak regularity 
assumption on the measure $\mu $. Assume first that $\mu =\Delta \phi $ is 
defined in a fixed $r$-\indep{} neighborhood of $\partial \Gamma $. For 
$D(z,t)$ contained in that \neigh{} we assume that as $t\to 0$, 
\ekv{cz.37}
{
W_z(t):=\mu (D(z,t))={\cal O}(t^{\rho _0}),
}
for some $0<\rho _0\le 2$. 
\begin{remark}\label{cz1.5}\rm
It is easy to see that this assumption on $\Delta \phi $ implies that $\phi $ 
is continuous near $\Gamma $. In the case $\rho _0>1$, we notice that as $r\to 0$,
\ekv{cz.38}
{
\mu (\gamma _r)={\cal O}(r^{\rho _0-1}).
}
(This is true also for $\rho _0\le 1$ but then of no interest.)
\end{remark}
\begin{lemma}\label{cz2}
Assume \no{cz.37} for some $\rho _0\in ]0,2]$. Then for every domain 
$\Omega \subset \gamma _r$ and every $z\in \Omega \cap \gamma _{(1-{1\over 
C})r}$, we have for $0<t\le r/2$:
\ekv{cz.39}
{
\int_{\Omega }-G_r(z,w)\mu (dw)\le {\cal O}(1)t^{\rho _0}\ln {r\over 
t}+{\cal O}(1)\ln({r\over t})\mu (\Omega ).
}
\end{lemma}
\begin{proof}
Write 
$$
\int_\Omega -G_r(z,w)\mu (dw)=\int_{D(z,t)\cap\Omega }-G_r(z,w)\mu 
(dw)+\int_{\Omega \setminus D(z,t)}-G_r(z,w)\mu (dw).
$$
For $\vert z-w\vert \ge t$, we have $-G_r(z,w)\le {\cal O}(1)\ln {r\over t}$ 
(cf \no{cz.5}), so the last integral is ${\cal O}(1)\ln ({r\over t})\mu 
(\Omega )$. For $w\in D(z,t)\cap\Omega $, we have 
$$
-G_r(z,w)\le {\cal O}(1)\ln {r\over \vert z-w\vert },
$$
hence
\begin{eqnarray*}
\int_{D(z,t)\cap \Omega }-G_r(z,w)\mu (dw)&\le &{\cal O}(1)\int_0^t \ln {r\over 
s}dW_z(s)\\
&=&{\cal O}(1)([\ln ({r\over s})W_z(s)]_0^t+\int_0^t {1\over s}W_z(s)ds)\\
&=&{\cal O}(1) t^{\rho _0}\ln {r\over t}.
\end{eqnarray*}
\end{proof}

\begin{coro}\label{cz3}
Under the same assumptions, we have for every $N\in{\bf N}$:
\ekv{cz.40}
{
\int_\Omega -G_r(z,w)\mu (dw)\le {\cal O}_N(1)(r^N+\ln ({1\over r})\mu 
(\Omega )).
}
\end{coro}
\begin{proof}
We just choose $t=r^M$, $0<M\in{\bf N}$ and use that $\ln r^{-M}=M\ln r\inv$
\end{proof}

If we assume \no{cz.38}, then the corollary allows us to simplify \no{cz.36} 
to 
\ekv{cz.41}
{
\int \chi (z)hn_u(dz)=\int \chi (z)\mu (dz)+{\cal O}(1){\epsilon \over 
r}+{\cal O}_N(r^N+\ln({1\over r})\mu (\gamma _r)).
}

\par Summing up the discussion, we have proved
\begin{prop}\label{cz4}
Let $\Gamma \Subset {\bf C}$ have smooth \bdy{} and let $\phi $ be a 
continuous subharmonic \fu{} defined near $\overline{\Gamma }$. 
Then we have the following 
result, valid \ufly{} for $0<\epsilon \ll 1$, $0<r\ll 1$, $0<h\ll 1$: 
Let $u$ be a \hol{} \fu{}, defined in $\Gamma +D(0,r)$ with $h\ln \vert 
u(z)\vert \le \phi (z)$, $z\in\partial \Gamma +D(0,r)$
and assume that there exist $z_1,...,z_N\in\partial \Gamma +D(0,{r\over 
2})$ such that 
\ekv{cz.41.5}
{\partial \Gamma +D(0,r)\subset \bigcup_1^N
D(z_j,2r),\ N\backsim {1\over r},\ h\ln \vert u(z_j)\vert \ge \phi 
(z_j)-\epsilon .}
Then with $\Omega _j=D(z_j,2r)\cap (\partial \Gamma +D(0,r))$, we 
have
\eekv{cz.42}
{
&&\vert \# (u\inv (0)\cap \Gamma )-{1\over 2\pi h}\int_\Gamma \Delta \phi 
L(dz)\vert \le {{\cal O}(1)\over h}\Big( {\epsilon \over r}+\mu (\gamma _r) 
}
{&&
+\sum_j (\sup_{z\in\Omega _j\cap 
(\partial \Gamma +D(0,{r\over 4}))}\int_{\Omega _j\cap (\partial \Gamma 
+D(0,{r\over 2}))} -G_{r\over 2}(z,w)\mu 
(dw)+\int_{\overline{\Omega } _j}-G_r(z_j,w)\mu (dw)\Big) .
}

\par If we assume also that \no{cz.37} holds for some $0<\rho _0\le 2$, 
then we have for every $N>0$:
\eekv{cz.43}
{
&&\vert \# (u\inv (0)\cap \Gamma )-{1\over 2\pi h}\int_\Gamma \Delta \phi 
L(dz)\vert \le
}
{&&
{{\cal O}(1)\over h}({\epsilon \over r}+{\cal O}_N(1)(r^N+\ln({1\over 
r})\mu (\partial \Gamma +D(0,r)))).
}\end{prop}
\begin{ex}\label{cz5}\rm
If $\phi $ is of class $C^2$ near the \bdy{}, then \no{cz.37} is satisfied
with $\rho _0=2$ and $\mu (\partial \Gamma +D(0,r))={\cal O}(r)$. We 
choose $N=1$
so that the \rhs{} of \no{cz.43} becomes
$$
{{\cal O}(1)\over h}({\epsilon \over r}+r\ln{1\over r}).
$$
If we choose $r=\sqrt{\epsilon }$, we get
$$
\vert \# (u\inv (0)\cap\Gamma )-{1\over 2\pi h}\int_\Gamma \Delta \phi 
(z)L(dz)\vert \le {{\cal O}(1)\over h}\sqrt{\epsilon }\ln{1\over \epsilon }.
$$
In this case we loose a factor $\ln \epsilon \inv$ compared to Proposition 
6.1 in 
\cite{Ha2}.
\end{ex}

\section{Spectral \asy{}s in a more general case} \label{gc}
\setcounter{equation}{0}

\par Let $\Gamma \Subset \widetilde{\Omega }$ be open with 
$C^\infty $ \bdy{}. For $z$ in a \neigh{} of $\partial \Gamma $ and 
$0<s,t\ll 1$, we put
\ekv{gc.1}
{
V_z(t)={\rm Vol\,}\{\rho \in{\bf R}^{2n};\, \vert p(\rho )-z\vert ^2\le 
t\},\ W_z(s)=V_z(s^2).
}

\par Recall that in any \bdd{} domain in phase space, the symbols $\vert 
p_z(\rho )\vert ^2=q_z(\rho )$ and $\vert p(\rho )-z\vert ^2$ are \ufly{} 
of the same order of magnitude. If we replace $\vert p(\rho )-z\vert ^2$ 
by $q_z(\rho )$ in \no{gc.1}, we get a new \fu{} $V_z^{\rm new}(t)$ such 
that 
\ekv{gc.2}
{
V_z({t\over C})\le V_z^{\rm new}(t)\le V_z(Ct).
}
For the purposes of this paper, we can therefore identify the two \fu{}s 
and resort to the second definition whenever we find it convenient. Also, 
when $z$ is fixed, we will sometimes write $V(t)$ instead of $V_z(t)$.

\par Our weak assumption, replacing \no{sa.1}, is
\ekv{gc.3}
{
\exists \kappa\in ]0,1], \hbox{ such that }V_z(t)={\cal O}(t^{\kappa}),\hbox{ \ufly{} for }z\in{\rm neigh\,}(\partial \Gamma ),\ 0\le t\ll 1.
} 

\begin{ex}\label{gc0}\rm
When \no{sa.1} holds for $z\in{\rm neigh\,}(\partial \Gamma )$, it is clear 
that \no{gc.3} is fulfilled with $\kappa=1$ and in particular this is 
the case when $\{ p,\overline{p}\}\ne 0$ on $p^{-1}(\partial \Gamma )$. 
If $z\in\partial \Sigma \setminus\Sigma _\infty $ then \no{sa.1} cannot 
hold, so if $\partial \Gamma \cap\partial \Sigma \ne 0$, the best we can 
hope for is that
\ekv{gc.3.1} 
{
\forall z \in p^{-1}(\partial \Gamma ) , \hbox{ either }
\{p,\overline{p}\}\ne 0 \hbox{ or } \{p,\{ p,\overline{p}\}\} \ne 0 .
}

This is the situation considered in the 1-dimensional case in \cite{Ha3} 
where deterministic upper bounds on the density of the \ev{}s were 
obtained. Following some arguments there, we shall see that \it if \no{gc.3.1} 
holds, then \no{gc.3} holds with $\kappa={3\over 4}$. \rm In fact, if 
we assume \no{gc.3.1} and if $p(\rho _0)=z_0\in\partial \Gamma $, we 
estimate the contribution to $W_{z_0}(\tau )$ from a \neigh{} of $\rho _0$ 
in the following way:

\par If $\vert \{p,\overline{p}\} (\rho _0)\vert\ge 1/C $, then $d\Re 
p,d\Im p$ are \indep{} near $\rho _0$ and the contribution is ${\cal 
O}(\tau ^2)$. If $\vert \{ p,\overline{p}\}(\rho _0)\vert $ is very small, 
we know that $\vert \{ p,\{ p,\overline{p}\}\}(\rho _0)\vert \ge 1/C$ and 
in order to fix the ideas we assume that $H_{\Re p}^2\Im p(\rho _0)\ge 
1/C$. This means that 
$$
H:=\{ \rho ;\, \{\Re p,\Im p\} (\rho )=0\}
$$
is a smooth hypersurface in a \neigh{} of $\rho _0$ and that $H_{\Re p}$ 
is transversal to $H$ there. A general point in a \neigh{} of $\rho _0$ 
can therefore we written 
$$
\rho =\exp (tH_{\Re p})(\rho '),\ t\in{\rm neigh\,}(0,{\bf R}),\, \rho 
'\in H.
$$
Then $\Re p(\rho )=\Re p(\rho ')$, $\Im p(\rho )=\Im p(\rho 
')+t^2g(t,\rho )$, $g>1/C$. Write $\Re p=s$ so that a general point $\rho '$ 
in $H$ is parametrized by $(s,\rho '')$ with $\rho ''\in {\rm 
neigh\,}(0,{\bf R}^{2n-2})$. 

\par Write $z_0=x_0+iy_0$. For every fixed $\rho ''$, if $\vert 
p(\rho )-z_0\vert \le \tau $ then $\vert s-x_0\vert \le \tau $ and $\vert \Im
p(\rho ')+t^2g(t,s,\rho '')-y_0\vert \le \tau $.  Then we are confined to
an interval of length $2\tau $ in the $s$-variable and, for every such
fixed $s$, to an interval of length ${\cal O}(\tau ^{1/2})$ in the
$t$-variable, or to the union of two such intervals.  By Fubini's theorem,
the contribution to the volume is therefore ${\cal O}(\tau ^{3/2})$.  Hence
$W_{z_0}(\tau )={\cal O}(\tau ^{3/2})$, so $V_{z_0}(t)={\cal O}(t^{3/4})$,
as claimed.
\end{ex}

\par In Section \ref{sa} we introduced the distribution $I(z)$ in \no{sa.2} and showed 
that $I(z)$ is subharmonic, satisfying \no{sa.3}. This implies that 
\ekv{gc.4}
{
\int_{D(z,s)}\Delta (I(w))L(dw)=2\pi  W_z(s),
}
and \no{gc.3} is equivalent to 
\ekv{gc.5}
{
\int_{D(z,t)}\Delta (I(w))L(dw)={\cal O}(t^{\rho _0}),\hbox{ \ufly{} for 
}z\in{\rm neigh\,}(\partial \Gamma ),\ 0\le t\ll 1,
}
with $\rho _0=2\kappa\in ]0,2]$. This is precisely the condition 
\no{cz.37} for $I=\phi $, $\mu =\Delta I$. In view of \no{gc.2}, the 
assumption \no{gc.3} is also equivalent to requiring \no{fc.44} 
to hold \ufly{} 
for $z\in{\rm neigh\,}(\partial \Gamma )$ with 
$q=q_z=|(\widetilde{p}-z)\inv (p-z)|^2$ .

\par Consider the \hol{} \fu{} 
\ekv{gc.6}
{
F_\delta (z;h)=\det P_\delta (z),\ z\in\widetilde{\Omega },
}
where we recall that $P_\delta (z)=(\widetilde{P}-z)\inv (P_\delta -z)$. 
Theorem \ref{lb2} and its proof give:
\begin{prop}\label{gc1}
Let $\delta $ satisfy \no{lb.20}. Then 
there exist constants $C,C_0,\widetilde{C}>0$ such that 
\smallskip
\par\noindent (a) With \proba{} $\ge 1-Ce^{-C_0h^{-2n}}$, we have
\ekv{gc.7}
{
\ln \vert F_\delta (z;h)\vert \le {1\over (2\pi h)^n}(I(z)+Ch^{\delta 
_0}\ln{1\over h}),
}
for all $z$ in some fixed \neigh{} of $\partial \Gamma $.
\smallskip
\par\noindent (b) For every $z\in{\rm neigh\,}(\partial \Gamma )$, $\epsilon \ge 
0$, we have 
\ekv{gc.8}
{
\ln \vert F_\delta (z;h)\vert \ge {1\over (2\pi h)^n}(I(z)-Ch^{\delta 
_0}(\ln {1\over h}+\ln {1\over \delta })-\epsilon ),
}
with \proba{} $\ge 1-Ce^{-\epsilon (2\pi 
h)^{-n}}-\widetilde{C}e^{-C_0h^{-2n}}$.
\end{prop}

\par For the upper bound \no{gc.7}, we recall that the upper bound \no{gp.15} 
was obtained when $\Vert Q\Vert _{{\rm HS}}$ satisfies the estimate \no{gp.4} 
and this event  is \indep{} of $z$.

\par We can now apply Proposition \ref{cz4}, with $\phi $ equal 
to $I+Ch^{\kappa}\ln {1\over h}$ and with $h$ there replaced by 
$(2\pi h)^n$, with $\epsilon $ in 
\no{cz.14} replaced by
$$
{\cal O}(1)(h^{\kappa}\ln{1\over h}+h^{\kappa}\ln{1\over \delta 
}+\epsilon ),
$$
and with $\epsilon $ in \no{gc.8} large enough, so that $\epsilon $ is the dominant term in 
the  last expression. In other words, we take 
\ekv{gc.23}
{
 \epsilon \gg h^{\kappa}\ln 
{1\over \delta },
}
using also that $\ln\delta \inv \ge\ln h\inv$.

\par For $0<r\ll 1$, choose $z_1,...,z_N$ and $N$ as in the first part
of \no{cz.41.5}. Then in view of (b) in the proposition, the last
estimate in \no{cz.41.5} (with $h$ there replaced by $(2\pi h)^n$)
holds for all $j$ with
a probability 
$$\ge 1-{C \over r}e^{-{\epsilon  \over 2}(2\pi h)^{-n}}-\widetilde{C}
e^{-C_0h^{-2n}}.
$$ 
The term 
$${1 \over 2\pi h}\int_\Gamma \Delta \phi L(dz) $$
in \no{cz.43} becomes after the substitutions $h\mapsto (2\pi h)^n$,
$\phi \mapsto I$:
$${1\over (2\pi h)^n}\int_\Gamma  {\Delta I(z)\over 2\pi }L(dz)=
{1\over (2\pi h)^n}{\rm Vol\,}(p\inv (\Gamma )),$$
where we also used \no{sa.3}.
\begin{theo}\label{gc2} Let $\delta $ satisfy \no{lb.20}. 
Assume \no{gc.3}, with $\kappa \in ]0,1]$. Let $N(P+\delta Q_\omega 
,\Gamma )$ be the number of \ev{}s of $P+\delta Q_\omega $ in $\Gamma $. 
Then for every fixed $K>0$ and for $0<r\ll 1$:
\eekv{gc.24}
{
&&\vert N(P+\delta Q_\omega ,\Gamma )-{1\over (2\pi h)^n}\iint_{p\inv 
(\Gamma )}dxd\xi \vert \le} {&&{C\over h^n}({\epsilon \over r}+C_K(r^K +\ln 
({1\over r})\iint_{p\inv (\partial \Gamma +D(0,r))}dxd\xi )),\ 0<r\ll 1,
}
with \proba{} 
\ekv{gc.25}
{
\ge 1-{C\over r}e^{-{\epsilon\over 2}(2\pi 
h)^{-n}}} provided that 
\ekv{gc.26}{
h^{\kappa}\ln{1\over \delta }
\ll \epsilon \ll 1,
}
or equivalently, 
$$
e^{-{\epsilon \over Ch^{\kappa}}}\le \delta ,\ C\gg 1,\ \epsilon \ll 1,
$$
implying that $\epsilon \ge \widetilde{C}h^{\kappa}\ln{1\over h}$, for 
some $\widetilde{C}>0$.
\end{theo}

\par In \no{gc.24} we want the \rhs{}  to be much smaller than $h^{-n}$ so 
it is natural to assume that
\ekv{gc.27}
{
\ln ({1\over r})\iint_{p\inv (\partial \Gamma +D(0,r))}dxd\xi ={\cal 
O}(r^{\alpha _0}),\ r\to 0 ,
}
for some $\alpha _0>0$. When $\kappa\in ]{1\over 2},1]$, we 
automatically have \no{gc.27} with any $\alpha _0\in ]0,2\kappa-1[$. 
 In the \rhs{} of \no{gc.24}, we first choose $N\ge 
\alpha_ 
0$ and we choose $r=\epsilon ^{1/(1+\alpha _0)}$, so that $\epsilon /r$, 
$r^{\alpha _0}$ $=$ ${\cal O}(\epsilon ^{{\alpha_0\over 1+\alpha _0}})$. 
Then the \rhs{} of \no{gc.24} becomes 
$$\le {C\over h^n}\epsilon ^{\alpha_0\over 1+\alpha _0}.$$

\par So, if $1\gg \epsilon \ge \widetilde{C}h^{\kappa}\ln {1\over h}$ 
with $\widetilde{C}$ \sufly{} large, and $\delta $ is as in the theorem, then 
\ekv{gc.30}
{
\vert N(P+\delta Q_\omega ,\Gamma )-{1\over (2\pi h)^n}\iint_{p\inv 
(\Gamma )}dxd\xi \vert \le {C\over h^n}\epsilon ^{\alpha _0\over 1+\alpha _0},
}
with \proba{}
\ekv{gc.31}
{
\ge 1-{C\over \epsilon ^{1\over 1+\alpha _0}}e^{-{\epsilon \over 2}(2\pi 
h)^{-n}}.
}
This expression is very 
close to 1 except possibly in the case $\kappa=1$, $n=1$. In that case, 
we replace $\kappa$ by a strictly smaller value and choose $\delta 
,\epsilon $ as above.

\begin{theo}\label{gc3}
Let ${\cal G}$
be a \fy{} of domains $\Gamma $ as in Theorem \ref{gc2} satifying the 
assumptions there \ufly{} (cf Theorem \ref{sa2}) and in particular we assume 
\no{gc.3} \ufly{} for all $z$ in a \neigh{} of the union of all the 
$\partial \Gamma $. Then we 
have \no{gc.24} 
with \proba{} $$\ge 1-{C\over r^2}e^{-{\epsilon \over (2\pi 
h)^n}}$$ provided that 
$$
h^{\kappa}\ln{1\over \delta }\ll \epsilon \ll 1.
$$

\end{theo}

\section{Appendix: Gaussian random variables in Hilbert spaces} \label{ap}
\setcounter{equation}{0}

\par In this appendix we review some generalities about Gaussian random variables 
in Hilbert spaces that seem to be quite standard to probabilists. 

\par Let $\alpha _1,\alpha _2,...$ be a sequence of \indep{} ${\cal 
N}(0,1)$-laws, and let ${\cal H}$ be a complex separable Hilbert space.
\begin{prop}\label{ap1} Let
$v_1,v_2,...\in{\cal H}$ be a 
sequence of vectors such that $\sum_1^\infty \Vert v_j\Vert ^2<\infty $, 
then if the sequence $n_1<n_2<...$ tends to $\infty $ \sufly{}
fast, we have that 
$$
\lim_{k\to \infty }\sum_1^{n_k}\alpha _j(\omega )v_j \hbox{ exists almost 
surely (a.s.)}.
$$
Let $S(\omega )$
denote the almost sure limit. If $\widetilde{n}_k$ is another increasing 
sequence tending to infinity, such that the limit 
$$\lim_{k\to \infty }\sum_1^{\widetilde{n}_k}\alpha _j(\omega 
)v_j=:\widetilde{S}(\omega )$$ 
exists almost surely, then $\widetilde{S}(\omega )=S(\omega )$
a.s.\end{prop}
\begin{proof} Let $(\Omega ,\mathsf P)$ be the underlying \proba{} space. Then 
$f_j:= \alpha _j(\omega )v_j$ can be viewed as elements of $L^2(\Omega 
,{\cal H})$ of norm $\Vert v_j\Vert $. They are mutually \og{} since the 
$\alpha _j$ are \indep{}. We thus have an \og{} sum $\sum_1^\infty \alpha 
_jv_j$ which converges in $L^2(\Omega ;{\cal H})$ and as usual, using the 
Chebyschev inequality, we deduce the existence of a sequence of partial sums 
that converges a.s. 
\end{proof}

\par Let 
$e_1,e_2,..$ and $f_1,f_2,...$ be two  \on{} bases in ${\cal H}$. 
Let $\alpha _1(\omega ),\alpha _2(\omega ),...$ be \indep{} complex 
${\cal N}(0,1)$-laws, and consider the formal vector $\sum_1^\infty  \alpha 
_j(\omega )e_j$. Almost surely, $\{ \alpha _j(\omega )\} _1^\infty $ is not 
in $\ell ^2$ so our vector is not in ${\cal H}$. However, if $v\in{\cal 
H}$, then a.s., we can define the scalar product 
\ekv{ap.2}
{
(\sum_1^\infty \alpha _j(\omega )e_j\vert v)=\sum_1^\infty \alpha 
_j(\omega )(e_j\vert v)
}
as in Proposition \ref{ap1}, since $\{ (e_j\vert v)\} _1^\infty \in \ell ^2$.

\par We now look for random variables $\beta _1(\omega ),\beta 
_2(\omega ),...$ such that 
\ekv{ap.3}
{
\sum_1^\infty \alpha _k(\omega )e_k=\sum_1^\infty \beta _j(\omega )f_j,
}
in the sense that the formal scalar products with $f_1,f_2,...$ are equal. 
This leads to the definition
\ekv{ap.4}
{
\beta _j(\omega )=\sum_{k=1}^\infty (e_k\vert f_j)\alpha _k(\omega ),
}
which is well-defined as in Proposition \ref{ap1}, since $k\mapsto 
(e_k\vert f_j)$ is in $\ell^2$. For every finite $N$, the variable 
\ekv{ap.5}
{
\sum_{k=1}^N (e_k\vert f_j)\alpha _k(\omega )
}
has the density 
$$
*_{k=1}^N {1\over \pi \vert (e_k\vert f_j)\vert ^2}e^{-\vert \alpha \vert 
^2/\vert (e_k\vert f_j)\vert ^2},
$$
where $*$ indicates convolution products. Hence the characteristic \fu{} 
(i.e. the \F{} \tf{}) is 
$$\exp \Bigl(-{1\over 4}(\sum_{k=1}^N \vert (e_k\vert 
f_j)\vert ^2)\vert \xi \vert ^2\Bigr),$$ 
so \no{ap.5} is a normal distribution 
${\cal N}(0,{\sum_{k=1}^N\vert (e_k\vert f_j)\vert ^2})$. The unitarity of the 
matrix $((e_k\vert f_j))$ then implies that $\beta _j(\omega )$
is a ${\cal N}(0,1)$-law.

\begin{prop}\label{ap2}
$\beta _j$ are \indep{} ${\cal N}(0,1)$-laws.
\end{prop}
\begin{proof}
We have already seen that $\beta _j$ are ${\cal N}(0,1)$-laws. To see that they are
\indep{}, we compute (using Proposition \ref{ap1}) the joint distribution 
of $\beta _1,\beta _2,...,\beta _N$. Write
$$
\beta ^{(N)}=\begin{pmatrix}\beta _1\cr \beta _2\cr ..\cr \beta 
_N\end{pmatrix}=\sum_1^\infty  \alpha _k\nu _k,
$$
where 
$$
\nu _k=\begin{pmatrix}\sigma _{1,k}\cr\sigma _{2,k}\cr \dots\cr \sigma _{N,k} \end{pmatrix},\ 
\sigma _{j,k}=(e_k\vert f_j).
$$
$\alpha _k(\omega )\nu _k$ is a \rv{} with values 
in ${\bf C}^N$ and with the characteristic \fu{}
\begin{eqnarray*}
\chi _{\alpha _k\nu _k}(\xi )&=&\int e^{-i\Re \alpha  (\nu _k\vert \xi 
)}e^{-\vert \alpha \vert ^2}{L(d\alpha )\over \pi }\\
&=& \exp \left( {-{1\over 4}\vert (\nu _k\vert \xi )\vert ^2} \right)\\
&=& \exp \left( -{1\over 4}\sum_{\ell =1}^N\sum_{m=1}^N \sigma _{\ell 
,k}\overline{\sigma _{m,k}}\overline{\xi }_\ell \xi _m \right).
\end{eqnarray*}
It follows that 
$$
\chi _{\sum_1^\infty \alpha _k\nu _k}(\xi )=\exp \left(-{1\over 4}
\sum_{\ell =1}^N\sum_{m=1}^N (\sigma _\ell \vert 
\sigma _m)\overline{\xi }_\ell \xi _m \right),
$$
where $\sigma _j=(\sigma _{j,k})_{k=1}^\infty \in\ell^2$. But the $\sigma _j$ 
form an \on{} system, so finally,
$$
\chi _{\sum_1^\infty \alpha _k\nu _k}(\xi )=\exp -{1\over 4}\vert \xi 
\vert ^2.
$$
This means that the joint distribution of $\beta _1,...,\beta _N$ is 
$$
{1\over \pi ^N}e^{-\vert \beta \vert ^2}L_{{\bf C}^N}(d\beta ),
$$
and that $\beta _1,...,\beta _N$ are \indep{}.
\end{proof}
\par The \rv{} \no{ap.2} is an ${\cal N}(0,\Vert v\Vert^2 )$-law. 
\par If $v\in{\cal H}$ is any finite linear combination of the $f_j$, we 
know by construction that 
$$
(\sum_1^\infty \alpha _j(\omega )e_j\vert v)=(\sum_1^\infty  \beta 
_j(\omega )f_j\vert v),\ {\rm a.s.}
$$
If $v\in{\cal H}$ is \aby{}, we write $v=v_\epsilon +r_\epsilon $, where 
$v_\epsilon $ is a finite linear combination of the $f_j$ and $\Vert 
r_\epsilon \Vert <\epsilon $. We conclude that almost surely,
$$
(\sum_1^\infty \alpha _j(\omega )e_j\vert v)=(\sum_1^\infty \beta 
_j(\omega )f_j\vert v)+(\sum_1^\infty \alpha _j(\omega )e_j\vert 
r_\epsilon )-(\sum_1^\infty \beta _j(\omega )f_j\vert r_\epsilon ).
$$
Here the last two terms are ${\cal N}(0,\Vert r_\epsilon \Vert ^2)$-laws and hence 
as small as we like with a \proba{} as close as we like to 1, 
when $\epsilon $
is small eneough. We conclude that 
\ekv{ap.6}
{
(\sum_1^\infty \alpha _j(\omega )e_j\vert v)=(\sum_1^\infty \beta 
_j(\omega )f_j\vert v)\ {\rm a.s.}
}
\begin{prop}\label{ap3}
Let ${\cal H}$, $\widetilde{{\cal H}}$ be two separable Hilbert spaces and 
let $T:{\cal H}\to \widetilde{{\cal H}}$ be a \hs{} \op{}. Let $\alpha 
_j(\omega )e_j$, $\beta _j(\omega )f_j$ be as above. Then 
$T(\sum_1^\infty  \alpha _j(\omega )e_j)$
is well defined a.s. and equal to $T(\sum_1^\infty \beta _j(\omega )f_j)$ 
a.s.
\end{prop}
\begin{proof}
We define $T(\sum_1^\infty  \alpha _j(\omega )e_j)$ as $\sum_1^\infty 
\alpha _j(\omega )Te_j$ in the sense of Proposition \ref{ap1}, using that 
$$
\sum \Vert Te_j\Vert _{\widetilde{{\cal H}}}^2=\Vert T\Vert _{{\rm HS}}^2<\infty .
$$
Notice also that for every $v\in\widetilde{{\cal H}}$ we have a.s.
\begin{eqnarray*}
(T(\sum_1^\infty \alpha _j(\omega )e_j)\vert v)&=&\sum_1^\infty  \alpha 
_j(\omega )(Te_j\vert v)\ {\rm a.s.}\\
&=& \sum_1^\infty  \alpha _j(\omega )(e_j\vert T^*v). 
\end{eqnarray*}
The same considerations apply to $T(\sum_1^\infty  \beta _j(\omega )f_j)$ 
so in view of \no{ap.6}, for every $v\in\widetilde{{\cal H}}$ we have 
$$
(T(\sum_1^\infty  \alpha _j(\omega )e_j)\vert v)=(T(\sum_1^\infty \beta 
_j(\omega )f_j)\vert v)\ {\rm a.s.}
$$
We get the same conclusion a.s. simultaneously for all $v$ in any countable
set, and letting $v=g_1,g_2,....$, where $g_j$
form and \on{} basis in $\widetilde{{\cal H}}$, we conclude that 
$$
T(\sum_1^\infty  \alpha _j(\omega )e_j)=T(\sum_1^\infty \beta _j(\omega 
)f_j)\ {\rm a.s.}
$$
\end{proof}

\par Now let ${\cal E},{\cal F}, {\cal G}, {\cal H}$ be separable Hilbert 
spaces and let $T:{\cal E}\to {\cal F}$, $S:{\cal G}\to {\cal H}$  be 
\hs{} \op{}s. If $f\in{\cal F}$, $g\in{\cal G}$, we also denote by $g,f$ 
the corresponding multiplication \op{}s ${\bf C}\ni z\mapsto zg,zf\in 
{\cal G}, {\cal F}$, so that $f^*u=(u\vert f)$. Then $gf^*u=(u\vert f)g$ 
defines an \op{} $:{\cal F}\to G$ which has the \hs{} norm $\Vert g\Vert 
\Vert f\Vert $. Let $f_j$, $j=1,2,...$, $g_j$, $j=1,2,...$
be \on{} bases in ${\cal F}$, ${\cal G}$. Then 
$\{g_jf_k^*\}_{j,k=1}^\infty $ is an \on{} basis for the space ${\rm HS}({\cal 
F},{\cal G}) $of \hs{} 
\op{}s ${\cal F}\to {\cal G}$. Now,
$$Sg_jf_k^*T=(Sg_j)(T^*f_k)^*,$$ 
and 
$$
\Vert Sg_jf_k^*T\Vert _{{\rm HS}}^2=\Vert Sg_j\Vert ^2\Vert T^*f_k\Vert ^2.
$$
It follows that 
$$
\sum_{j,k}\Vert Sg_jf_k^*T\Vert _{{\rm HS}}^2=\Vert S\Vert _{{\rm HS}}^2\Vert T\Vert 
_{{\rm HS}}^2,
$$
and we conclude that
$$
{\rm HS}({\cal F},{\cal G})\ni A\mapsto SAT\in {\rm HS}({\cal E},{\cal H})
$$
is a \hs{} \op{}. The earlier discussion can therefore be applied:
\begin{prop}\label{ap4}
Let 
$\alpha _{j,k}(\omega )$ be \indep{} ${\cal N}(0,1)$ laws. Then 
$$
S\sum_{j,k}\alpha _{j,k}(\omega )g_jf_k^*T=\sum_{j,k}\alpha _{j,k}(\omega 
)Sg_jf_k^*T
$$ 
is almost surely defined as a \hs{} \op{}. Moreover, if 
$\widetilde{g}_j$, $\widetilde{f}_k$ are new \on{} bases in ${\cal G}$, 
${\cal F}$, then there exists a new set of \indep{} ${\cal N}(0,1)$-laws $\beta 
_{j,k}(\omega )$ such that 
\ekv{ap.7}
{
S\circ (\sum_{j,k}\alpha _{j,k}(\omega )g_jf_k^*)\circ T=
S\circ (\sum_{j,k}\beta _{j,k}(\omega )\widetilde{g}_j\widetilde{f}_k^*)\circ 
T\ {\rm a.s.}
}
\end{prop}

\end{document}